\documentclass[reqno,11pt]{amsart}
\usepackage{amsmath,amssymb,amsthm,amsfonts,mathrsfs,physics}

\usepackage{dsfont} 
\usepackage[utf8]{inputenc}
\usepackage[T1]{fontenc}
\usepackage[english]{babel}
\usepackage{comment}

\usepackage{fullpage}
\usepackage{enumitem}
\usepackage[toc]{appendix}
\usepackage{esint}
\usepackage{url}
\usepackage{graphics,graphicx}
\usepackage{ifthen}

\usepackage[fixamsbugs]{mathtools}

\usepackage{palatino}

\allowdisplaybreaks

\usepackage{caption}  
\usepackage{tikz} 

\newcommand{\one}{\mathbf{1}}

\usepackage[backref=page]{hyperref}
 \hypersetup{
 	plainpages=false,
 	unicode=false,          
 	pdftoolbar=true,        
 	pdfmenubar=true,        
 	pdffitwindow=true,      
 	pdftitle={My title},    
 	pdfauthor={},     
 	pdfsubject={},   
 	pdfcreator={},   
 	pdfproducer={}, 
 	pdfkeywords={}, 
 	pdfnewwindow=true,      
 	colorlinks=true,       
 	linkcolor=blue,          
 	citecolor=blue,        
 	filecolor=blue,      
 	urlcolor=blue,          
 	urlbordercolor=0 1 1
 }

\usepackage{xcolor}
\definecolor{Blue}{rgb}{0.,0.,1}

\newcommand{\tV}{\widetilde{V}}

\usepackage{setspace}

\usepackage[pass]{geometry}

\newcommand{\xC}{{\rm C}}
\newcommand{\del}{\Delta t}

\newcommand{\R}{\mathbb R}
\newcommand{\N}{\mathbb N}

\newcommand{\cA}{\mathcal{A}}
\newcommand{\cH}{\mathcal{H}}

\newcommand{\cC}{\mathcal{C}}

\newcommand{\cL}{\mathcal{L}}

\newcommand{\cM}{\mathcal{M}}

\newcommand{\cT}{\mathcal{T}}

\newcommand{\G}{{\mathbb G}_{d,n}}
\newcommand{\V}{\|V\|}
\newcommand{\W}{\|W\|}

\newtheorem{theo}{Theorem}[section]
\newtheorem*{theo*}{Theorem}
\newtheorem{prop}[theo]{Proposition}
\newtheorem{lemma}[theo]{Lemma}
\newtheorem{dfn}[theo]{Definition}
\newtheorem*{dfn*}{Definition}
\newtheorem{cor}[theo]{Corollary}

\newtheoremstyle{rmdotless}{}{}{\upshape}{}{\bfseries}{.}{0.5em}{}
\theoremstyle{rmdotless}
\newtheorem{remk}[theo]{Remark}

\DeclareMathOperator{\mdiv}{div}

\DeclareMathOperator*{\supp}{spt}

\DeclareMathOperator*{\lip}{Lip}

\renewcommand{\phi}{\varphi}
\renewcommand{\epsilon}{\varepsilon}
\newcommand{\e}{\epsilon}
\renewcommand{\G}{G_{d,n}}

\newcommand{\res}{\mathop{\hbox{\vrule height 7pt width .5pt depth 0pt
\vrule height .5pt width 6pt depth 0pt}}\nolimits} 

\newcommand{\sleq}{<}

\newcommand{\dt}{dt}

\usepackage{xcolor}

\def\ds{\displaystyle}

\thanks{B. Buet acknowledges support from the French National Research Agency (ANR) under grant ANR-21-CE40-0013-01 (project GeMfaceT). G.P.Leonardi acknowledges support from the Italian Agency for University and Research (Ministero dell'Università e della Ricerca), specifically under the following grants: PRIN 2017TEXA3H ``Gradient flows, Optimal Transport and Metric Measure Structures''; PRIN 2022PJ9EFL ``Geometric Measure Theory: Structure of Singular Measures, Regularity Theory and Applications in the Calculus of Variations'' financed by European Union - Next Generation EU, Mission 4, Component 2 - CUP:E53D23005860006. S. Masnou and A. Sagueni acknowledge support from the French National Research Agency (ANR) under grants ANR-19-CE01-0009-01 (project MIMESIS-3D) and ANR-24-CE40-XXX (project STOIQUES). Part of this work was also supported by the LABEX MILYON (ANR-10-LABX-0070) of Universit\'e de Lyon, within the program ”Investissements d’Avenir” (ANR-11-IDEX- 0007) operated by the French National Research Agency (ANR), and by the European Union Horizon 2020 research and innovation programme under the Marie Sklodowska-Curie grant agreement No 777826 (NoMADS).}

\author{Blanche Buet}
\address{Universit\'e Paris-Saclay, Inria, CNRS, LMO UMR8628, 91400 Orsay, France}
\email{blanche.buet@universite-paris-saclay.fr}

\author{Gian Paolo Leonardi}
\address{Dipartimento di Matematica, via Sommarive 14, IT-38123 Povo - Trento, Italy}
\email{gianpaolo.leonardi@unitn.it}

\author{Simon Masnou}
\address{Université Claude Bernard Lyon 1, CNRS, Ecole Centrale de Lyon, INSA Lyon, Université Jean Monnet, ICJ UMR5208, 69622 Villeurbanne, France}
\email{masnou@math.univ-lyon1.fr}

\author{Abdelmouksit Sagueni}
\address{Université Claude Bernard Lyon 1, CNRS, Ecole Centrale de Lyon, INSA Lyon, Université Jean Monnet, ICJ UMR5208, 69622 Villeurbanne, France}
\email{sagueni@math.univ-lyon1.fr}

\date{\today}
\title{Approximate mean curvature flows of a general varifold,\\ and their limit spacetime Brakke flow}

\keywords{Geometric measure theory; varifolds; approximate mean curvature flow; Brakke flow.}

\begin{document}

\begin{abstract}
We propose a construction of mean curvature flows by approximation for very general initial data, in the spirit of the works of Brakke and of Kim \& Tonegawa based on the theory of varifolds. Given a general varifold, we construct by iterated push-forwards an approximate time-discrete mean curvature flow depending on both a given time step and an approximation parameter. We show that, as the time step tends to $0$, this time-discrete flow converges to a unique limit flow, which we call the approximate mean curvature flow. An interesting feature of our approach is its generality, as it provides an approximate notion of mean curvature flow for very general structures of any dimension and codimension, ranging from continuous surfaces to discrete point clouds. We prove that our approximate mean curvature flow satisfies several properties: stability, uniqueness, Brakke-type equality, mass decay. By coupling this approximate flow with the canonical time measure, we prove convergence, as the approximation parameter tends to $0$, to a spacetime limit measure whose generalized mean curvature is bounded. Under an additional rectifiability assumption, we further prove that this limit measure is a spacetime Brakke flow.
\end{abstract}

\maketitle
\tableofcontents
\section{Introduction}
Let $\cM$ be a $d$-dimensional manifold and $F_0: \cM\to\R^n$ a smooth embedding. The mean curvature flow starting from $\cM_0=F_0(\cM)$ is, see e.g.~\cite{mantegazzabook,belbook}, a smooth time-dependent family of embeddings $F: \cM\times[0,T)\to\R^n$ satisfying
\begin{equation*}
\left\{\begin{array}{lll}
 \ds\frac{\partial F}{\partial t}(x,t) &=& H(x,\cM_t)\\
 \ds F(x,0)&=&F_0(x)\end{array}\right.
\end{equation*}
where $H(x, \cM_t)$ is the mean curvature vector of $\cM_t=  F(\cM,t)$ at $x_t = F(x,t)$ defined by
\begin{equation}\label{intro_mcdfn}
 H(x, \cM_t)= -\sum\limits_{j=1}^{n-d} (\mdiv_{\cM_t}\nu_j) \nu_j,
\end{equation}
with $\lbrace \nu_j\rbrace_j$ any orthonormal basis of the normal space $T_{x_t} \cM_t^{\perp}$ at $x_t$ and $\mdiv_{\cM_t}$ the tangential divergence on $\cM_t$. The time evolution equation above is equivalent to  \(\frac{\partial F}{\partial t}(x,t) = \Delta_{\cM_t} F(x,t)\),
where $\Delta_{\cM_t}$ is the Laplace-Beltrami operator on $\cM_t$ associated to the metric induced by $F(\cdot,t)$. By the theory of pseudo-parabolic PDEs, this equation has a smooth solution defined on a nontrivial time interval~\cite{eidel}.

Grayson proved that the mean curvature flow shrinks closed curves in $\R^2$ into single points in finite time~\cite{gra87}. The result extends to mean convex hypersurfaces of $\R^n$, by a result due to Huisken~\cite{hui}. But there are well known examples of smooth hypersurfaces that may develop singularities, other than a limit point, in finite time, see for instance the construction by Grayson in~\cite{gra87} of a dumbbell flowing by mean curvature, whose neck pinches off before the two bells shrink.

To extend the definition of the mean curvature flow beyond singularities, several approaches have been proposed that yield weak notions of mean curvature flow:
\begin{itemize}
 \item The Brakke flow~\cite{brakke,kt,ton,Ilmanen}, defined in the setting of rectifiable varifolds and, at least for the original construction of Brakke, in arbitrary codimension;
 \item The level set formulation, based on an implicit representation of the evolving interface and the PDE satisfied by this representation~\cite{evsp1,evsp2,evsp3,evsp4,cgg,AmSo};
 \item Phase fields methods, based on diffuse representations of evolving interfaces and associated reaction-diffusion PDEs~\cite{allen,Ilmanen2,henlau};
 \item Diffusion generated motions, based on convolution-thresholding schemes~\cite{MBO,LauxOtto};
 \item De Giorgi's method of minimal barriers~\cite{DeG1,AmSo,BelPao,BeNo99};
 \item The minimizing movement approach, based on time-discrete variational approximations \cite{lust,ATW}.
 \end{itemize}
Further references and details can be found in~\cite{belbook,smoczyk,Ecker, amb-mat98, Ambrosio2000}.

\medskip
The starting point of our work was the need for a weak notion of mean curvature flow that is adapted to structured data (such as interfaces) and unstructured data (such as point clouds) in arbitrary dimension and codimension (including codimension $0$ when dealing with volumetric data).  Among the above mentioned approaches, only Brakke flows, level set flows, diffusion generated motions, and flows based on De Giorgi's barriers can be considered in higher codimension. None of them, though, has been explicitly adapted to handle in a consistent way both unstructured data such as point clouds and volumetric data. We propose in this paper a construction {\it à la} Brakke that achieves this goal.

\medskip The core notion in Brakke's~\cite{brakke} and Kim \& Tonegawa's~\cite{kt} constructions is the notion of varifold. 
Recall that a $d$--varifold $V$ in $\R^n$ is a non-negative Radon measure on $\R^n \times \G$, where $\G$ is the Grassmannian manifold of $d$-dimensional vector subspaces of $\R^n$~\cite{simon}. $V_d(\R^n)$  denotes the space of $d$--varifolds in $\R^n$. A varifold $V$ is associated with a mass measure $\|V\|$ characterized by its action on every Borel set $A\subset \R^n$: $\|V\|(A)=V(A\times \G)$. $V$ is called rectifiable if it decomposes as $V=\theta \cH^d\res\cM\otimes \delta_{T_{x}\cM}$ where $\cM \subset\R^n$ is $d$--rectifiable, $\theta \in L_{loc}^1(\cH^d\res\cM)$ is nonnegative, and $T_x\cM$ is the approximate tangent space of $\cM$ at $x$. If, in addition, $\theta \in \N $ a.e., then $V$ is called integral. 

The  \textit{first variation} of a varifold $V\in V_d(\R^n)$ of finite mass is the map 
\[
\delta V:\; X \in C^1(\R^n, \R^n)\longmapsto \int_{\R^n \times \G} \mdiv_S(X)(x)\, dV(x,S).
\]
If $\delta V$ is locally bounded, it can be represented by a vector Radon measure thanks to the Riesz representation theorem. Then the generalized mean curvature $H(\cdot, V)$ of $V$ is (minus) the Radon-Nikodym derivative of $\delta V$ with respect to $\|V\|$.

Let $(\cM_t)_{t\in[0,T)}$ be a smooth mean curvature flow as at the beginning of this paper. Denote $(M_t)_{t\in[0,T)}$ the associated family of integral varifolds, i.e. $M_t = \cH^d\res\cM_t\otimes \delta_{T_x\cM_t}$, $\forall t\in [0,T)$. The mean curvature flow is fully characterized by the following Brakke integral equality~\cite{ton}: for all test function $\phi \in \xC_c([0,T) \times \R^n,\R_+)$, and for any pair $(t_1,t_2)$ such that $0\leq t_1<t_2<T$,
\begin{multline}\label{intro_eq:weak_mcf}
 \|M_{t}\|( \phi(t_2,\cdot)-\|M_{t}\|( \phi(t_1,\cdot) = \\\int_{t_1}^{t_2} \int_{\cM_t} \Big[ - \phi(t,x) |H(x,\cM_t)|^2 +(T_x\cM_t)^{\perp}\bigl( \nabla\phi(t,x) \bigr) \cdot H(x,\cM_t)  +\partial_t \phi(t,x)\Big]\,d\cH^d(x) dt.
 \end{multline}
This leads to the notion of Brakke flow in the integrated sense of Kim \& Tonegawa~\cite{kt}:
\begin{dfn}\label{intro_def_brakke-V2} A one-parameter family $(V_t)_{t\in[0,T)}$ in $V_d(\R^n)$ is called a Brakke flow if
\begin{enumerate}
 \item For a.e. $t\in[0,T)$, $V_t$ is integral.
 \item For any compact $K\subset\R^n$ and $t< T$, $\sup_{s\in [0,t)}\| V_s \| (K) < \infty$. 
 \item For a.e. $t\in[0,T)$, $V_t$ has locally bounded first variation and $\delta V_t  <\!< \| V_t \| $.
 \item $H(\cdot,V_t):=-\delta V_t/ \| V_t \|  \in L_{loc}^2(\dt \otimes \|V_t\|)$.
 \item $(V_t)_{t\in[0,T)}$ satisfies the Brakke inequality, i.e., for $0 \leq t_1 \leq t_2 < T$ and $\phi \in \xC_c^1([0,T) \times \R^n , \R_+)$, 
 \begin{multline}\label{intro_eq:brakke_inequality-intro}
  \| V_{t_2}\|(\phi(t_2,\cdot))- \| V_{t_1}\|(\phi(t_1,\cdot)) \leq \\
  \int_{t_1}^{t_2} \int_{\R^n} \Big[-\phi(t,x) |H(x,V_t)|^2 + \nabla \phi(t,x) \cdot H(x,V_t) + \partial_t \phi(t,x)\Big] \, d \|V_t\|(x) dt. 
  \end{multline}
\end{enumerate}
\end{dfn}
The reason an inequality is required in~\eqref{intro_eq:brakke_inequality-intro}, rather than an equality as in the smooth case~\eqref{intro_eq:weak_mcf}, is to allow for sudden mass loss or topological changes, which is necessary to obtain a consistent definition of weak mean curvature flow. 

Starting from an initial $(n-1)$-varifold associated with an open partition of $\R^n$, Kim \& Tonegawa's construction of a Brakke flow is obtained as a continuous limit of an iterated two-step scheme. The first step involves a desingularization map which is essential to go beyond singularities (but different choices for the map may yield different flows, so the uniqueness cannot be guaranteed). The second step uses a push-forward map involving an approximate mean curvature depending on a scale $\e$. Brakke's own construction is more general, for it can move by mean curvature integral varifolds of arbitrary codimension. However, it does not exclude the triviality of the flow, i.e. $V_t = 0, \, \forall t>0$, even starting from a smooth set. In contrast, Kim \& Tonegawa's construction guarantees the nontriviality of the flow starting from a smooth $(n-1)$-set.

Neither Brakke's approach nor Kim \& Tonegawa's approach can handle point clouds. For the latter, it is because isolated points do not define a proper open partition of $\R^n$. As for Brakke's construction, the flow starting from a point cloud is trivial.
More precisely, in both approaches, for given parameters $\Delta t$ (time step) and $\epsilon$ (regularization scale), alternating a time-discrete Euler scheme applied to an $\epsilon$--regularization of the curvature with a desingularization step defines a time-discrete approximate flow. As a second step, a diagonal extraction procedure yields a flow in the joint limit $(\Delta t, \epsilon) \to 0$, but this flow is unfortunately trivial when starting from a point cloud varifold. In this paper, we propose a construction that first lets $\Delta t \to 0$ alone, yielding a time-continuous $\epsilon$-approximate flow, see Theorem~\ref{damcfconvergence}. This limiting flow is well-defined starting from any varifold with compact support, and nontrivial in the setting of point clouds. We also study the limit $\epsilon \to 0$, and naturally the limit flow is trivial when starting from a point cloud. However, our construction provides some control: if one considers as initial data a sequence of point cloud varifolds $W_k$ converging to a submanifold $\cM$, our construction allows us to quantify, at a given regularization scale $\epsilon$, the discrepancy between the approximate flow starting from $\cM$ and the approximate flows starting from the $W_k$'s. In particular, there exists a sequence of scales $\epsilon_k \to 0$ ensuring that this discrepancy vanishes as $k \to \infty$, see Remark~\ref{remkPointCloudApprox}.

A key advantage of the construction we propose in this paper is that approximate flows can be defined starting from {\it any} varifold with compact support and arbitrary codimension.  However, the absence of a counterpart to Brakke's or Kim \& Tonegawa's desingularization steps implies that stationary varifolds, whose regularized mean curvature also vanishes at any scale $\epsilon > 0$, do not evolve under the flow. For instance, two non-parallel lines will not evolve under the approximate flow at any scale $\epsilon$.

%
%
\par
The paper is organized as follows: 
\begin{enumerate}[nosep,topsep=-\parskip]
\item Preliminary section~\ref{sec:preliminaries} contains the notions necessary to our construction.
 \item In Section~\ref{def_time_discrete}, we recall the approximate notion of mean curvature for varifolds associated with an approximation scale $\e>0$ defined in~\cite{kt} after~\cite{brakke}. Given a varifold $V$ with finite mass, we build a time-discrete approximate mean curvature flow starting from $V$. The construction relies on iterated push-forwards of V by diffeomorphisms of the form $\operatorname{id}+ \Delta t\, h_\e$, where $\Delta t$ is a given time step and $h_\e$
  is the approximate mean curvature. We exhibit stability properties of this time-discrete flow with respect to time subdivisions and with respect to the initial varifold.
  \item We let the time step $\Delta t$ go to $0$ in Section~\ref{def_time_step_to_0}, and we show that the time-discrete flow constructed in  (2)  has a consistent limit. We consider this limit as an approximate mean curvature flow for the approximation scale $\e$, and we prove that it fulfills several properties (uniqueness, stability, Brakke-type equality).
 \item In Section~\ref{sec:spacetime-Brakke-flows}, we consider the measure defined in (3) coupled with the time measure $dt$, and we exhibit a limit as the scale of approximation $\e$ tends to $0$. Under a rectifiability assumption on the limit measure, we prove that it satisfies a spacetime Brakke inequality. This limit measure can be interpreted as a spacetime track of a generalized Brakke flow. 
 \item The appendix collects a few useful lemmas.
\end{enumerate}

\section{Preliminaries}\label{sec:preliminaries}

\subsection{Notations}
Throughout the paper, we let $d,n \in \N$ be such that $1\leq d \leq n$, $2\leq n$, and we adopt the following notations:
\begin{itemize}
\item $\cL^n$, $\cH^d$ denote the $n$--dimensional Lebesgue measure and the $d$--dimensional Hausdorff measure in $\R^n$, respectively.
\item $B_r(x)$ denotes the open ball of radius $r>0$ and center $x \in \R^n$. We set $B_r=B_r(0)$. For closed balls, $B$ is replaced by $\overline{B}$.
\item For $k\in\N$, $\omega_k$ denotes the volume of the $k-$dimensional unit ball.
\item For a set $A$ and $\delta>0$, $A^{\delta} :=\bigcup_{x\in A}B_{\delta}(x)=\lbrace y \in \R^n, d(y,A)\sleq \delta \rbrace$.
\item $\mathcal{M}_{p,q}$ is  the space of real matrices with $p$ rows, $q$ columns
\item The default matrix norm $\| \cdot \|$ considered in $\mathcal{M}_{p,q}$ is the operator $2$-norm associated with the Euclidean norms $| \cdot |$ in $\R^p$ and $\R^q$. We also consider the norm $| \cdot |_\infty$ defined as $| M |_\infty = \max_{\substack{i = 1 \ldots p \\ j = 1 \ldots q}} |M_{ij}|$, for $M \in \cM_{p,q}$ and we recall the classical relation:
\begin{equation} \label{eq:maxEuclNorms}
\forall M \in \cM_{p,q},\qquad | M |_\infty \leq \| M \| \leq \sqrt{p q}\, | M |_\infty \: .
\end{equation}
\item The space $\mathcal{M}_{p}=\mathcal{M}_{p,p}$ will be for some calculations equipped with the scalar product
\(M:N = \tr(MN^t)=\tr(NM^t)\), and the associated Frobenius norm.
\item $\G$ is the Grassmannian manifold of $d$-dimensional vector subspaces of $\R^n$. We identify $S\in\G$ with its orthogonal projection on the $d$--subspace $S \in \cM_n(\R)$. The distance between $S,T \in \G$ is $\| S - T \|$, where $\| \cdot \|$ is the matrix norm introduced above.
\item The functions and vectors involved may depend on both space and time, we conventionally use $\nabla$ for space derivation and $\partial_t$ for time derivation.
\item Given two topological spaces $X$ and $Y$, $\xC_c (X,Y)$ denotes the space of continuous and compactly supported functions $f : X \rightarrow Y$.
\item Given an open subset $U$ of $\R^n$, $\xC^{k}(U, \R^m)$ denotes the space of functions $u : U \rightarrow \R^m$ that are of class $\xC^k$ and, for $u \in \xC^{k}(U, \R^m)$, $\|u\|_{\xC^k} = \sum_{i=0}^k \| D^i u \|_\infty $.
\item For $a,b\in\R$ with $a < b$ and $m\geq1$, $\cT=\lbrace t_i \rbrace_{i=0}^m$ is called a subdivision of $[a,b]$ if $a=t_0 < t_1 \dots < t_m=b$. We denote $\delta(\cT) := \max_{i \in \lbrace 1, \dots, m \rbrace}t_i - t_{i-1}$.%
\end{itemize}

\subsection{Varifolds}\label{intro_sec:varifolds}

We recall below basic definitions and results concerning varifolds, see~\cite{simon}.
\begin{dfn}(Varifolds)\label{varifolds}${ }_{ }$\\
A $d$--varifold in $\R^n$ is a nonnegative Radon measure on $\R^n \times \G$. The space of $d$--varifolds in $\R^n$ is denoted as $V_d(\R^n)$.
Every varifold $V$ is associated with its mass measure $\|V\|$, a Radon measure on $\R^n$ defined for every Borel set $A\subset \R^n$ by $\|V\|(A)=V(A\times \G)$.
\end{dfn}

\noindent We collect below classical examples of varifolds, we refer to~\cite{brakke,menne-survey} for more general, less trivial, examples.
\begin{enumerate}
\item {\em Rectifiable varifolds}: $V\in V_d(\R^n)$ is rectifiable if it decomposes as $V = \theta \cH^d\res\cM\otimes \delta_{T_x\cM}$ where $\cM \subset\R^n$ is a $d$--rectifiable set~\cite{simon,afp}, $T_x\cM$ is the approximate tangent space of $\cM$ at $x\in\cM$, and $\theta \in L_{loc}^1(\cH^d\res\cM)$ is a nonnegative function called the multiplicity of $V$.  We denote $V=\mathbf{v}(\cM,\theta)$. The associated mass measure is $\V = \theta \cH^d\res\cM$.
 \item {\em Integral varifolds}: a rectifiable $d$-varifold $V={\mathbf v}(\cM,\theta)$ is called {\it integral} if $\theta(x) \in \N $ for $\cH^d$--almost every $x \in \cM$.
 \item {\em Point cloud varifolds}: to a finite collection of points $\lbrace x_j \rbrace_{j=1}^{N}$ in $\R^n$, $d$-planes $\lbrace P_j \rbrace_{j=1}^{N}$ in $\G$ and masses $\lbrace m_j \rbrace_{j=1}^{N}$ in $\R_+$, we may associate the varifold $V=\sum_{j=1}^{N} m_j\delta_{x_j}\otimes\delta_{P_j}$.
The associated mass measure is $\V=\sum_{j=1}^{N} m_j\delta_{x_j}$.
\end{enumerate}

\begin{remk}
From now on, given a smooth submanifold $\cM$, we denote the canonical associated varifold $\mathbf{v}(\cM,1)$ as $M: =\mathcal{H}^d\res\cM\otimes \delta_{T_y\cM}$ with mass measure $\| M\| : =\mathcal{H}^d\res\cM$.
\end{remk}

Each varifold in the three classes above decomposes as a generalized tensor product of measures. This property is actually true for any varifold, as stated in the result below that is a consequence of a general disintegration result, see for instance \cite[Theorem 2.28]{afp}.
\begin{prop}[Disintegration] \label{prop:disintegration}
Let $V$ be a $d$--varifold in $\R^n$. There exists a family $(\nu_x)_x$ of probability measures in $\G$, defined for $\| V \|$--a.e. $x \in \R^n$, such that $V = \| V \| \otimes \nu^x$ in the following sense 
\begin{equation*}
    \forall \phi \in \xC_c(\R^n \times \G), \quad \int \phi \: dV = \int_{x \in \R^n} \int_{S \in \G} \phi(x,S) \: d \nu^x(S) \: d \| V \|(x) \: .
\end{equation*}
\end{prop}

We recall now the definition of the bounded Lipschitz distance for Radon measures, that will be used to compare varifolds regardless of their type (rectifiable, point clouds...).
\begin{dfn}(Bounded Lipschitz distance)
Let $(X,d)$ be a locally compact separable metric space. The bounded Lipschitz distance between two finite Radon measures $\mu$, $\nu$ on $X$ is
\begin{equation}\label{dfnboundedlip}
\Delta(\mu,\nu) := \sup\limits \Big\lbrace \Big| \int_X \phi(x) d\mu(x) - \int_X \phi(x) d\nu(x) \, \Big|,\; \phi \in \xC^{0}(X,\R_+),\; \max\lbrace \|\phi\|_{\infty}, \lip(\phi) \rbrace \leq1 \Big\rbrace
\end{equation}
\end{dfn}
The notion of convergence we will be dealing with throughout this paper is the weak-$*$ convergence of Radon measures:
\begin{dfn}[Weak-$*$ convergence]
 Let $(X,d)$ be a locally compact and separable metric space and $(\mu_i)_{i\in\N}$, $\mu$ Radon measures on $X$. We say that
 $(\mu_i)_i$ converges weakly-$*$ to $\mu$, and we write $\mu_i \xrightharpoonup[i \to \infty]{\ast} \mu$, if
\begin{equation}
\forall  \phi \in \xC_c(X,\R_+),\quad  \int_X\phi \, d\mu_i \rightarrow \int_X\phi \, d \mu.
\end{equation}
\end{dfn}
The bounded Lipschitz distance provides a local metrization of the weak-$*$ convergence, as stated in the following result, see~\cite[Thm 5.9]{st}.
\begin{prop}\label{metrization1}
 Let $(X,d)$ be a locally compact separable metric space and $(\mu_i)_{i \in \N}$, $\mu$ finite Radon measures with support included in a compact subset of $X$. 
Then 
 \begin{center}
  $(\mu_i)$ \text{converges weakly-$*$ to} $\mu$ \quad $\Longleftrightarrow$  $\quad \Delta(\mu_i,\mu) \xrightarrow[i\rightarrow \infty]{} 0.$
 \end{center}
\end{prop}

The notion of push-forward of a varifold~\cite[Sec. 1.4]{ton} is crucial for our construction:
\begin{dfn}(Push-forward of a varifold)\label{dfn:varifoldpush-forward}
Let $V$ be a $d$-varifold in $\R^n$ and $f$ a $\xC^1$ diffeomorphism of $\R^n$. The push-forward of $V$ by $f$ is the varifold $ f_{\#}V$ defined for every $\phi \in C_{c}(\R^n\times\G, \R)$ by
\begin{equation} \label{eq:varifoldpush-forward}
  f_{\#}V(\phi):=\int_{\R^n\times\G} \phi(f(x),Df(x)(S))J_Sf(x)\, dV(x,S),
\end{equation}
where $Df(x)(S)$ is the image of $S$ in $\G$ by the linear isomorphism $Df(x)$, and the tangential Jacobian $J_Sf(x)$ is the determinant of the isomorphism $Df(x)$ from $S$ to $Df(x)(S)$, defined as follows: let $\tilde{S}= \left( s_1 | \dots | s_d \right)^t \in \cM_{d,n}$, with $\lbrace s_i \rbrace_{i=1}^d$ an orthonormal basis of $S$, and $Y = Df(x)\tilde{S}^t$,  where $Df(x)$ is identified with its $n\times n$ matrix in the canonical basis of $\R^n$, then
\begin{equation} \label{eq:JSf}
J_Sf(x) := \det\left( Y^tY \right)^{\frac{1}{2}}.
\end{equation}
\end{dfn}
\begin{remk}
Using the notations above, the orthogonal projection $P$ onto the space $Df(x)(S)$ satisfies (see \cite[p. 184]{stra})
\begin{equation} \label{eq:DfS}
 P = Y(Y^tY)^{-1}Y^t.
\end{equation}
\end{remk}

\begin{remk}
Let $\cM \subset \R^n$ be a $d$--submanifold and $f: \R^n \rightarrow \R^n$ a $\xC^1$ diffeomorphism of $\R^n$. By the area formula, 
\begin{equation}\label{intro_def:push_mass}
\int_{f(\cM)} \phi(y) \, d\cH^d(y) = \int_{\cM} \phi(f(x)) J_{T_x\cM}f(x) \, d\cH^d(x) \quad \forall \phi \in C^1_0(\R^n,\R_+),
\end{equation}
where $J_{T_x\cM}f(x)$ is the tangential Jacobian of $f$ with respect to $T_x\cM$. Note that if $V_\cM = v(\cM,1) = \cH^d \res \cM \otimes \delta_{T_x \cM}$ denotes the unit multiplicity rectifiable $d$--varifold associated with $\cM$ then, according to Definition~\ref{dfn:varifoldpush-forward} and \eqref{intro_def:push_mass},
\[
 f_\# V_{\cM} = V_{f(\cM)} = \cH^d \res f(\cM) \otimes \delta_{T_y f(\cM)} \: . 
\]
More generally, if a nonnegative multiplicity function $\theta \in L_{loc}^1(\cH^d\res\cM)$ is used, then
\[                                                                                                                                                                                                                                                    
f_\# v(\cM, \theta) = \theta ( f^{-1} (y)) \cH^d \res f(\cM) \otimes \delta_{T_y f(\cM)}.
\]
In particular, 
\[
 \left\| f_\# v(\cM, \theta) \right\| = (\theta \circ f^{-1}) \cH^d \res f(\cM) = v(f(\cM), \theta \circ f^{-1}).
\]
%
\end{remk}

The following lemma characterizes the composition of two push-forwards.
\begin{lemma}\label{basiclemma}
Let $V \in V_d(\R^n)$ and let $f$, $g$ be two $\xC^1$ diffeomorphisms of $\R^n$. Then
$$ f_{\#}(g_{\#}V) = (f\circ g)_{\#}(V).$$
\end{lemma}
\begin{proof}
Let $\phi \in \xC_c(\R^n\times \G,\R_+)$. Then
\begin{equation*}
\begin{split}
\int \phi \: d(f_{\#}(g_{\#}V)) & = \int \phi \left(f(y),Df(y)(T) \right) \: J_Tf(y) \: d(g_{\#}V)(t,y) \\& = \int \phi\left((f(g(x)),Df(g(x))(Dg(x)(S))\right) \: J_{Dg(x)(S)}f(g(x)) \: J_Sg(x) \: dV(x,S) \\
& = \int \phi\left( f \circ g(x) , D(f\circ g)(x)(S) \right) \: J_S(f\circ g)(x) \: dV(x,S)
\end{split}
\end{equation*}
where we used that for $(x,S) \in \R^n \times \G$,
\begin{equation*}
 Df(g(x))(Dg(x)( S)) = D(f\circ g)(x)(S) \quad \text{and} \quad J_{Dg(x)(S)}f(g(x))J_Sg(x) = J_S(f\circ g)(x)
\end{equation*}
thanks to the multiplicative property of the determinant.
\end{proof}
Given $X \in C^1(\R^n, \R^n)$ and $t> 0$, we define
\begin{equation*}
 f_t={\rm id}+tX \: .
\end{equation*}
For $t$ small enough, $f_t$ is a $\xC^1$ diffeomorphism and we can examine the infinitesimal change of the mass measure of a varifold $V\in V_d(\R^n)$ pushed by $f_t$.
We have the following formula~\cite{simon}:
\begin{equation*}
 \partial_t \|(f_{t})_{\#}V\|(\R^n)_{|t=0}  = \int_{\R^n \times \G} \mdiv_S X(x) \, dV(x,S) \,
\end{equation*}
where $\mdiv_SX := \tr(S \,DX)$ is the tangential divergence. This computation motivates the definition of the first variation of a varifold.
\begin{dfn}[First variation]
The \textit{first variation} of a varifold $V\in V_d(\R^n)$ of finite mass is the map 
\begin{equation}\label{firstvar}
\delta V:\; X \in C^1(\R^n, \R^n)\longmapsto \int_{\R^n \times \G} \mdiv_S X(x)\, dV(x,S).
\end{equation}
\end{dfn}
\noindent Given $\phi \in \xC^1(\R^n, \R)$ and $X \in \xC^1(\R^n, \R^n)$,  we denote: 
\begin{equation}\label{wfirstvar}
     \begin{split}\delta (V,\phi)(X) &:= \int_{\R^n \times \G} \phi(x) \,\mdiv_S X(x) \, dV(x,S) + \int_{\R^n\times\G} \nabla\phi(x) \cdot X(x) \,dV(x,S)
     \\& = \delta V (\phi X) + \int_{\R^n \times \G } X(x) \cdot S^{\perp} (\nabla\phi(x)) \, dV(x,S).
\end{split}\end{equation}
where $S^{\perp}$ denotes the orthogonal projection on the space orthogonal to $S$. 
\begin{dfn}[Weighted first variation]
Let $V\in V_d(\R^n)$ of finite mass. The map $\delta (V,\cdot)(\cdot)$  from $(\xC^1(\R^n,\R), \xC^1(\R^n,\R^n))$ onto $\R$ defined above is called the \textit{weighted first variation} of $V$.
\end{dfn}
\begin{remk} \label{remk:wfirstvar}
Note that $\delta (V,\cdot)(\cdot)$ is bilinear and
\begin{equation*}
 \left|\delta (V,\phi)(X) \right| \leq n \| X \|_{\xC^1} \| V \| (\R^n) \: \| \phi \|_{\xC^1} \: , \quad \text{for all } \phi \in C^1(\R^n,\R), \: X \in \xC^1(\R^n , \R^n).
\end{equation*}
\end{remk}

\noindent Given a closed (compact and without boundary) smooth $d$-manifold $\cM$, we recall the notation $M: =\mathcal{H}^d\res\cM\otimes \delta_{T_x\cM}$. For $X \in C^1(\R^n, \R^n)$,
\begin{equation}\label{intro_mcchar}
\delta M (X) = \int_{\R^n \times \G} \hspace{-15pt} \mdiv_S X(x) \,dM(x,S)=
\int_{\cM} \mdiv_{T_x\cM} X (x) \,d\cH^d(x)=-\int_{\cM} H(x,\cM) \cdot X(x) \,d\mathcal{H}^d(x),
\end{equation}
where $H(\cdot,\cM)$ is the mean curvature vector of $\cM$ defined in~\eqref{intro_mcdfn}.
 
\par
When the first variation $\delta V$ of a varifold $V\in V_d(\R^n)$ is locally bounded, the Riesz representation theorem implies that $\delta V$ can be represented by a vector Radon measure, which we continue to denote by $\delta V$ for simplicity. Then, the Radon-Nikodym decomposition theorem implies the existence of a vector $H(\cdot,V) \in (L^1(\R^n,\V))^n$ such that 
\begin{equation}
 \delta V = - H(\cdot,V) \|V\| + \delta^s V,\label{gmc}
\end{equation}
where $ \delta^s V$ is a vector measure singular with respect to $\|V\|$. The analogy with  \eqref{intro_mcchar} motivates the following definition.

\begin{dfn}[Generalized mean curvature]
The vector $H(\cdot,V)$ defined in \eqref{gmc} for a varifold $V\in V_d(\R^n)$ with locally bounded first variation is called the generalized mean curvature of $V$ \end{dfn}

Not every varifold $V$ has locally bounded first variation, for instance if $V$ is a point cloud varifold. For such a varifold, the generalized mean curvature is not defined. In this situation, the alternative is to define an approximate mean curvature, as in~\cite{blm1} or, with better regularity properties (which is what we shall need here), as in~\cite{brakke,kt}. The latter definition of approximate mean curvature is recalled at the beginning of the next section.
%

%

\section{Definition and stability of time-discrete approximate mean curvature flows}\label{def_time_discrete}

In this section, we define, for a given time subdivision and approximation scale, a time-discrete approximate mean curvature flow starting from any varifold $V$ with finite mass. The construction relies on iterated push-forwards of $V$  by diffeomorphisms of the form ${\rm id} + \tau h_\epsilon$, where $\e>0$ is the approximation scale, $\tau$ is the local time step, and $h_\epsilon$ is the approximate mean curvature of $V$ at scale $\e$ introduced in \cite{kt} after \cite{brakke}, obtained using suitable mollifications by a regularization kernel. The definition and properties of the regularization kernel are recalled in Subsection~\ref{subsec:kernel}, and the definition of $h_\e$ in Subsection~\ref{subsec:basicHeps}. Subsections~\ref{subsec:timeDiscreteMCF}, \ref{subsec:stabilityInitialData}
and \ref{subsec:stabilityTimeSubdivision} are devoted to the definition of the time-discrete approximate mean curvature flow, and to the study of its stability properties with respect to initial data and time subdivision.

\subsection{Preliminaries on the regularization kernel}
\label{subsec:kernel}
Let $\psi\in C^{\infty}({\mathbb R}^{n})$ be a radially symmetric function such that:
\begin{equation}\label{defpsi}
\begin{split}
&\psi(x)=1\mbox{ for }|x|\leq 1/2,\qquad \psi(x)=0\mbox{ for }|x|\geq 1, \\
&0\leq \psi(x)\leq 1,\quad  |\nabla\psi(x)|\leq 3,\quad \|\nabla^2\psi(x)\|\leq 9,\qquad \mbox{ for all }x\in
{\mathbb R}^{n}.
\end{split}
\end{equation}
Define for each $\e\in(0,1)$:
\begin{equation}\label{defkernel}
\hat\Phi_{\e}(x):=\frac{1}{(2\pi\e^2)^{\frac{n}{2}}}\exp\Big(-\frac{|x|^2}{2\e^2}\Big)\qquad\mbox{and}\qquad
\Phi_{\e}(x):=c(\e)\psi(x)\hat\Phi_{\e}(x),
\end{equation}
where $\displaystyle c(\e):=\frac{1}{\int_{\R^n} \psi(x)\hat{\Phi}_{\e}(x)dx}$ so that, because $\ds\int_{\R^n}\hat{\Phi}_{\e}(x)dx=1$, we have that
\begin{equation} \label{eq:PhiEpsL1}
\int_{{\mathbb R}^{n}}\Phi_{\e}(x)\, dx=1.
\end{equation}
Then, using $\psi \leq 1$,
\begin{equation*}
 c(\e)=\frac{1}{\int_{\R^n} \psi(x)\hat{\Phi}_{\e}(x)dx} \geq \frac{1}{\int_{\R^n}\hat{\Phi}_{\e}(x)dx}=1.
\end{equation*}
Also, as $\psi =1$ on $[0,\frac12]$,
\begin{equation*}
 \int_{\R^n} \psi(x)\hat{\Phi}_{\e}(x)dx \geq \int_{B_\frac12(0)}\hat{\Phi}_{\e}(x)dx.
\end{equation*}
By the change of variables $y = \e^{-1} x $ we obtain 
\begin{equation}\label{defc}
\int_{B_\frac12(0)}\hat{\Phi}_{\e}(x)dx=\int_{B_\frac{1}{2\e}(0)}\hat{\Phi}_{1}(y)dy \geq \int_{B_\frac{1}{2}(0)} \hat{\Phi}_{1}(y)dy =: c^{-1}.
\end{equation}
$c$ is a constant depending only on $n$, and by definition, $1 \leq c(\e) \leq c$.

The kernel $\Phi_\e$ has a remarkable property: its derivatives are bounded by a power of $\e$ times the kernel plus an exponentially small term (see \cite[Lemmas 4.13 and 4.14]{kt}). This property is the key ingredient for the computations of \cite[Section 5]{kt}. 

In the following lemma we list some of the properties of the kernel $\Phi_{\e}$, similar to the estimates in \cite[Lemma 4.13]{kt}.
\begin{lemma}[Kernel properties]\label{kerprop}
Let $\e\in(0,1)$ and $\Phi_{\e}$ be defined as in \eqref{defkernel}. There exists a constant $c_0$ depending only on $n$ such that
\begin{equation}\label{kerprop1}
 |\nabla\Phi_{\e}|\leq \e^{-2}\Phi_{\e}+c_0\chi_{B_1(0)},
\end{equation}
\begin{equation}\label{kerprop2}
 | \nabla^2\Phi_{\e} | \leq 2\e^{-4}\Phi_{\e}+2c_0\chi_{B_1(0)}.
\end{equation}
As a consequence,
\begin{equation}\label{kerprop7}
\| \nabla\Phi_{\e} \|_{L^1} \leq (1+\omega_n c_0) \e^{-2}  \quad \text{and} \quad \| \nabla^2\Phi_{\e} \|_{L^1} \leq 2(1+\omega_n c_0) \e^{-4} \: ,
\end{equation}
and
\begin{equation}\label{lipker}
 \lip(\Phi_{\e})\leq (c(2\pi)^{-\frac{n}{2}} + c_0 ) \e^{-n-2} \, ,  \quad \lip(\nabla\Phi_{\e}) \leq 2(c(2\pi)^{-\frac{n}{2}} + c_0 ) \e^{-n-4}.
\end{equation}
\end{lemma}
\begin{proof}
Define
\begin{equation}\label{defc_0}
c_0 := \sup\limits_{\e\in(0,1)}c(\e)\frac{9\e^{-2-n}}{(2\pi)^{n/2}}\exp\left( -\frac{1}{8\e^2}\right) < \infty.
\end{equation}
By \eqref{defkernel} for  all $x\in \R^n$
\begin{equation}\label{kerprop3}
  \nabla\Phi_{\e}(x)= -\e^{-2}\Phi_{\e}(x)x+c(\e)\hat{\Phi}_{\e}(x)\nabla\psi(x)
\end{equation} 
$\Phi_{\e}$ is supported on $B_1(0)$ therefore for all $x\in\R^n $
$$\displaystyle \big|\e^{-2}\Phi_{\e}(x)x\big|\leq \e^{-2}\Phi_{\e}(x)$$
also by construction $\nabla\psi=0$ on $[0,\frac12] \cup [1,\infty)$ and $|\nabla\psi| \leq 3$ on $[\frac{1}{2},1]$, this yields for all $x\in\R^n$
\begin{equation}\label{kerprop4}
\begin{split}
|c(\e)\hat{\Phi}_{\e}(x)\nabla\psi(x)| & \leq   3c(\e) \frac{\e^{-n}}{(2\pi)^{\frac{n}{2}}}\sup\limits_{\frac12\leq|x|\leq1}\exp\left(-\frac{|x|^2}{\e^2} \right) \chi_{B_1(0)}(x)
\\& \leq  3c(\e)\frac{\e^{-n}}{(2\pi)^{\frac{n}{2}}}\exp\left(-\frac{1}{8\e^2}\right) \chi_{B_1(0)}(x) \leq c_0\chi_{B_1(0)}(x) 
\end{split}
\end{equation}
this proves \eqref{kerprop1}. For \eqref{kerprop2}, differentiating \eqref{kerprop3} gives for all $x\in\R^n$
\begin{equation*}
\begin{split}
 \nabla^2\Phi_{\e}(x)&=-\e^{-2}x\otimes\left(-\e^{-2}\Phi_{\e}(x)x+c(\e)\hat{\Phi}_{\e}(x)\nabla\psi(x)\right) -\e^{-2}\Phi_{\e}(x)I_{n}
\\& +c(\e)\left(\hat{\Phi}_{\e}(x)\nabla^2\psi(x)-\hat{\Phi}_{\e}(x)\e^{-2}\nabla\psi(x)\ \otimes x\right)
\\& = \e^{-4}\Phi_{\e}(x) x \otimes x-2\e^{-2}c(\e)\hat{\Phi}_{\e}(x) x\otimes\nabla\psi(x)-\e^{-2}\Phi_{\e}(x)I_{n}
\\&+c(\e)\hat{\Phi}_{\e}(x)\nabla^2\psi(x).
 \end{split}\end{equation*}
We know that $\|v \otimes w \| \leq \|v\|w\|$ for every two vectors $v$ and $w$, the fact that $\Phi_{\e}$ and $\psi$ are supported on $[0,1]$ implies, for all $x\in\R^n$
\begin{equation*}
\begin{split}
\|\nabla^2\Phi_{\e}(x)\| 
&\leq \e^{-4}\Phi_{\e}(x)+2\e^{-2}|c(\e)\hat{\Phi}_{\e}(x)\nabla\psi(x)|+\Phi_{\e}(x)+\|c(\e)\hat{\Phi}_{\e}(x)\nabla^2\psi(x)\|
\\& \leq 2\e^{-4}\Phi_{\e}(x)+2\e^{-2}|c(\e)\hat{\Phi}_{\e}(x)\nabla\psi(x)|+\|c(\e)\hat{\Phi}_{\e}(x)\nabla^2\psi(x)\|.
\end{split}
\end{equation*}
Similarly to \eqref{kerprop4}, we have for all $x \in \R^n$
\begin{equation}\label{kerprop5}
\begin{split}
2\e^{-2}|c(\e)\hat{\Phi}_{\e}(x)\nabla\psi(x)| & \leq   c(\e)6\frac{\e^{-2-n}}{(2\pi)^{\frac{n}{2}}}\sup\limits_{\frac12\leq|x|\leq1}\exp\left(-\frac{|x|^2}{\e^2} \right)\chi_{B_1(0)}(x)
\\& \leq  c(\e)6\frac{\e^{-2-n}}{(2\pi)^{\frac{n}{2}}}\exp\left(-\frac{1}{8\e^2} \right) \chi_{B_1(0)}(x) \leq c_0 \chi_{B_1(0)}(x) 
\end{split}
\end{equation}
and
\begin{equation}\label{kerprop6}
 \begin{split}
  \|c(\e)\hat{\Phi}_{\e}(x)\nabla^2\psi(x)\| & \leq   c(\e)9\frac{\e^{-n}}{(2\pi)^{\frac{n}{2}}}\sup\limits_{\frac12\leq|x|\leq1}\exp\left(-\frac{|x|^2}{\e^2} \right)\chi_{B_1(0)}(x)
\\& \leq  c(\e)9\frac{\e^{-n}}{(2\pi)^{\frac{n}{2}}}\exp\left(-\frac{1}{8\e^2} \right)\chi_{B_1(0)}(x) \leq c_0\chi_{B_1(0)}(x).
 \end{split}
\end{equation}
Finally, \eqref{kerprop5} and \eqref{kerprop6} imply
\begin{equation*}
 \|\nabla^2\Phi_{\e}\|\leq 2\e^{-4}\Phi_{\e} + 2c_0\chi_{B_1(0)}
\end{equation*}
and this concludes the proof of \eqref{kerprop2}. For the $L^1$-estimates, we write using \eqref{kerprop1}
\begin{equation}\begin{split}
 \|\nabla\Phi\|_{L^1} &= \int_{\R^n} |\nabla\Phi_{\e}(x)| dx \leq  \int_{\R^n}\left(\e^{-2}\Phi_{\e}(x)+ c_0 \chi_{B_1(0)}(x)\right)dx
 \\& \leq \e^{-2}\int_{\R^n}\Phi_{\e}(x) \,dx + c_0 \int_{B_1(0)}1 \,\,dx
 \\& \leq \e^{-2}+\omega_nc_0 \leq (1+\omega_nc_0)\e^{-2}.
\end{split}\end{equation}
Similarly, we prove $\|\nabla^2\Phi\|_{L^1}\leq 2(1+c_0\omega_n)\e^{-4}$, and this completes the proof of \eqref{kerprop7}.\\
For the Lipschitz constant, we use the definition of $\Phi_{\e}$ (see \eqref{defkernel} to get: 
\begin{equation*}
 \lip(\Phi_{\e})\leq (c(2\pi)^{-\frac{n}{2}} + c_0 ) \e^{-n-2} \, ,  \quad \lip(\nabla\Phi_{\e}) \leq 2(c(2\pi)^{-\frac{n}{2}} + c_0 ) \e^{-n-4}
\end{equation*}
and this finishes the proof of \eqref{lipker} and Lemma \ref{kerprop}.
\end{proof}

\subsection{Approximate mean curvature: definition and properties}
\label{subsec:basicHeps}

We can now define the approximate mean curvature vector at scale $\e>0$ for any varifold $V$.
\begin{dfn}[Approximate mean curvature]
The approximate mean curvature of a varifold $V\in V_d(\R^n)$ at $x\in \R^n$ for the approximation scale $\e>0$ is 
\begin{equation}\label{dfnhepsilon}
h_{\e}(x,V)=(\Phi_{\e}\ast \tilde{h}_{\e}(\cdot,V))(x)\,, \quad \text{where} \,\,\, \tilde{h}_{\e}(y,V)=-\frac{(\delta V \ast \Phi_{\e})(y)}{(\|V\| \ast \Phi_{\e})(y)+\e} \,\, \text{for any $y \in\R^n$}.
\end{equation}
\end{dfn}
The double convolution guarantees the decay of the mass (up to a small error), as we will see in \eqref{massdecrease6}, and later in Remark \ref{remk:massdecrease4}. It also reduces the calculation of $\|h_{\e}\|_{C^2}$ to that of $\|\tilde{h}_{\e}\|_{\infty}$ and $\|\Phi_{\e}\|_{C^2}$, thereby avoiding the differentiation of the fraction $\tilde{h}_{\e}$ (see Proposition \ref{hepsilonbound} for details).

The following property is a mere adaptation of \cite[Lemma 5.1]{kt}: we bound the $\xC^2$--norm of the approximate mean curvature for $\e\in(0,1)$, but we impose here no smallness requirement on $\e$ contrary to the original statement.
\begin{prop}[$\xC^2$ boundedness of $h_{\e}$] \label{hepsilonbound}
There exists a constant $c_1 \geq 2$ depending only on $n$ with the following property: for any $\e\in(0,1)$ and $M \in [1, +\infty )$, if $V \in V_d(\R^n)$ is a $d$--varifold with total mass $\V(\R^n) \leq M$, then
\begin{equation}\label{hepsilonbound1}
 \| \tilde{h}_{\e}(\cdot,V) \|_{\infty} \leq c_1 M \e^{-2}, \quad  \| h_{\e}(\cdot,V) \|_{\infty} \leq c_1 M \e^{-2},
\end{equation}
\begin{equation}\label{hepsilonbound2}
 \| Dh_{\e}(\cdot,V) \|_{\infty} \leq c_1 M \e^{-4},
\end{equation}
\begin{equation}\label{hepsilonbound3}
  \| D^2h_{\e}(\cdot,V) \|_{\infty} \leq c_1 M \e^{-6}.
\end{equation}
\end{prop}
\begin{proof}
Let $\e \in (0,1)$, $M \geq 1$ and let $V$ be a $d$--varifold satisfying $\V(\R^n) \leq M$. We start with the proof of \eqref{hepsilonbound1}, setting $c_1=2(1+\omega_n c_0)(1+c_0)$. 
For any $\phi \in \xC(\R^n,\R)$, we recall that $\Phi_\e$ is radial and then
\begin{equation*}
 (\Phi_{\e}  \ast \V)(\phi) = \V(\Phi_{\e}  \ast \phi) = \int_{\R^n}\int_{\R^n}\Phi_{\e} (y-x)\, \phi(y) \, dx \, d\|V\|(y)
\end{equation*}
hence we can associate the convoluted measure $\Phi_{\e}  \ast \V$ with the function 
\begin{equation}\label{convmass}
 x\in\R^n\mapsto (\Phi_{\e}\ast \V)(x) := \int_{\R^n}\Phi_{\e}(y-x) \, d\|V\|(y).
\end{equation}
Similarly, 
\begin{equation*}
\begin{split}
 (\Phi_{\e}\ast \delta V)(\phi) 
 = \delta V(\Phi_{\e}\ast\phi) &=\int_{\R^n\times\G} \int_{\R^n} S(\nabla\Phi_{\e}(y-x))\phi(y) \, \, dxdV(y,S)\\& = \int_{\R^n\times\G} \left(\int_{\R^n} S(\nabla\Phi_{\e}(y-x))\phi(y) \, dV(y,S)\right)\,dx.
 \end{split}
\end{equation*}
hence the convoluted first variation $\Phi_{\e}\ast \delta V$ can be associated with the function 
\begin{equation}\label{convvar}
x\in\R^n\mapsto (\Phi_{\e}\ast \delta V)(x) := \int_{\R^n} S(\nabla\Phi_{\e}(y-x)) \, dV(y,S) dx.
\end{equation}

Using that, for all $S \in \G$, $\| S \| \leq 1$, we have
\[
\left| (\Phi_{\e} \ast \delta V)(x) \right| = \left| \int_{\R^n \times \G} S(\nabla\Phi_{\e}(x-y)) \: dV(y,S) \right| \leq  \int_{\R^n}|\nabla\Phi_{\e}(x-y)| \: d \V (y) \: .                                                                                                 
\]

Therefore, applying \eqref{kerprop1} and then \eqref{convmass} we obtain
\begin{equation} \label{eqhepsilonbound6}
 |(\Phi_{\e}\ast \delta V)(x)| \leq \int_{\R^n} \bigl( \e^{-2}\Phi_{\e}(x-y) +  c_0\chi_{B_1(0)}(x-y) \bigr) d \V(y)
\leq \e^{-2} (\Phi_{\e}\ast \V)(x) + c_0 M \: . 
\end{equation}
 It remains to write the definition \eqref{dfnhepsilon} of $\tilde{h}_{\e}(\cdot,V)$ and apply \eqref{eqhepsilonbound6} to infer 
 \begin{equation} \label{eqhepsilonbound7}
 \| \tilde{h}_{\e}(\cdot,V) \|_{\infty}=\sup\limits_{x \in \R^n}\frac{|(\Phi_{\e}\ast \delta V)(x)|}{\Phi_{\e}\ast \V(x)+ \e} \leq \e^{-2} + \e^{-1}c_0 M \leq (1+c_0) M \e^{-2}  \: .
\end{equation}
We can now use $h_{\e}:=\Phi_{\e} \ast \tilde{h}_{\e}$ (see  \eqref{dfnhepsilon}), $\| \Phi_\e \|_{L^1} = 1$, and \eqref{eqhepsilonbound7} to obtain
\begin{equation} \label{eqhepsilonbound8}
\left\| h_{\e}(\cdot,V) \right\|_{\infty}  \leq  \|\Phi_{\e} \|_{L^1} \: \| \tilde{h}_{\e}(\cdot,V) \|_{\infty} \leq (1+c_0) M \e^{-2} \: ,
\end{equation}
and noting that $(1+c_0) \leq c_1$ concludes the proof of \eqref{hepsilonbound1}.\\
We similarly have both $Dh_{\e}(\cdot,V)=\nabla\Phi_{\e}\ast \tilde{h}_{\e}$ and
$D^2h_{\e}(\cdot,V)=\nabla^2\Phi_{\e}\ast \tilde{h}_{\e}$ so that applying \eqref{kerprop7} together with \eqref{eqhepsilonbound7} concludes the proof of Proposition~\ref{hepsilonbound} as follows:
\begin{align*}
& \left\| Dh_{\e}(\cdot,V) \right\|_{\infty}  \leq \|\nabla\Phi_{\e} \|_{L^1} \: \|\tilde{h}_{\e}(\cdot,V) \|_{\infty} \leq (1+c_0\omega_n)(1+c_0)M \: \e^{-4} \leq c_1 M \: \e^{-4} \: ,\\
& \| D^2h_{\e}(\cdot,V) \|_{\infty}  \leq \|\nabla^2\Phi_{\e} \|_{L^1} \: \|\tilde{h}_{\e}(\cdot,V) \|_{\infty} \leq 2(1+c_0\omega_n)(1+c_0)M \: \e^{-6} \leq c_1 M \: \e^{-6} \:.
\end{align*}
\end{proof}
\subsection{Definition of a time-discrete approximate mean curvature flow}
\label{subsec:timeDiscreteMCF}
The goal of this subsection is to define a time-discrete approximate mean curvature flow (see Definition~\ref{damcf}) starting from an initial varifold $V_0 \in V_d(\R^n)$, for a given time subdivision $\cT$ and approximation parameter~$\e$. Such a definition relies on iterating push-forwards, starting from the initial varifold $V_0$, with velocity equal to the approximate mean curvature vector.
To this end, we first investigate the effect of a single push-forward:
in Proposition~\ref{mcfmotion1}, we  derive an expansion with respect to $\Delta t$ of the push-forward of the mass of a varifold under the map $f={\rm id}+\Delta t \, h_{\e}$. Computations rely on the Taylor expansion of the tangential Jacobian (see Lemma~\ref{detprop}), from which follows estimate \eqref{appbrakke}.
It is then possible to prove that the mass of the push-forward varifold decays up to a small error $\del$ (see \eqref{massdecrease6}), 
hence allowing to iterate push-forwards for suitable time steps (see condition \eqref{smallstep}) and resulting in Definition~\ref{damcf}.

Given $\e \in (0,1)$ and $V \in V_d(\R^n)$, we introduce the notation
\begin{equation*}
f_{\e,V} = {\rm id} + \Delta t \: h_\e (\cdot , V) \: ,
\end{equation*}
Depending on the context, we may drop the $\e$ or $V$ index dependency for simplicity.

\begin{prop}\label{mcfmotion1}
Let $\e, \, \del \in (0,1)$ and $M \geq 1$. Let $V \in V_d(\R^n)$ satisfy $\| V \|(\R^n) \leq M$ and $S \in \G$. There exists a constant $c_3 > 0$ that depends only on $n$ and such that, if
\begin{equation}\label{smallstep0}
c_3 M \Delta t \leq \e^4 
\end{equation}
then $f = {\rm id} + \Delta t \: h_\e (\cdot , V)$ is a diffeomorphism, and
\begin{equation}\label{eq:JSfbounds}
\forall (x,S) \in \R^n \times \G, \quad J_S f(x)  \in \Big[ \frac12 , \frac32 \Big] \cap  \big[ 1 -  c_3 \Delta t \| Dh_{\e}\|_{\infty}, 1 +  c_3 \Delta t \| Dh_{\e} \|_{\infty}  \big] \: ,
\end{equation} 
\begin{equation}\label{appbrakke}
\forall \phi \in \xC^2 (\R^n, \R_+), \quad    \big| \|f_{\#}V\|(\phi) - \V(\phi) - \Delta t \, \delta (V,\phi)(h_{\e}(\cdot,V)) \big| \leq c_3 \|\phi\|_{\xC^2} M^3 (\Delta t)^2 \e^{-8} \: ,
\end{equation}
and
\begin{equation}\label{massdecrease1}
\delta V(h_\e(\cdot, V)) = - \int_{\R^n}  \frac{|(\Phi_{\e}\ast \delta V)(y)|^2}{(\Phi_{\e}\ast\V)(y)+\e} dy \leq 0 \: .
\end{equation} 
If we furthermore assume $c_3 M^3 \del \leq \e^8$ (that implies \eqref{smallstep0}) then, 
\begin{equation}\label{massdecrease6}
\| f_{\#}V \| (\R^n) \leq \| V \|(\R^n) + \del \: .
\end{equation}
\end{prop}


%
\begin{proof}
Throughout the proof, we may increase $c_3 > 0$ to meet the requirements, provided that we respect the exclusive dependency in $n$.
Let $\e, \, \del \in (0,1)$, $M \geq 1$ and $V \in V_d(\R^n)$ satisfy $\V (\R^n) \leq M$. As $V$ is fixed, we write $h_\e$ for $h_\e (\cdot, V)$ hereafter as well as $f = {\rm id} + \Delta t \, h_\e$.

\noindent {\bf Step 1:}
We first prove that $f$ is a diffeomorphism under the condition \eqref{smallstep0}. To do so,
we only need to check the hypothesis of Lemma \ref{pre_lem:diff} (see below) with $h=h_{\e}$.
From \eqref{hepsilonbound1} and \eqref{smallstep0} we can infer that 
\begin{equation*} \del \, \| h_{\e}\|_{\infty} \leq c_1 M \del \e^{-2} \leq  \frac{c_1}{c_3} < 1.
\end{equation*}
From \eqref{hepsilonbound2} and \eqref{smallstep0}
\begin{equation} \label{eq:diffeoProof1}
\Delta t \, \| Dh_{\e} \|_\infty \leq  c_1 M \del \e^{-4} \leq \frac{c_1}{c_3} < 1,
\end{equation}
and we can then apply \eqref{eq:detprop1} (with $k=n$ and $Q =\Delta t \: Dh_{\e}$) together with \eqref{eq:maxEuclNorms}, \eqref{smallstep0} and \eqref{hepsilonbound2} so that for all $x\in\R^n$,
\begin{equation*}
 \left| Jf(x) -1 \right| = \left| \det (I_n + \Delta t \: Dh_{\e}(x)) - \det(I_n) \right|  \leq c_2 \Delta t | Dh_{\e}(x)|_\infty  \leq \frac{c_1 c_2}{c_3}< 1
\end{equation*}
and $f$ is a diffeomorphism of $\R^n$ thanks to Lemma~\ref{pre_lem:diff}.

\noindent {\bf Step 2:} Let $(x,S) \in \R^n \times \G$,
we now prove \eqref{eq:JSfbounds} and
\begin{equation} \label{eq:JSfboundsOrder2}
 |J_Sf(x) -1 - \Delta t \, \mdiv_S(h_{\e}(x))| 
\leq c_2  (\Delta t \, \| Dh_{\e} \|_{\infty})^2.
\end{equation}
Let us write $\tilde{S}=\left( \tau_1 | \dots |\tau_d \right)^t \in \cM_{d,n}$ where $\lbrace \tau_i \rbrace_{i=1}^d$ is an orthonormal basis of $S$.
We recall that we also denote by $S$ the orthogonal projector onto the subspace $S$. For clarity, we denote by $\circ$ the matrix product in the calculations below. By construction,
\begin{equation*} 
\tilde{S}\circ \tilde{S}^t = I_d \in \cM_d \quad \text{and} \quad \tilde{S}^t \circ \tilde{S} = S \in \cM_n.
\end{equation*}
We recall that by definition of tangential Jacobian \eqref{eq:JSf},
\begin{equation*}
 J_S f(x) = \det\left( (( I_n + \del \, Dh_{\e}(x))\circ \tilde{S}^t)^t \circ ((I_n +\del \, Dh_{\e}(x))\circ \tilde{S}^t) \right)^\frac{1}{2}
\end{equation*}
and we can apply \eqref{detprop4} with $R=\Delta t \, Dh_{\e}(x)$ and $L=\tilde{S}$, indeed, $c_2 |R|_\infty \leq \frac{c_1 c_2}{c_3} \leq 1$ thanks to \eqref{eq:diffeoProof1}. We obtain, again using \eqref{eq:diffeoProof1},
\begin{equation} \label{eq:JSfboundsBis}
 |J_S f(x) - 1 | \leq c_2 \Delta t | Dh_{\e}(x) |_\infty \leq c_3 \Delta t \| Dh_{\e} \|_{\infty} \leq \frac{1}{2} \: ,
\end{equation}
hence proving \eqref{eq:JSfbounds}.\\
Similarly to the proof of \eqref{eq:JSfbounds}, we are allowed to use \eqref{detprop5} with $R=\Delta t \,  Dh_{\e}(x)$ and $L=\tilde{S}$. Noting that
\begin{equation*}
 \tr(Dh_{\e}(x) \circ \tilde{S}^t \circ \tilde{S}) = \tr( Dh_{\e}(x) \circ S ) = \mdiv_S(h_{\e}(x))
\end{equation*}
we can infer that
\begin{equation*}
 |J_Sf(x) -1 - \Delta t \, \mdiv_S(h_{\e}(x))| \leq c_2 (\Delta t  \, |Dh_{\e}(x)|_\infty)^2
\leq c_2  (\Delta t \, \| Dh_{\e} \|_{\infty})^2 \: .
\end{equation*}

\noindent {\bf Step 3:} We now prove \eqref{appbrakke}.
Let $\phi \in C^2(\R^n,\R_+)$ and assume $\| \phi \|_{\xC^2} < \infty$ (otherwise there is nothing to prove). 
Coming back to the definitions of push-forward varifold (Definition~\ref{dfn:varifoldpush-forward}) and weighted first variation \eqref{wfirstvar}, we have
\begin{multline} \label{eq:weightedFVproof}
\|f_{\#}V\|(\phi) - \V(\phi) - \Delta t \: \delta (V,\phi)(h_{\e}) \\ = \int_{\R^n \times \G} \Big[ \phi(f(x)) J_Sf(x) - \phi(x)- \Delta t \, \phi(x) \mdiv_S(h_{\e}(x)) - \Delta t \, \nabla\phi(x) \cdot h_{\e}(x) \Big] \: dV(x,S).
\end{multline}
Let $(x,S) \in \R^n \times \G$. We first recall that $f(x) - x = \Delta t \, h_\e(x)$ so that 
\[
 |f(x) - x| \leq \Delta t \, \| h_\e \|_\infty \leq c_1 \Delta t M \e^{-2}
\]
thanks to \eqref{hepsilonbound1}.
We can then apply Taylor's inequality to $\phi$ between $x$ and $f(x)$ to obtain 
\begin{equation} \label{eq:weightedFVproof1}
 |\phi(f(x))-\phi(x)|\leq |f(x)-x | \|\nabla\phi\|_{\infty} \leq c_1 \| \phi \|_{\xC^2}  M \Delta t \, \e^{-2} 
 \end{equation}
 and
 \begin{align}
 |\phi(f(x))-\phi(x)- \Delta t \, h_\e(x) \cdot \nabla\phi(x)| & = |\phi(f(x))-\phi(x)-(f(x)-x) \cdot \nabla\phi(x)| \nonumber \\
 & \leq \frac{1}{2} |f(x)-x|^2\|\nabla^2\phi\|_{\infty} \leq c_1^2 \| \phi \|_{\xC^2} M^2 \Delta t^2 \, \e^{-4} \: . \label{eq:weightedFVproof2}
\end{align}
Now rewriting the integrand in the right-hand side of \eqref{eq:weightedFVproof} and using \eqref{eq:diffeoProof1} to \eqref{eq:weightedFVproof2}
\begin{align*}
\left| \phi(f(x)) J_Sf(x) \right. & \left. - \phi(x)- \Delta t \: \phi(x) \mdiv_S(h_{\e}(x)) - \Delta t \: \nabla\phi(x) \cdot h_{\e}(x) \right| \\
\leq & \left|  \phi(f(x)) - \phi(x) \right| \left| J_Sf(x) -1 \right| + \phi(x) \left| J_S f(x) - 1 - \Delta t \: \mdiv_S(h_{\e}(x)) \right| \\
& + \left|  \phi(f(x)) - \phi(x) - \Delta t \: h_\e (x) \cdot \nabla\phi(x)  \right| \\
\leq & c_1  \| \phi \|_{\xC^2} M \Delta t \, \e^{-2} \, c_3 M \Delta t \: \e^{-4} + \| \phi \|_\infty c_2 (c_1 M \Delta t \: \e^{-4} )^2 + c_1^2 M^2 \| \phi \|_{\xC^2} \Delta t^2 \: \e^{-4} \\
\leq & c_3 \|  \phi \|_{\xC^2} M^2 \Delta t^2 \: \e^{-8}
\end{align*}
and integrating the previous inequality together with \eqref{eq:weightedFVproof} leads to \eqref{appbrakke}.

\noindent {\bf Step 4:}
By definition \eqref{dfnhepsilon}, $h_\e = \Phi_\e \ast \tilde{h}_\e$ and thus, for all $S \in \G$, $\displaystyle \mdiv_S(h_{\e}) = S(\nabla\Phi_{\e})\ast\tilde{h}_{\e}$.
Then, by definition of $\delta V$ and \eqref{convvar}, we obtain \eqref{massdecrease1}:
\begin{align*}
\delta V(h_{\e})& =\int_{\R^n\times \G} \mdiv_S(h_{\e}(x)) \: dV(x,S) = \int_{\R^n\times \G} \int_{\R^n} S(\nabla\Phi_{\e}(y-x))\cdot\tilde{h}_{\e}(y) \: dy \:  dV(x,S)
\\& = \int_{\R^n} \int_{\R^n\times \G}  S(\nabla\Phi_{\e}(y-x)) \: dV(x,S)\cdot \tilde{h}_{\e}(y) \: dy 
\\& =\int_{\R^n} (\Phi_{\e}\ast \delta V)(y) \cdot \tilde{h}_{\e}(y) \: dy  = - \int_{\R^n}  \frac{|(\Phi_{\e}\ast \delta V)(y)|^2}{(\Phi_{\e}\ast\V)(y)+\e} \: dy \leq 0.
\end{align*}
We are left with the proof of \eqref{massdecrease6} and we assume $c_3 \del M^3 \leq \epsilon^8$ then $\Delta t$ in particular satisfies \eqref{smallstep0}, assumption under which the map $f$ is a diffeomorphism of $\R^n$ and \eqref{eq:JSfboundsOrder2} holds. Consequently, applying Definition~\ref{dfn:varifoldpush-forward} of push-forward varifold and using \eqref{eq:diffeoProof1}, \eqref{eq:JSfboundsOrder2} and \eqref{massdecrease1}, we obtain
\begin{align*}
 \|f_{\#}V \|(\R^n) & = \int_{\R^n \times \G} J_S f(x) \: d V(x,S) \\
 & = \int_{\R^n \times \G} 1 +  \Delta t \: \mdiv_S(h_{\e}(x)) + \left( J_Sf(x) -1 - \Delta t \: \mdiv_S(h_{\e}(x)) \right) \: dV(x,S) \\
 & \leq \|V \|(\R^n) + \Delta t \: \delta V(h_\e) +  c_2M(\Delta t \|Dh_{\e} \|_{\infty})^2  \\
 & \leq \|V \|(\R^n)+ c_1^2c_2 M^3 \Delta t^2 \e^{-8}\\
 & \leq \|V \|(\R^n) + \Delta t \:
\end{align*} 
hence concluding the proof of \eqref{massdecrease6}.
\end{proof}
\begin{lemma}\label{pre_lem:diff}
Let $f\in \xC(\R^n,\R^n)$ such that $f := {\rm id} + \del \, h, \, \del >0, \, h\in \xC^1(\R^n,\R^n)$, assume that 
\begin{equation*}
\max \lbrace \del \|h\|_{\infty}, \del \|Dh\|_{\infty}, \|Jf -1 \|_{\infty} \rbrace < 1.  
\end{equation*}
Then, $f$ is a diffeomorphism of $\R^n$.
\end{lemma}
\begin{proof}
For any $x\in \R^n$, we have $Jf(x) \neq 0$, therefore $f$ is a local diffeomorphism by the inverse mapping theorem. We now prove that $f$ is injective, indeed, for any $x,y \in \R$ one has 
\begin{equation*}
 |f(x)-f(y)| = |x-y-\del \left( h(x)-h(y) \right)| \geq \big| |x-y| - \del |h(x)-h(y)| \big| \geq |x-y| \big| 1- \del \, \lip(h)\big| >0,
\end{equation*}
this proves the injectivity of $f$.\\
Up to now, we have checked that $f$ is injective and a local diffeomorphism at every point therefore $f$ is a global diffeomorphism from $\R^n$ onto $f(\R^n)$; as $f(\R^n)$ is open, it remains to show that $f(\R^n)$ is closed. We have by assumption that $\|f-{\rm id}\|_{\infty}=\del \|h\|_{\infty} < 1$, this implies that $f$ is proper, by \cite{palais} it is closed, therefore $f(\R^n)$ is closed in $\R^n$, recalling that it was open we have $f(\R^n)=\R^n$ and $f$ is a diffeomorphism of $\R^n$.
\end{proof}

Given $M \geq 1$ and a $d$-varifold $V_0$ satisfying $\| V_0 \|(\R^n)  \leq M$, Proposition~\ref{mcfmotion1} gives the condition $c_3 \Delta t < M^{-3} \epsilon^8$ allowing to define $V_1 = {f_0}_\# V_0$ with $f_0 = f_{\epsilon, V_0} = {\rm id} + \Delta t \: h_\epsilon(\cdot, V_0)$.
We would like to iterate on several time steps and thus push the varifold $V_1$ by the map $f_1 = f_{\e, V_1}={\rm id}+\Delta t \:  h_{\e}(\cdot, V_1)$. However, note that $\|{f_0}_{\#}V_0\|(\R^n) \leq M + \Delta t$ and not necessarily $\|{f_0}_{\#}V_0\|(\R^n) \leq M$:
the choice of $\Delta t$ is no longer suitable. To rule out this issue, we can initially choose $\Delta t$ satisfying
\begin{equation*}\label{defc_3}
c_3  \Delta t < (M+1)^{-3}\e^{8} ,
\end{equation*}
and \eqref{massdecrease6} 
thus ensuring
$\|V_1\|(\R^n)= \|{f_0}_{\#}V_0\|(\R^n) \leq \|V_0 \|(\R^n) + \Delta t \leq M+1$,
and we can iterate the process as long as the mass remains less than $M+1$, thus at least  $ \lfloor 1/ \Delta t \rfloor $ times when considering uniform time discretizations of $[0,1]$. Considering a possibly non uniform time discretization $(\Delta t_i)_{i=1 \ldots m} \in (0,1)$ of $[0,a]$ for $a \leq 1$: $\sum\limits_{i=1}^m \Delta t_i = a \leq 1$, one can iterate the process $m$ times with $\Delta t_i$ being the time step at step $i$, this justifies the following definition.
\begin{dfn}[Time-discrete approximate mean curvature flow]
\label{damcf} 
Let $M \geq 1$, $\e\in(0,1)$ and $a\in(0,1]$. Consider a subdivision
$\cT=\lbrace t_i \rbrace_{i=0}^m$ of $[0,a]$ and assume
\begin{equation}\label{smallstep}
c_3 \delta(\cT)  \leq (M+1)^{-3}\e^{8} 
\end{equation}
where $\Delta t_i = t_i - t_{i-1}$ for $i = 1, \ldots, m$ and $\delta(\cT) = \max_{1 \leq i \leq m} \Delta t_i$.\\
Let $V_0 \in V_d(\R^n)$ satisfy $\|V_0\|(\R^n)\leq M$. Define $\left( V_{\e,\cT} (t_i) \right)_{i = 0 \ldots m}$ by $V_{\e,\cT}(0):= V_0$ ($t_0 = 0$) and, for $i = 1, \ldots, m$,
\begin{equation*}
 V_{\e,\cT}(t_{i}) := {f_i}_\# V_{\e,\cT}(t_{i-1}) \quad \text{with} \quad 
 f_i = {\rm id}+ \del_i \,h_{\e}(\cdot,V_{\e,\cT}(t_{i-1})) \: .
\end{equation*}
We then define the family $\left( V_{\e,\cT}(t) \right)_{t \in [0,a]}$ by linear interpolation between the points of the subdivision, and we call it a \emph{time-discrete approximate mean curvature flow}: 
\begin{equation*}
 V_{\e,\cT}(t) := \left[ {\rm id}+(t-t_i)h_{\e}(\cdot,V_{\e,\cT}(t_i))\right]_{\#}V_{\e,\cT}(t_{i}) \quad \text{if} \quad t \in[t_{i},t_{i+1}].
\end{equation*}
\end{dfn}
\begin{remk}\label{rem:uniformbound}
We note that under the assumptions of Definition~\ref{damcf} (and using the same notations),
we have 
\begin{equation}\label{uniformbound}
  \|V_{\e,\cT}(t)\|(\R^n) \leq  M+1, \quad \forall t \in [0,a] \: ,
\end{equation}
and we will use \eqref{uniformbound} extensively throughout the chapter.\\
Moreover, if we assume that
there exists $R_0 > 0$ such that $\supp V_0 \subset B_{R_0}(0) \times \G$, then 
\begin{equation*}
\forall t \in [0,a], \quad \supp V_{\e, \cT} \subset B_{R_0 + c_1 (M+1) \e^{-2}}(0) \times \G \: .
\end{equation*} 
Indeed, thanks to Proposition \ref{hepsilonbound} and \eqref{uniformbound}, for $t \in [0,a] \subset [0,1]$, 
\begin{equation*}
 \|h_{\e}(\cdot,V_{\e,\cT}(t))\|_{\infty} \leq c_1 (M+1)\e^{-2}
\end{equation*}
and therefore, 
$\supp V_{\e,\cT}(t) \subset B \left(0, R_0 + c_1 \: t \: (M+1) \e^{-2} \right) \times \G$. Such compactness property will be used when letting $\Delta t$ go to $0$ (for a fixed $\e$) to define a limit flow in Section~\ref{def_time_step_to_0}.
\end{remk}

\begin{remk}[Piecewise constant flow]
\label{remk:PWflow}
Note that in Definition~\ref{damcf}, we first define $V_{\e,\cT}(t_i)$ at the points $t_i$ of the subdivision $\cT$ and we then define $V_{\e,\cT}(t)$ for $t \in [t_i,t_{i+1}]$ by a linear interpolation between $t_i$ and $t_{i+1}$. It is possible to consider an alternative definition of the flow between $t_i$ and $t_{i+1}$, simply taking the following piecewise constant extension: for $i \in \{0,1,\ldots,m-1\}$,
\begin{equation*}
 V_{\e,\cT}^{pc}(t) := V_{\e,\cT}(t_{i}) \quad \text{if} \quad t \in(t_i,t_{i+1}) \: .
\end{equation*}
As we will see in Proposition~\ref{coincidence}, both $V_{\e,\cT}$ and $V_{\e, \cT}^{pc}$ lead to the same limit flow $V_\e$ when the size of the subdivision tends to zero. We consequently restrict our study to only one of the two flows and we choose to investigate $V_{\e, \cT}$ introduced in Definition~\ref{damcf}.
\end{remk} 
Hereafter, $M \geq 1$ and $\epsilon \in (0,1)$ are fixed, all subdivisions we consider satisfy \eqref{smallstep} and we define time-discrete approximate mean curvature flow starting from a varifold of mass less than $M$.
\subsection{Stability with respect to initial data}
\label{subsec:stabilityInitialData}
When investigating a discrete scheme for a flow, the stability arises as a crucial issue. More precisely, we consider in Proposition~\ref{stability} two time-discrete approximate mean curvature flows $(V(t))_t$, $(W(t))_t$ starting from $V_0$, $W_0$, respectively, and we prove that the stability holds in terms of bounded Lipschitz distance: $\Delta( V(t) , W(t) ) \leq \exp (\Lambda t) \Delta (V_0, W_0)$, where $\Lambda \sim \epsilon^{-n-7}$. Up to a constant, $\Lambda$ is an upper bound of the Lipschitz constant of $V \mapsto h_\epsilon (\cdot, V)$  with respect to the $\xC^1$--norm, as established in Lemma~\ref{hepsilonstability}. In Remark~\ref{remk:ODEnum}, we draw a parallel with the classical time discretization of ODEs showing that $\Lambda$ is the expected constant in our setting.

\begin{lemma}
\label{hepsilonstability}
Let $\e \in (0,1)$ and $M \geq 1$. Let $V$ and $W$ be two varifolds of $ V_d(\R^n)$ satisfying $\V(\R^n) \leq M$, $\W(\R^n)\leq M$. There exists $c_4>0$ only depending on $n$ such that
\begin{equation*}
\big\|  h_{\e}(\cdot,V) -   h_{\e}(\cdot,W)   \big\|_{\infty} \leq c_4 M \e^{-n-5} \Delta(V,W) \text{ and } \big\|  D h_{\e}(\cdot,V) -  D h_{\e}(\cdot,W)  \big\|_{\infty} \leq c_4 M \e^{-n-7} \Delta(V,W). 
\end{equation*}
\end{lemma}
\begin{proof}
As previously, we may increase $c_3 > 0$ throughout the proof to meet the requirements, provided that we respect the exclusive dependency in $n$.
Let $\e \in (0,1)$ and $M \geq 1$. Let $V$ and $W$ be two varifolds of $ V_d(\R^n)$ satisfying $\V(\R^n) \leq M$, $\W(\R^n)\leq M$.
We first show that
\begin{equation} \label{eq:hepsLip}
\begin{array}{ll}
 \left\| \Phi_\e \ast \V - \Phi_\e \ast \W \right\|_\infty  \leq (c(2\pi)^{-\frac{n}{2}}+c_0)\e^{-n-2} \Delta(V,W), \, \text{and} \\
\left\| \Phi_\e \ast \delta V - \Phi_\e \ast \delta W \right\|_\infty \leq 2(c(2\pi)^{-\frac{n}{2}}+c_0)\e^{-n-4}  \Delta (V,W)\:.
\end{array}
\end{equation}
We have for any $x\in \R^n$ by \eqref{lipker} and the definition \eqref{defkernel} of $\Phi_{\e}$ 
\begin{equation*}
\begin{split}
 \Big| \Phi_\e \ast \V(x) - \Phi_\e \ast \W(x) \Big| 
 &= \Big|\int_{\R^n} \Phi_{\e}(x-y) d \V(y) -\int_{\R^n} \Phi_{\e}(x-y) d \W(y)\Big|
 \\& 
= \Big| \V\left( \Phi_{\e}(\cdot-x) \right) -  \W\left( \Phi_{\e}(\cdot-x) \right)  \Big|
 \\& \leq \max\lbrace  \| \Phi_{\e} \|_{\infty} , \lip(\Phi_{\e}) \rbrace \Delta(\V, \W)
\\& \leq \max \lbrace  c(2\pi)^{-\frac{n}{2}} \e^{-n} , (c(2\pi)^{-\frac{n}{2}} + c_0) \e^{-n-2} \rbrace \Delta(\V, \W) 
\\& \leq (c(2\pi)^{-\frac{n}{2}} + c_0) \e^{-n-2}  \Delta(V, W) \: ,
\end{split}
\end{equation*}
since $\Delta(\V, \W) \leq \Delta(V, W)$, this gives the first estimate of \eqref{eq:hepsLip}. For the second estimate we first recall that for $x \in \R^n$, $\displaystyle \Phi_\e \ast \delta V(x) = \int_{\R^n \times \G} S \nabla \Phi_\e (x -y) \: dV(y,S)$ and
we thus compute the Lipschitz constant of the map $\Theta : (y,S) \mapsto S (\nabla\Phi_{\e}(y))$ (the map $y \mapsto x-y$ being an isometry), we have for $(y,S), (z,T) \in \R^n \times \G$, using $\|S \| =1$
\begin{equation*}
\begin{split}
\left| \Theta (y,S) - \Theta (t,y) \right| & = \left| S(\nabla\Phi_{\e}(y)) - T(\nabla\Phi_{\e}(z)) \right| \\
& \leq \|S\| \left| \nabla\Phi_{\e}(y) - \nabla\Phi_{\e}(z) \right| + \| S- T\| \left| \nabla \Phi_{\e}(z)\right|
\\& \leq \lip(\nabla\Phi_{\e}) |y-z| + \| \nabla \Phi_{\e} \|_{\infty} \|S-T\|
\\& \leq 2(c(2\pi)^{-\frac{n}{2}}+c_0)\e^{-n-4} \quad \text{thanks to Lemma~\ref{kerprop}.}
\end{split}\end{equation*}
Therefore $\lip(\Theta) \leq 2(c(2\pi)^{-\frac{n}{2}}+c_0)\e^{-n-4}$, also from \eqref{kerprop1} we have $\|\Theta\|_{\infty}\leq (c(2\pi)^{-\frac{n}{2}}+c_0)\e^{-n-2}$. We can now carry on with the proof the the second inequality of \eqref{eq:hepsLip}: for $x \in \R^n$,
\begin{equation*}
\begin{split}
\Big|\Phi_\e \ast \delta V(x) - \Phi_\e \ast \delta W(x) \Big|
&= \Big|\int_{\R^n} S(\nabla\Phi_{\e})(x-y) d V(y,S) -\int_{\R^n} S(\nabla\Phi_{\e})(x-y) d W(y,S)\Big|
\\& \leq \max\lbrace  \| \Theta \|_{\infty} , \lip(\Theta) \rbrace \Delta(V, W)
\\& \leq 2(c(2\pi)^{-\frac{n}{2}} + c_0) \e^{-n-4}  \Delta(V, W) \: ,
\end{split}
\end{equation*}
which gives the desired result.
%
From \eqref{eq:hepsLip} and \eqref{hepsilonbound1}, we have for $x \in \R^n$,
\begin{align}\label{eq:tildeHepsLip}
\Big|   \tilde{h}_{\e}(x,V) -  &  \tilde{h}_{\e}(x,W)   \Big| =  \Big| \frac{\Phi_{\e}\ast \delta V(x)}{\Phi_{\e}\ast\V(x)+\e} - \frac{\Phi_{\e}\ast \delta W(x)}{\Phi_{\e}\ast\W(x)+\e} \Big| \nonumber
\\& 
\leq  \frac{| \Phi_{\e}\ast \delta V(x) - \Phi_{\e}\ast \delta W(x)|}{\Phi_{\e}\ast\V(x)+\e} + \left| \frac{\Phi_{\e}\ast \delta W(x)}{\Phi_{\e}\ast\W(x)+\e}\right| \frac{\big|\Phi_{\e}\ast\V(x)-\Phi_{\e}\ast\W(x)\big|}{\Phi_{\e}\ast\V(x)+\e} \nonumber
\\
& \leq \frac{1}{\e} \| \Phi_{\e}\ast \delta V - \Phi_{\e}\ast \delta W \|_\infty + \| \tilde{h}_\e (\cdot,W) \|_\infty \frac{1}{\e} \left\| \Phi_{\e}\ast\V -\Phi_{\e}\ast\W  \right\|_\infty \nonumber
\\
& \leq  \frac{1}{\e} \left(  2(c(2\pi)^{-\frac{n}{2}}+c_0)\e^{-n-4} + c_1 M \e^{-2} (c(2\pi)^{-\frac{n}{2}}+c_0)\e^{-n-2} \right) \Delta(V,W) \nonumber \\
& \leq  (c(2\pi)^{-\frac{n}{2}} + c_0 )  \left( 2 + c_1 M  \right) \e^{-n-5} \Delta(V,W)   \\
& \leq c_4 M  \e^{-n-5} \Delta(V,W) \: . \nonumber
\end{align}
We recall that $h_{\e}=\Phi_{\e} \ast \tilde{h}_{\e}$ and we obtain thanks to \eqref{eq:PhiEpsL1} and \eqref{eq:tildeHepsLip}:
\begin{align*}
\|   h_{\e}(\cdot,V) -  h_{\e}(\cdot,W)   \|_{\infty} 
&\leq \|\Phi_{\e}\|_{L^1} \|\tilde{h}_{\e}(\cdot,V) -  \tilde{h}_{\e}(\cdot,W) \|_{\infty} 
\leq  c_4 M  \e^{-n-5} \Delta(V,W)  \: .
\end{align*}
Similarly $Dh_{\e}=\nabla\Phi_{\e}\ast\tilde{h}_{\e}$ and using \eqref{kerprop7} and \eqref{eq:tildeHepsLip}, we obtain
\begin{align*}
\|Dh_{\e}(\cdot,V)-Dh_{\e}(\cdot,W) \|_{\infty} & \leq \|\nabla\Phi_{\e} \|_{L^1} \| \tilde{h}_{\e}(\cdot,V)-\tilde{h}_{\e}(\cdot,W) \|_{\infty}
 \\
& \leq \e^{-2}(1+c_0\omega_n) \: (c(2\pi)^{-n/2} +c_0)(2+c_1 M) \: \Delta(V,W) \e^{-n-5}
\\
& \leq c_4 M \e^{-n-7} \Delta(V,W) ,
\end{align*}
hence concluding the proof.
\end{proof}
In Proposition~\ref{damcfstability}, we investigate the evolution of the bounded Lipschitz distance between two varifolds $V$ and $W$ through one step of the time-discrete approximate flow introduced in Definition~\ref{damcf}. The proof relies on Lemma~\ref{StildetoS}, Lemma \ref{hepsilonstability} and on careful estimates of the Lipschitz constant of the map $(x,S,V) \in \R^n \times \G \times V_d(\R^n) \mapsto J_Sf_{\e,V}(x)$.
\begin{prop}\label{damcfstability}
Let $\e, \del \in (0,1)$ and $M \geq 1$. Let $V, \, W \in V_d(\R^n)$ satisfy $\| V \|(\R^n) \leq M$, $\| W \|(\R^n) \leq M$. Let $g \in C^1 (\R^n,\R^n)$ be such that $\|Dg- I_n \|_{\infty} \leq c_1 M \del \e^{-4}$ and recall the notation $f_{\e,V} = I_n + \Delta t \: h_{\e}(\cdot, V)$.
There exists $c_5 > 0$ depending only on $n$ such that if $\del$ satisfies
\begin{equation}\label{smallstep2}
 c_5 M^3 \del \leq \e^{8} \: ,
\end{equation}
then, for any $(x,S), (t,y) \in \R^n \times \G$ 
\begin{align}
& \big| J_S f_{\e,V}(x) - J_T g(y) \big| \leq c_5  M \del \left( \e^{-4}\| S - T\| +  \e^{-6}  |x-y| \right)  + c_5 \| Df_{\e,V} - Dg \|_{\infty} ,\label{damcfstability2} \\
& \big| J_S f_{\e,V}(x) - J_T f_{\e,W}(y) \big| \leq c_5  M \del  \left( \e^{-4} \| S-T \| + \e^{-6} |x-y| + \e^{-n-7} \Delta(V,W) \right). \label{damcfstability3}
\end{align}
and
\begin{align}
& \Delta((f_{\e,V})_{\#}V, g_{\#}V) \leq M \left( c_5 \| f_{\e,V}- g \|_{C^1}  + \| J_{\cdot}f_{\e,V}- J_{\cdot}g\|_{\infty} \right), \label{damcfstability5} \\
& \Delta \left( (f_{\e,V})_{\#} V , (f_{\e,W})_{\#} W \right) \leq (1+c_5  M^2 \del \e^{-n-7})\Delta(V,W). \label{damcfstability6}
\end{align}
\end{prop}
\begin{proof}
Note that throughout the proof, we adapt i.e. increase $c_5 > 0$ with the constraint that $c_5$ only depends on $n$ (like other $c_i$) whenever needed.
As previously, $\e$ is fixed throughout the proof and we can write $f_{V}$ (resp. $f_{W}$, $h$) instead of $f_{\e,V}$ (resp. $f_{\e,W}$, $h_\e$). We fix $(x,S), (t,y) \in \R^n \times \G$ and thanks to Lemma~\ref{StildetoS}, we choose $\tilde{S}, \tilde{T} \in \cM_{d,n}$ such that, as orthogonal projectors, $S = \tilde{S}^t \circ \tilde{S}$,  $T = \tilde{T}^t \circ \tilde{T}$, where $\circ$ denotes the matrix product for clarity. Moreover, 
\begin{equation*}
 \tilde{S} \circ \tilde{S}^t = \tilde{T} \circ \tilde{T}^t = I_d \quad \text{and} \quad \| \tilde{S} - \tilde{T} \| \leq 2 \| S - T \| \: .
\end{equation*}
We introduce some additional notations: we write $G = Dg - I_n$, $F = Df_V - I_n = \del Dh( , V)$, and we recall that, by hypothesis on $g$, \eqref{smallstep2} and Proposition~\ref{hepsilonbound}, one has (noting that $\del M \e^{-4} \leq \del M^3 \e^{-8}$ and taking $c_5$ large enough),
\begin{equation}\label{G-bound}
 \| G(y) \| \leq c_1 M \del \e^{-4} \leq \frac{c_1}{c_5} \leq 1 \quad \text{and} \quad \|F(x)\| \leq  c_1 M \del \e^{-4} \leq 1 \: .
\end{equation}
We also use the notations
\begin{equation*}
P = \tilde{S} \circ Df_V(x)^t \circ Df_V(x) \circ \tilde{S}^t \quad \text{and}  \quad N = \tilde{T}\circ Dg(y)^t \circ Dg(y) \circ \tilde{T}^t \:.
\end{equation*}
and we recall that $J_S f_ V(x) = \det (P)^\frac{1}{2}$ and similarly $J_T g(y) = \det (N)^\frac{1}{2}$.

\smallskip
\noindent {\bf Step 1:} We first prove that
\begin{equation}\label{damcfstability1}
\| P - N \| \leq 3 \left( 4 c_1 M \del \e^{-4}\| S - T \| +   c_1 M \del \e^{-6}  |x-y| +  \| Df_{V} - Dg \|_{\infty} \right).
\end{equation}
Let us decompose $P$ and $N$ as follows, using $Df_V(x) = I_n + F(x)$,
\begin{equation} \label{eqDecP}
\begin{split}
P  = \tilde{S} \circ (I_n + F(x)) \circ (I_n + F(x))^t \circ \tilde{S}^t 
 = I_d +  \tilde{S} \circ \left( F(x)^t + F(x) \right) \circ \tilde{S}^t +  \tilde{S}\circ F(x)^t \circ F(x) \circ \tilde{S}^t 
 \end{split}
\end{equation}
and similarly
\begin{equation} \label{eqDecN}
N = I_d +  \tilde{T} \circ \left( G^t(y) + G(y) \right) \circ \tilde{T}^t  +  \tilde{T}\circ G^t(y) \circ G(y) \circ \tilde{T}^t \: .
\end{equation}
On one hand, a simple computation (see Lemma~\ref{lem:intro_trivial_lemma}) together with \eqref{G-bound} give
\begin{align}
\label{af-ag}
   \| \tilde{S} & \circ \left( F(x)^t + F(x) \right) \circ \tilde{S}^t   - \tilde{T} \circ \left( G(x)^t + G(x) \right) \circ \tilde{T}^t \| \nonumber \\
   & \leq  2 \left( \| F(x)^t + F(x)\| + \|G(y)^t + G(y)\| \right) \| \tilde{S} - \tilde{T} \| + 2 \| F(x)^t + F(x) - G(y)^t - G(y) \| \nonumber \\
   & \leq  2 \left( \| F(x)\| + \|G(y)\| \right) \| \tilde{S} - \tilde{T} \| + 2 \| F(x) - G(y) \| 
\end{align}
Similarly using Lemma~\ref{lem:intro_trivial_lemma} and \eqref{G-bound}, we obtain 
\begin{align}\label{bf-bg}
 \| \tilde{S} \circ  F(x)^t  \circ F(x) & \circ \tilde{S}^t  - \tilde{T} \circ  G(x)^t \circ G(x)  \circ \tilde{T}^t \| \nonumber \\
& \leq  \left( \| F(x)^t \circ F(x) \| + \| G(y)^t \circ G(y) \| \right) \| \tilde{S} - \tilde{T} \| +  \| F(x)^t \circ F(x) - G(x)^t \circ G(x) \| \nonumber \\
& \leq \left( \| F(x)\| + \| G(y) \| \right) \| \tilde{S} - \tilde{T} \|  + 2 \| F(x) - G(y) \| 
\end{align}
Consequently, from \eqref{eqDecP} to \eqref{bf-bg}, we infer
\begin{equation}
    \| P - N \| \leq 3 \left( \| F(x)\| + \| G(y) \| \right) \| \tilde{S} - \tilde{T} \|  + 3 \| F(x) - G(y) \| 
\end{equation}
so that we can conclude the proof of \eqref{damcfstability1} (Step 1) recalling that $F(x) = \del Dh_{\e}(x,V)$ and applying \eqref{hepsilonbound2} and \eqref{hepsilonbound3}:
\begin{equation} \label{eq:Fx-Gy}
\begin{split}
 & \| F(x) - G(y) \|  \leq \| F(x) - F(y) \| + \| F(y) - G(y) \| \leq c_1 M \del \e^{-6}  |x-y| + \| Df_V - Dg \|_{\infty} , \\
 & \| F(x) \| + \| G(y) \|  \leq 2 c_1 M \del \e^{-4} \| \tilde{S} - \tilde{T} \| \leq  4 c_1 M \del \e^{-4} \| S - T \| \: .
 \end{split}
\end{equation}

\smallskip
\noindent {\bf Step 2:} We now prove \eqref{damcfstability2} and \eqref{damcfstability3}.

\noindent Let us show that
\begin{equation} \label{eq:boundinvPinvN}
 \| P - I_d \| \leq \frac{1}{4}, \quad \| N - I_d \| \leq \frac{1}{4} \quad \text{and } P, \, N \text{ are invertible with } \|P^{-1} \| \leq 2 \text{ and } \|N^{-1} \| \leq 2 \: . 
\end{equation}
To this end, we apply Lemma~\ref{detprop} with $L=\tilde{S}$ and $R = F(x) = \del \, Dh(x,V)$: using \eqref{hepsilonbound2} and \eqref{smallstep2},
\begin{equation*}
|F(x)|_\infty = \Delta t \, \left| Dh_{\e}(x,V) \right|_\infty \leq   c_1 \del M \e^{-4} \leq \frac{c_1}{c_5} \leq 1 \quad \text{and} \quad c_2 |F(x)|_\infty \leq \frac{c_1 c_2}{c_5} \leq 1 \: ,
\end{equation*}
which allows to apply \eqref{detprop3} so that
\begin{equation} \label{eq:proofMId}
\begin{split}
\| P- I_d \| \leq d \left| P - I_d \right|_\infty 
& \leq d\Big| \tilde{S} \circ \left( F(x)^t + F(x) \right) \circ \tilde{S}^t +  \tilde{S}\circ F(x)^t \circ F(x) \circ \tilde{S}^t  \Big|_{\infty}
\\& \leq n c_2  \left| F(x) \right|_\infty \leq nc_2\frac{c_1}{c_5}  < \frac14,
\end{split}
\end{equation}
where $\circ$ denotes the matrix product for clarity.
As $c_2 |F(x)|_\infty  \leq 1$, we can also apply \eqref{detprop4} to conclude that $\left| \det (P)^\frac{1}{2} - 1\right| \leq c_2 |F(x)|_\infty \leq \frac{c_1 c_2}{c_5} < 1 $ and thus $P$ is invertible. Furthermore, using \eqref{eq:proofMId}, we have
\begin{equation}\label{inversebound}
\| P^{-1} \| \leq  \| P^{-1} - I_d \| + \| I_d \| \leq \| P^{-1}\|  \| I_d- P \| + 1 \leq \frac12 \| P^{-1} \| + 1 \quad \Rightarrow \quad \| P^{-1} \| \leq 2 \: .
\end{equation}
We similarly proceed with $G$, we recall \eqref{G-bound} 
$|G(y)|_{\infty} \leq \|G(y)\| \leq 2c_1 M \del \e^{-4} \leq \frac{2c_1}{c_5} \leq 1$ and we can apply \eqref{detprop3} to obtain
\begin{equation}\label{damcfstabilityprf2}
\begin{split}
\left\| N - I_d \right\| & \leq d \left| N - I_d \right|_\infty = d \left| \tilde{T} \circ (  G(y)^t +  G(y)) \circ \tilde{T}^t+   \tilde{T}  \circ  G(y)^t  \circ G(y) \circ \tilde{T}^t \right|_\infty \\& \leq n c_2 \left| G(y) \right|_\infty \leq nc_2 \frac{2c_1}{c_5} \leq \frac{1}{4} \: .
\end{split}
\end{equation}
Note that similarly to \eqref{inversebound}, we also have that $N$ is invertible and $\| N^{-1} \| \leq 2 $.

\noindent We recall that $\frac14 \leq | \det(P) | = J_Sf_V(x)^2 \leq 4$ thanks to \eqref{eq:JSfbounds} (which applies since \eqref{smallstep0} holds for $c_5$ large enough).
We can now show that $P$ and $N$ satisfy the condition given in \eqref{claimcondition}. Indeed, from \eqref{eq:proofMId}, \eqref{inversebound} and \eqref{damcfstabilityprf2}, we have
\begin{equation*}
 \| P^{-1} \| \| P-N \|  \leq 2 \| (P-I_d) -  (N-I_d) \| \leq 2 \left( \frac{1}{4} + \frac{1}{4} \right) \leq 1
\end{equation*}
and thus applying \eqref{claimcondition} leads to
\begin{equation}\label{eq:detStability0}
 \left| \det(P) - \det(N) \right| \leq  c_2 \left| \det(P) \right| \: \| P^{-1} \| \: \| P - N \| \leq 8 c_2 \: \| P- N \| \: .
\end{equation}
We note that for $a\geq \frac12$, $b\geq0$ 
\begin{equation*}
 |a-b| = \frac{|a^2-b^2|}{a+b} \leq 2|a^2-b^2| \: ,
\end{equation*}
and we apply it with $a = \det(P)^{\frac12}=J_Sf_V(x)\geq \frac12$ and $b = \det(N)^{\frac12} = J_T g(y)$ so that, using \eqref{eq:detStability0}, we obtain
\begin{equation}\label{eq:detStability}
 \left|  J_S f_V(x) - J_T g(y) \right| = \left|  \det(P)^{\frac12} - \det(N)^{\frac12} \right|\leq  2\left| \det(P) - \det(N) \right| \leq 16 c_2 \| P-N \|.
\end{equation}
We conclude the proof of \eqref{damcfstability2} gathering \eqref{eq:detStability} and \eqref{damcfstability1}. Then, applying \eqref{damcfstability2} with $g = f_W$ allows to prove \eqref{damcfstability3}. Indeed, $f_W \in C^1(\R^n, \R^n)$ and from \eqref{hepsilonbound2} we have $\|Df_W  - I_n \|_{\infty} = \del \| Dh(\cdot, W) \|_\infty \leq c_1 M \del \e^{-4}$ so that 
\begin{equation*}\label{damcfstability3prf1}
 \big| J_S f_{V}(x) - J_T f_W(y) \big| \leq c_5 M \del \left( \e^{-4}\| S-T \| +   \e^{-6}  |x-y| \right) + c_5 \| Df_{V} - Df_W \|_{\infty} 
\end{equation*}
and by Lemma~\ref{hepsilonstability}:
\begin{equation*}\label{damcfstability3prf2}
\begin{split}
 \| Df_{V} - Df_W \|_{\infty} \leq \del \| Dh (\cdot,V) -  Dh (\cdot,W) \|_{\infty} \leq c_4 M \del \e^{-n-7} \Delta(V,W) 
 \end{split}
\end{equation*}
hence concluding the proof of \eqref{damcfstability3} and Step 2.

\smallskip
\noindent {\bf Step 3:} We now study the Lipschitz constant of the map $(x,S,f) \mapsto Df(x)(S) \in \G$, which will be crucial in proving the remaining estimates of the proposition. We namely prove 
\begin{equation}\label{damcfstability4}
 \| Df_{V}(x)(S) - Dg(y)(T) \| \leq ( 1 + c_5 M\del \e^{-4} ) \| S-T \| + c_5 M\del \e^{-6} |x-y| + c_5 \| Df_{V} - Dg \|_{\infty} \: . 
\end{equation}
We recall (see \eqref{eq:DfS}) that we can write
\begin{multline*}
 Df_V(x)(S) = Y(Y^tY)^{-1}Y^t = Y \circ P^{-1} \circ Y^t \quad \text{and} \quad  Dg(y)(T) = Z(Z^tZ)^{-1}Z^t = Z \circ N^{-1} \circ Z^t 
 \\ \quad \text{with} \quad Y = Df(x)\circ \tilde{S}^t = \tilde{S}^t + F(x) \circ \tilde{S}^t \quad \text{and} \quad  Z = Dg(y)\circ \tilde{T}^t = \tilde{T}^t + G(y) \circ \tilde{T}^t \: ,
\end{multline*}
so that we can further decompose $Df_V(x)(S)$ and $Dg(y)(T)$ as follows.
By the formulas $\tilde{S}^t \circ \tilde{S} = S$ and $\tilde{T}^t \circ \tilde{T} = T$, if we set $\tilde{P} = P^{-1} - I_d$ and $\tilde{N} = N^{-1} - I_d $ we obtain
\begin{equation*}
\begin{split}
 & Df_V(x)(S) = \left( \tilde{S}^t + F(x)\circ \tilde{S}^t \right) \circ P^{-1} \circ \left( \tilde{S} + \tilde{S}\circ F(x)^t \right) \\
 & \quad = \underbrace{\tilde{S}^t \circ P^{-1} \circ \tilde{S} }_{C_f}
 + \underbrace{F(x)\circ \tilde{S}^t \circ P^{-1} \circ \tilde{S} +  \tilde{S}^t \circ P^{-1}\circ\tilde{S} \circ  F(x)^t}_{R_f}+ \underbrace{F(x) \circ \tilde{S}^t \circ P^{-1} \circ  \tilde{S}\circ F(x)^t }_{E_f} \\
 & \quad  = C_f + R_f + E_f
 \end{split}
\end{equation*}
and similarly
\begin{equation*}
\begin{split}
 & Dg(y)(T) = \left( \tilde{T}^t + G(y) \circ \tilde{T}^t \right) \circ N^{-1} \circ \left( \tilde{T} + \tilde{T}\circ G(y)^t \right)\\
 & \quad  = \underbrace{\tilde{T}^t \circ N^{-1} \circ \tilde{T}}_{C_g} + \underbrace{ G(y) \circ \tilde{T}^t \circ N^{-1} \circ \tilde{T} + \tilde{T}^t \circ N^{-1} \circ \tilde{T}\circ G(y)^t}_{R_g} +  \underbrace{G(y) \circ \tilde{T}^t \circ N^{-1} \circ \tilde{T} \circ G(y)^t}_{E_g} \\
 & \quad = C_g + R_g + E_g \: ,
 \end{split}
\end{equation*}
so that
\begin{equation}\label{projector-diff}
 \| Df_V(x)(S) - Dg(y)(T)  \| \leq  \| C_f - C_g \|  + \|  D_f - D_g \|  + \| E_f - E_g \| \: ,
\end{equation}
We first prove that 
\begin{equation}\label{cf-cg}
 \| C_f - C_g \| \leq (1 + c_5 M \del \e^{-4})\| S -T \| + c_5 M \del  \e^{-6} |x-y| + c_5 \| Df_V - Dg \|_{\infty}  \: .
\end{equation}
Since $\tilde{S}^t \circ \tilde{S} = S$ and $\tilde{T}^t \circ \tilde{T} = T$, we have $C_f - C_g = S - T + \tilde{S}^t \circ (P^{-1} - I_d) \circ \tilde{S} - \tilde{T}^t \circ (N^{-1} - I_d) \circ \tilde{T}$ so that applying Lemma~\ref{lem:intro_trivial_lemma} and~\ref{StildetoS}, we obtain
\begin{align} \label{eq:cf-cg-1}
    \| C_f - C_g \| &  \leq \| S - T \| + \left( \| P^{-1} - I_d \| + \| N^{-1} - I_d \| \right) \| \tilde{S} - \tilde{T} \| + \| (P^{-1} - I_d) - (N^{-1} - I_d) \| \nonumber \\
    & \leq \left( 1 + 2 \left( \| P^{-1} - I_d \| + \| N^{-1} - I_d \| \right)  \right) \| S - T \| + \| P^{-1} - N^{-1}  \|
\end{align}
From \eqref{eq:proofMId}, \eqref{hepsilonbound2} and $\| P^{-1} \| \leq 2$ on one hand and from \eqref{damcfstabilityprf2}, the assumption on $g$ and $\| N^{-1} \| \leq 2$ on the other hand, we infer
\begin{equation} \label{eq:cf-cg-2}
\begin{split}
\| P^{-1} - I_d \| & \leq 2 n c_2 |F(x)|_\infty \leq 2n c_1 c_2 M \del \e^{-4} \\
\| N^{-1} - I_d \| & \leq 2 n c_2 |G(y)|_\infty \leq 2n c_1 c_2 M \del \e^{-4} \: .
\end{split}
\end{equation}
Furthermore, as $P^{-1} - N^{-1} = P^{-1} (N - P) N^{-1}$, we have
\begin{equation}\label{tildeM-tildeN}
 \| P^{-1} - N^{-1} \|  \leq \| P^{-1} \| \|N^{-1} \| \| P - N \| 
 \leq 4 \| P - N \| \: .
\end{equation}
We conclude the proof of \eqref{cf-cg} putting together \eqref{eq:cf-cg-1}, \eqref{eq:cf-cg-2}, \eqref{tildeM-tildeN} and \eqref{damcfstability1}.

\noindent We now prove
\begin{equation}\label{df-dg}
\| R_f - R_g \| \leq  c_5 M \del \left( \e^{-4} \| S -T \| +  \e^{-6} | x-y| \right) + c_5 \| Df_V - Dg \|_{\infty}.
\end{equation}
We note that by definition $P^t=P$ and then $(P^{-1})^t = P^{-1}$, and similarly, $(N^{-1})^t=N^{-1}$ so that (using $\| P^{-1} \| \leq 2 $ and $ \| N^{-1} \| \leq 2$),
\begin{equation*}
 \begin{split}
 & \| R_f - R_g \| \\
 & \quad = \| F(x)\circ \tilde{S}^t \circ P^{-1}  \circ \tilde{S} + \tilde{S}^t \circ P^{-1} \circ \tilde{S} \circ F(x)^t 
 - G(y) \circ \tilde{T}^t \circ N^{-1} \circ \tilde{T} - \tilde{T}^t \circ N^{-1} \circ \tilde{T} \circ G(y)^t \| \\
 & \quad \leq 2 \| F(x)\circ \tilde{S}^t \circ P^{-1} \circ \tilde{S} - G(y) \circ  \tilde{T}^t \circ N^{-1} \circ \tilde{T} \|
\end{split}
\end{equation*}
and decomposing the term above into $4$ terms and recalling $\| \tilde{S} - \tilde{T} \| \leq 2 \| S - T \|$ and \eqref{tildeM-tildeN}, we obtain
\begin{equation*}
 \begin{split}
 \| F(x)\circ & \tilde{S}^t \circ P^{-1} \circ \tilde{S}  - G(y) \circ  \tilde{T}^t \circ N^{-1} \circ \tilde{T} \| \\
 \leq & \|  F(x) \circ \tilde{S}^t \circ P^{-1} \circ (\tilde{S} -  \tilde{T} ) \|  + \|  F(x) \circ \tilde{S}^t \circ ( P^{-1} - N^{-1}) \circ \tilde{T} \| \\
&  + \| F(x) \circ (\tilde{S}^t - \tilde{T}^t) \circ N^{-1} \circ \tilde{T} \| +
\| (F(x) - G(y)) \circ \tilde{T}^t \circ N^{-1} \circ \tilde{T} \| \\
\leq &  2 \| F(x) \|  \| \tilde{S} - \tilde{T} \| +  \| F(x) \|  \| P^{-1} - N^{-1} \|  + 2 \|F(x)\| \|\tilde{S}-\tilde{T} \|+ 2 \| F(x) - G(y) \| \\
\leq &  \| F(x) \| \left( 8\| S - T \| + 4 \| P - N \| \right) +  2 \| F(x) - G(y) \| \: .
 \end{split}
\end{equation*}
The estimate \eqref{df-dg} then follows from \eqref{G-bound}, \eqref{damcfstability1} and \eqref{eq:Fx-Gy}.

\noindent We are left with checking that
\begin{equation}\label{ef-eg}
\| E_f - E_g \| \leq  c_5 M \del \left( \e^{-4} \| S -T \| +  \e^{-6} | x-y| \right) + c_5 \| Df_V - Dg \|_{\infty} \: .
\end{equation}
Indeed, recalling
$\|P^{-1}\|\leq 2$, $\| N^{-1}\|\leq 2$ (see \eqref{eq:boundinvPinvN}), $ \| F(x) \| \leq c_1 M \del \e^{-4} \leq 1$ and  $\| G(y) \| \leq 2 c_1 M \del \e^{-4} \leq 1$ (see \eqref{G-bound}), \eqref{eq:Fx-Gy} and using Lemma \ref{StildetoS} one has
\begin{align} \label{eq:ef-eg-0}
 & \| E_f - E_g \| =  \| F(x)\circ \tilde{S}^t \circ P^{-1} \circ \tilde{S} \circ F(x)^t - G(y) \circ \tilde{T}^t \circ N^{-1} \circ \tilde{T} \circ G(y)^t \| \nonumber \\
 & \quad \leq \| F(x) \circ \left( \tilde{S}^t \circ P^{-1}\circ \tilde{S} - \tilde{T}^t \circ N^{-1}\circ \tilde{T} \right) \circ F(x)^t \| 
 + \| (F(x) - G(y)) \circ \tilde{T}^t \circ N^{-1}\circ \tilde{T} \circ F(x)^t  \|\nonumber \\
 & \qquad  + \| G(y) \circ \tilde{T}^t \circ N^{-1} \circ \tilde{T} \circ (F(x)^t- G(y)^t  \| \nonumber \\
 & \quad \leq \| F(x) \|^2 \| \tilde{S}^t \circ P^{-1}\circ \tilde{S} - \tilde{T}^t \circ N^{-1}\circ \tilde{T} \| + 2 (\|F(x) \| + \| G(y) \|)  \| F(x) - G(y) \| \nonumber \\
 & \quad \leq  c_1 M \del \e^{-4} \| \tilde{S}^t \circ P^{-1}\circ \tilde{S} - \tilde{T}^t \circ N^{-1}\circ \tilde{T} \|  + 4 \left( c_1 M \del \e^{-6}  |x-y| + \| Df_V - Dg \|_{\infty} \right) \: .
\end{align}
We can then apply Lemma~\ref{lem:intro_trivial_lemma}, \eqref{tildeM-tildeN} and $\| \tilde{S} - \tilde{T} \| \leq 2 \| S - T \|$ to obtain
\begin{equation} 
\begin{split} \label{eq:ef-eg-1}
    \| \tilde{S}^t \circ P^{-1}\circ \tilde{S} - \tilde{T}^t \circ N^{-1}\circ \tilde{T} \| & \leq \left( \| P^{-1} \| + \| N^{-1} \| \right) \| \tilde{S} - \tilde{T} \| + \| P^{-1} - N^{-1} \| \\
    & \leq 8 \| S - T \| + 4 \| P - N \| \: .
\end{split}
\end{equation}
Therefore, \eqref{ef-eg} follows from \eqref{eq:ef-eg-0}, \eqref{eq:ef-eg-1} and \eqref{damcfstability1}.

\smallskip
\noindent{\bf Step 4:} We now prove \eqref{damcfstability5}.

\noindent Let $\phi \in \xC(\R^n \times \G, \R)$ satisfying $\|\phi\|_{\infty} \leq 1$ and  $\lip(\phi) \leq 1$. For $(x,S)\in \R^n \times \G$, we recall \eqref{eq:JSfbounds} $|J_S f_V(x)| \leq 2$ and thanks to \eqref{damcfstability4} we have 
\begin{align} \label{eq:pushf-g-0}
  & \big| \phi(f_V(x), Df_V(x)(S)) J_Sf_V(x)  - \phi(g(x),Dg(x)(S)) J_Sg(x) \big| \nonumber \\
  & \quad \leq \lip(\phi) \left( |f_V(x)-g(x)|+\|Df_V(x)(S)-Dg(x)(S)\|\right) |J_S f_V(x)|  +\|\phi\|_{\infty} |J_Sf_V(x)-J_Sg(x)| \nonumber \\
  & \quad \leq \| f_V - g \|_{\infty} + 2 c_5 \| Df_V - Dg \|_{\infty} + \| J_{\cdot}f_V - J_{\cdot}g\|_{\infty} \nonumber \\
  & \quad \leq c_5 \| f_V - g \|_{\xC^1} + \| J_{\cdot}f_V - J_{\cdot}g\|_{\infty}
\end{align}
where $\|\cdot\|_{\infty} $ is taken over all $(x,S) \in \R^n \times \G$ (and the last inequality holds up to doubling $c_5$). Integrating \eqref{eq:pushf-g-0} over $\R^n \times \G$ one has by definition of push-forward varifold (see Definition~\ref{dfn:varifoldpush-forward})
\begin{equation*}
\big| f_{V\#}V(\phi) - g_{\#}V(\phi) \big| \leq M \left( c_5 \| f_V- g \|_{\xC^1} + \|J_{\cdot}f_V-J_{\cdot}g\|_{\infty} \: \right).
\end{equation*}
Taking the supremum with respect to such Lipschitz functions $\phi$ leads to \eqref{damcfstability5} by definition of $\Delta$.

\smallskip
\noindent{\bf Step 5:} We are left with the proof of \eqref{damcfstability6}. We first apply \eqref{damcfstability5} with $V=W$, $f=f_{V}$ and $g=f_{W}$ so that
\begin{align}
 \Delta & \left( (f_{V})_{\#}W , (f_{W})_{\#}W \right)  \leq M \left(c_5\|f_{V}-f_{W}\|_{C^1} + \| J_{\cdot}f_{V}-J_{\cdot} f_{W}\|_{\infty} \right) \nonumber \\
 & \leq c_5 M \Delta t \| h_\e (\cdot, V) - h_\e(\cdot, W) \|_{\xC^1} + c_5 M^2 \del \e^{-n-7} \Delta(V,W) \text{ (by \eqref{damcfstability3} with $x=y$, $S=T$)}\nonumber \\
 & \leq (2c_4+ c_5)M^2 \Delta t \:  \e^{-n-7} \Delta (V,W) \text{ by Lemma~\ref{hepsilonstability}} \label{eq:stabilityDeltaPush2}
\end{align}
We now prove that
\begin{equation} \label{eq:stabilityDeltaPush0}
\Delta \left( (f_{V})_{\#} V, (f_{V})_{\#}W \right)  \leq \left( 1 + c_5 M \Delta t \e^{-6} \right) \Delta (V,W) \: .
\end{equation}
Let $\phi \in \xC_c^{0} (\R^n \times \G,\R)$ be a Lipschitz function satisfying $\|\phi\|_{\infty} \leq 1$ and $\lip(\phi) \leq 1$, then, coming back to the definition of $\Delta$ (Definition~\ref{dfnboundedlip}), we consider
$
\psi : (x,S) \mapsto \phi(f_{V}(x), Df_{V}(x)(S)) \: J_{S}f_{V}(x)$. As $f_V$ is a $\xC^1$- diffeomorphism, we have $\psi \in \xC_c (\R^n \times \G)$ and by definition of varifold push-forward,
\begin{equation} \label{eq:stabilityDeltaPush1}
\left| \int \phi \: d (f_V)_\# V - \int \phi \: d (f_V)_\# W \right| = \left| \int \psi \: dV - \int \psi \: dW \right| \leq \max (\|\psi\|_\infty, \lip(\psi)) \Delta (V,W) \: .
\end{equation}
One has to pay attention to the fact that the $\| \cdot \|_{\infty}$ and $\lip(\cdot)$ refer to both variables $(x,S) \in \R^n \times \G$.
Introducing the notations $\psi_1 : (x,S) \mapsto Df_V(x)(S) = \left(I_n + \Delta t \: D h_\e (x,V) \right) (S)$ and $\psi_2 : (x,S) \mapsto J_S f_V(x)$, we have
\begin{multline} \label{eq:stabilityLip}
 \lip(\psi_1) \leq (1 + c_5M \del \e^{-6})  \: , \quad \lip(\psi_2) \leq c_5  M \Delta t \: \e^{-6}  \: , \\
 \| \psi_2 \|_\infty \leq 1 + c_5 M \Delta t \: \e^{-4} \quad \text{and} \quad \lip(f_V) \leq 1 + c_5 M \Delta t \: \e^{-4} \: .
\end{multline}
\noindent Indeed, let $(x,S)$, $(t,y) \in \R^n \times \G$, and apply \eqref{damcfstability4} and \eqref{damcfstability3} with $f_V$ and $g=f_V$ so that
\begin{align*}
  | \psi_1(x,S) - \psi_1(t,y) | & \leq ( 1 + c_5 M \del \e^{-4} ) \| S-T \|  + c_5 \del \e^{-6} |x-y| \\
  & \leq ( 1 + c_5 M \del \e^{-6} ) ( \| S- T \| + |x-y| ) \\
 | \psi_2 (x,S) - \psi_2(t,y) | & \leq c_5 M \Delta t \: \left(\e^{-4} \| S-T \| + \e^{-6} |x-y|  \right) \: .
\end{align*}
Furthermore, by \eqref{eq:JSfbounds} and \eqref{hepsilonbound2}, $\displaystyle \| \psi_2 \|_\infty = \| J_\cdot f_V  \|_\infty \leq 1 + c_2 \Delta t \: \| D h_\e (\cdot,V) \|_\infty \leq 1 + c_1 c_2 M \Delta t \: \e^{-4}$ and $\lip(f_V) \leq 1 + \Delta t \: \lip(h_\e) \leq 1 + c_1 M \Delta t \: \e^{-4}$. 
With \eqref{eq:stabilityLip} in hand, we already note that
\begin{equation}\label{eq:maxPsi}
 \| \psi \|_\infty \leq \|\phi\|_\infty \|\psi_2\|_\infty \leq \|\psi_2\|_\infty \leq\left( 1 + c_5 M \Delta t \e^{-6} \right)
\end{equation}
and then, we can estimate $\lip(\psi) $ as follows:

\begin{align}
| \psi(x,S) & - \psi(t,y) | \nonumber \\
\leq & \left| \phi (f_V(x) , \psi_1 (x,S)) - \phi (f_V(y) , \psi_1 (t,y)) \right| \| \psi_2 \|_\infty + \| \phi \|_\infty \lip(\psi_2) \left( |x-y| + \| S - T \| \right) \nonumber \\
\leq & \left[ \lip(\phi) \left( \max\{ \lip(f_V) , \lip(\psi_1) \} \right) \| \psi_2 \|_\infty + \lip(\psi_2) \right] \left( |x-y| + \| S - T \| \right) 
 \nonumber \\
\leq & \left( 1 + c_5 M \del \e^{-6} \right) \big( 1 + c_5 M \Delta t \: \e^{-4}  \big) + c_5 M \Delta t \: \e^{-6}  
\leq  \left( 1+ 4 c_5 M \Delta t \e^{-6} \right) \: , \label{eq:lipPsi}
\end{align}
recalling that $c_5 M \Delta t \: \e^{-4} \leq 1$ by \eqref{smallstep2}. 
We conclude the proof of \eqref{eq:stabilityDeltaPush0} thanks to \eqref{eq:stabilityDeltaPush1}, \eqref{eq:maxPsi} and \eqref{eq:lipPsi} (up to increasing $c_5)$.
Combining \eqref{eq:stabilityDeltaPush2} and \eqref{eq:stabilityDeltaPush0} we conclude the proof of \eqref{damcfstability6}, and subsequently the proof of Proposition~\ref{damcfstability}.
\end{proof}
%
 
Iterating Proposition~\ref{damcfstability} leads to the following stability result on the time-discrete approximate mean curvature flow.
\begin{prop}[Stability with respect to initial data]
\label{stability}
Let $\e \in (0,1)$ and $M \geq 1 $. Let $V_0$ and $W_0$ be two varifolds in $V_d(\R^n)$ with $\|V_0\|(\R^n) \leq M,\|W_0\|(\R^n)\leq M$.\\
Let $\mathcal{T}=\lbrace t_i \rbrace_{i=0}^m$ be a subdivision of $[0,1]$ satisfying \eqref{smallstep}. Denote by $V_{\e,\cT}(t)$ (resp. $W_{\e,\cT}(t)$) the time-discrete approximate mean curvature flow with respect to $\cT$ starting from $V_0$ (resp. $W_0$) as introduced in Definition~\ref{damcf}. Then,
for any $t\in[0,1]$, one has
\begin{equation} \label{eq:stabilityFlowInitialData}
 \Delta\left(V_{\e,\cT}(t),W_{\e,\cT}(t) \right) \leq \exp(c_{5,M}t \: \e^{-n-7})\Delta(V_0,W_0),
\end{equation}
where $c_{5,M}=c_5(M+1)^2$ and $c_5$ was introduced in Proposition~\ref{damcfstability}.
\end{prop}
\begin{proof}
As $\e \in (0,1)$ and the subdivision $\cT$ are fixed, we write $V(t)$ (resp. $W(t)$) for $V_{\e,\cT}(t)$ (resp. $W_{\e,\cT}(t)$) hereafter.
From ~\eqref{damcfstability6} applied with $V=V(t_{i-1}), W=W(t_{i-1})$ and $\Delta t  = d_i = t_i - t_{i-1} $, and noting that $\| V(t_{i-1})\| (\R^n) \leq M+1$ and $\| W(t_{i-1})\| (\R^n) \leq M+1$ (see Remark~\ref{rem:uniformbound}), we infer that for any $i\in\lbrace 1 \dots,m\rbrace$ we have: 
 \begin{equation*}
  \Delta\left(V(t_i),W(t_i) \right) \leq (1+ c_{5,M} d_i \: \e^{-n-7})\Delta(V(t_{i-1}),W(t_{i-1})) \: .
 \end{equation*}
By iteration of the previous inequality for $k\in \lbrace 1 \dots,i\rbrace$ and applying the inequality $1 + a \leq \exp(a)$ in $\R$, we obtain
\begin{equation*}
 \Delta\left( V(t_i),W(t_i) \right) \leq \prod\limits_{k=1}^i (1+  c_{5,M} d_k \: \e^{-n-7})\Delta(V(0),W(0)) \underbrace{\leq}_{\sum_{k=1}^i d_k = t_i } \exp(c_{5,M} t_i \: \e^{-n-7}) \Delta(V_0,W_0) \: .
\end{equation*}
Let now $t \in (t_i,t_{i+1}]$ and apply once again Proposition~\ref{damcfstability} (with $\Delta t = t - t_i$) so that
\begin{align*}
 \Delta\left( V(t),W(t) \right) & \leq (1+c_{5,M}(t-t_i)\: \e^{-n-7})\Delta\left( V(t_i),W(t_i) \right) \\
 &\leq (1+c_{5,M}(t-t_i)\e^{-n-7})\exp(c_{5,M} t_i \: \e^{-n-7}) \Delta(V_0,W_0) 
 \\& \leq \exp(c_{5,M} t \: \e^{-n-7}) \Delta(V_0,W_0) \: ,
\end{align*}
thus ending the proof of the stability of the time-discrete approximate mean curvature flow with respect to initial data.
\end{proof}
\begin{remk}[Analogy with ODE discretization]
\label{remk:ODEnum}
The construction of the time-discrete approximate mean curvature flow defined in this paper can be compared to the discretization of the classical Cauchy problem in $[0,T]$:
\begin{equation*}
\left\{ \begin{array}{ll}
 y'(t)=f(t,y),\\
 y(0)=y_0.
\end{array} \right.
\end{equation*}
It is known that a stability constant (with respect to the supremum norm on $[0,T]$) for the explicit Euler discretization of the ODE is $\exp(LT)$ with $L=\max\limits_{t \leq T} \lip(f(\cdot,t))$ (see for instance \cite[Section 2.4]{ahs}).
Comparing with the stability estimate \eqref{eq:stabilityFlowInitialData} obtained in
Proposition~\ref{stability}, we observe that $c_{5,M} \e^{-n-7}$ is indeed a bound on the Lipschitz constant of $V \mapsto H_\e (\cdot , V)$ when $V_d(\R^n)$ is endowed with the bounded Lipschitz distance $\Delta$, and $\xC^1 (\R^n , \R^n)$ is endowed with $\| \cdot \|_{\xC^1}$, see Lemma~\ref{hepsilonstability}.
\end{remk}

%
%
\subsection{Stability with respect to time subdivision}
\label{subsec:stabilityTimeSubdivision}

In this section, we investigate the robustness of the time-discrete approximate mean curvature flow (introduced in Definition~\ref{damcf}) with respect to the choice of the subdivision $\cT$. It is a natural property to expect for a numerical scheme and it is furthermore crucial in order to take the limit "$\delta(\cT) \to 0$'' and obtain a well-defined ``time-continuous'' approximate mean curvature flow as subsequently done in Theorem~\ref{damcfconvergence}.
We establish in Proposition~\ref{stabilitysub} that time-discrete approximate mean curvature flows are stable with respect to subdivisions. The proof of Proposition~\ref{stabilitysub} is split into several steps:
the section starts with two lemmas (Lemma~\ref{introlemma1} and \ref{introlemma2}) aiming to compare two flows corresponding to a fine subdivision and the trivial subdivision of small time interval $[0,\delta]$. Then, in Lemma ~\ref{preconvergence} we extend the comparison to the case of two nested subdivisions of the interval $[0,1]$ from which Proposition~\ref{stabilitysub} can be inferred straightforwardly.
\begin{prop}[Stability with respect to time subdivision]\label{stabilitysub}
Let $\e\in(0,1)$ and  $M\geq 1$. Let $V_0\in V_d(\R^n)$ with  $\|V_0\|(\R^n)\leq M$, let $\cT_1=\lbrace t_i \rbrace_{i=0}^m$ and $\cT_2=\lbrace s_j \rbrace_{j=0}^{m'}$ be two subdivisions of $[0,1]$ satisfying \eqref{smallstep}. Let $V_{\e,\cT_1}(t)$ (resp.$V_{\e,\cT_2}(t))$ be the time-discrete approximate mean curvature flow with respect to $\cT_1$ (resp. $\cT_2$) starting from $V_0$. We set: 
\begin{equation*}
    \delta = \max \big\lbrace  \delta(\cT_1) , \delta(\cT_2) \big\rbrace.
\end{equation*}
Then, for all $t\in [0,1]$, one has: 
\begin{equation*}
\Delta(V_{\e,\cT_1}(t),V_{\e,\cT_2}(t)) \leq  c_{7,M} \, t \, \delta \,  \e^{-n-11} \exp(c_{5,M} \, t \, \e^{-n-7}).
\end{equation*}
where $c_{7,M}=c_7 (M+1)^5$, $c_{5,M} = c_5 (M+1)^2$ and $c_7$, $c_5$ are constants depending only on $n$ introduced in Proposition~\ref{stability} and Lemma~\ref{introlemma2}.
\end{prop}
Before proving Proposition~\ref{stabilitysub}, we shall introduce some preliminary lemmas. 

\par We first estimate in the following lemma how far the push-forward operation is from satisfying the semigroup property. In practice, we measure how far apart are two time-discrete approximate mean curvature flows constructed with respect to two subdivisions including the trivial subdivision.
\begin{lemma}\label{introlemma1}
Let $\e, \delta \in(0,1)$ and $M\geq1$ such that $c_3 (M+1)^3 \delta  \e^{8} \leq 1$. Let $V_0 \in V_d(\R^n)$ satisfy $\|V_0\|(\R^n)\leq M$. Consider $\cT=\lbrace t_i \rbrace_{i=0}^m$ a given subdivision of $[0,\delta]$ and $\cT'$ the trivial subdivision of $[0,\delta]$. For $i \in \lbrace 1, \dots m \rbrace$, we introduce
\begin{equation*}
d_i=t_i-t_{i-1} \quad \text{and} \quad \tilde{f}_i = \left( {\rm id} +d_i \: h_{\e}(\cdot,V_0)\right) \: .
\end{equation*}
We then consider two different flows:
\begin{itemize}
 \item $\left( V_{\e,\cT'}(t_i) \right)_{i = 0 \ldots m}$ where $ V_{\e,\cT'}$ is the time-discrete approximate mean curvature flow of $V_0$ with respect to $\cT'$ according to Definition~\ref{damcf},
 \item $\left( \tV_{\e,\cT}(t_i) \right)_{i = 0 \ldots m}$ is defined as follows:
 \begin{equation}\label{auxflow}
\left\{
 \begin{array}{ll}
\tV_{\e,\cT}(0) :=V_{0}\\
\tV_{\e,\cT}(t_{i}) :=(\tilde{f}_{i})_{\#} \tV_{\e,\cT}(t_{i-1}) \quad \forall  i\in \lbrace 1, \dots m \rbrace.
 \end{array} \right.
\end{equation}
\end{itemize}
%
Then, there exists $c_6 > 0$ that depends only on $n$ and such that
\begin{equation*}
  \Delta \left(\tV_{\e,\cT}(t_i),V_{\e,\cT'}(t_i)\right) \leq c_6 M^3t_i^2 \e^{-10} \quad \forall \,\, i \in\lbrace 0, \dots m \rbrace \: .
\end{equation*}
%
\end{lemma}

It is important to note that, using the previous notation with $\Delta t = d_i$, we have $\tilde{f}_i = f_{\e, V_0}$, while the definition of time-discrete approximate mean curvature flow would involve $f_{\e, V_{\e,\cT}(t_{i-1})}$ instead. The velocity $h_\e$ is taken with respect to the initial varifold all along the subdivision when defining $\tilde{V}_{\e,\cT}$.


\begin{proof}
As previously, the constant $c_6 > 0$ is adapted throughout the proof, provided that it depends only on $n$.
We introduce the following notations: $g_0 = \tilde{g}_0 = f_0 = \tilde{f_0} = {\rm id}$ and
$$ \forall  i\in \lbrace 1, \dots m \rbrace \quad g_i={\rm id}+t_ih_{\e}(\cdot,V_0) \quad \text{and} \quad \tilde{g}_i=\tilde{f}_i \circ \dots \circ \tilde{f}_1 \: .$$
We first prove that for all $i \in \{ 1, \ldots, m \}$,
\begin{equation}\label{ggtildec1}
 \|\tilde{g}_{i}- g_i\|_{\xC^1} \leq c_6 \, M^2 t_i^2 \e^{-10} \: .
\end{equation}
%
Indeed, let $i \in \lbrace 1, \dots m \rbrace$ and $x \in \R^n$, we have by definition
\begin{align}\label{gtilde-id}
    \tilde{g}_i(x) & =
    \tilde{f}_i \circ \ldots \circ \tilde{f}_1(x)  = \tilde{f}_i \left( \tilde{g}_{i-1} (x) \right) =  \tilde{g}_{i-1}(x)+d_i \: h_{\e}(\tilde{g}_{i-1}(x),V_0) 
    = x + \sum_{k=1}^i d_k h_\e \left( \tilde{g}_{k-1}(x) , V_0 \right) \: , \nonumber \\
    & \text{hence} \quad \left| \tilde{g}_i (x) - x \right| \leq \sum_{k=1}^i d_k \left|h_\e\left( \tilde{g}_{k-1}(x) , V_0 \right) \right| \leq t_i \left\| h_\e ( \cdot, V_0 ) \right\|_\infty, \\
    g_i(x) & = x + t_i \: h_\e (x , V_0) = x + \sum_{k=1}^i d_k \: h_\e \left( x , V_0 \right) \nonumber
    \: ,
\end{align}
and applying the mean value theorem we infer
\begin{align}
|\tilde{g}_{i}(x)- g_i(x) | & \leq \sum_{k=1}^i d_k \left|  h_\e \left( \tilde{g}_{k-1}(x) , V_0 \right)  - h_\e \left( x , V_0 \right) \right| \leq \sum_{k=1}^i d_k \| Dh_\e (\cdot ,V_0) \|_\infty \left|  \tilde{g}_{k-1}(x) - x \right| \nonumber \\
& \leq t_i \: t_{i-1} \| Dh_\e (\cdot ,V_0) \|_\infty \| h_\e (\cdot , V_0) \|_\infty \nonumber \\
& \leq c_1^2 \, t_i^2  M^2 \e^{-6} \quad \text{thanks to Proposition~\ref{hepsilonbound}}. \label{eq:gitilde_gi_sup}
\end{align}
\noindent We proceed similarly to bound the derivatives but we have to handle the term $D \tilde{g}_i(x)$ arising from the chain rule applied to $x \mapsto h_\e (\tilde{g}_i(x), V_0)$: recalling that $\tilde{g}_{i} = ({\rm id} + d_i \: h_\e (\cdot, V_0)) \circ \tilde{g}_{i-1}$, we infer for $x \in \R^n$ and $i \in \{1, \ldots , m \}$,
\begin{align}
\| D \tilde{g}_i (x) \| & \leq \| {\rm id} + d_i \: D h_\e (\tilde{g}_i(x), V_0) \| \| D \tilde{g}_{i-1} (x) \| \leq (1 + d_i \: c_1 M \e^{-4} ) \| D \tilde{g}_{i-1} (x) \| \text{ using \eqref{hepsilonbound2}} \nonumber \\
& \leq \prod_{k=1}^i \left( 1 + d_k \: c_1 M \e^{-4} \right) \| D \tilde{g}_{0} (x) \| \leq \prod_{k=1}^i \exp \left( d_k \: c_1 M \e^{-4} \right) \nonumber \\
& \leq \exp \left( t_i c_1 M \e^{-4} \right) \leq 2\label{eq:DtildegiBound}
\end{align}
where we used $\tilde{g}_0 = {\rm id}$, $1 + s \leq \exp(s)$ for $s \in \R$ and 
\begin{equation*}
 \exp \left( t_i c_1 M \e^{-4} \right) \leq \exp \left( \delta c_3 (M+1)^3 \e^{-8} \right)\leq e \leq 3.
\end{equation*}

\noindent We can now expand $D\tilde{g}_i$ and $D g_i$ as we did previously for $\tilde{g}_i$ and $g_i$:
\begin{equation*}
    D \tilde{g}_i(x) 
    = x + \sum_{k=1}^i d_k \: D h_\e \left( \tilde{g}_{k-1}(x) , V_0 \right) \circ D \tilde{g}_{k-1}(x) \: , 
\end{equation*}
so that \eqref{hepsilonbound2} and \eqref{eq:DtildegiBound} imply
\begin{align}\label{Dgtilde-id}
    \left\| D \tilde{g}_i (x) - {\rm id} \right\| \leq \sum_{k=1}^i d_k \left\| D h_\e\left( \tilde{g}_{k-1}(x) , V_0 \right) \right\| \| D \tilde{g}_{k-1}(x) \| \leq 3 c_1 t_i M \e^{-4} \: , 
\end{align}
and
\begin{equation*}
D g_i(x)  = {\rm id} + t_i \: D h_\e (x , V_0) = {\rm id} + \sum_{k=1}^i d_k \: D h_\e \left( x , V_0 \right). 
\end{equation*}
We can then apply the mean value theorem to $Dh_\e(\cdot, V_0)$, together with \eqref{gtilde-id}, \eqref{Dgtilde-id} and Proposition \ref{hepsilonbound} to infer 
\begin{align}
\| D\tilde{g}_{i}(x)- D g_i(x) \| & \leq \sum_{k=1}^i d_k \left\|  D h_\e \left( \tilde{g}_{k-1}(x) , V_0 \right) \circ D \tilde{g}_{k-1}(x) - D h_\e \left( x , V_0 \right) \right\| \nonumber \\
& \leq \sum_{k=1}^i d_k \| Dh_\e (\cdot ,V_0) \|_\infty \left\|  D \tilde{g}_{k-1}(x) - {\rm id} \right\| + \| D^2 h_\e (\cdot , V_0 ) \|_\infty \left| \tilde{g}_{k-1}(x) -  x \right| \nonumber \\
& \leq 3 t_i \: t_{i-1}  (c_1 M \e^{-4})^2  + c_1 M \e^{-4} t_i t_{i-1} c_1 M \e^{-6} \nonumber \\
& \leq 4 c_1^2 M^2 \: t_i^2 \: \e^{-10}\: \label{eq:Dgitilde_gi_sup} \: .
\end{align}
We conclude the proof of \eqref{ggtildec1} thanks to \eqref{eq:gitilde_gi_sup} and \eqref{eq:Dgitilde_gi_sup}.

\noindent By definition of $\tV_{\e,\cT}(t_i)$ (see \eqref{auxflow}) and Lemma~\ref{basiclemma} we have
$$
 \tV_{\e,\cT}(t_i) = {(\tilde{f}_i)}_\# \left( (\tilde{f}_{i-1})_\# \ldots_\# \left( (\tilde{f_1})_\# V_0 \right) \right) = (\tilde{g}_{i})_\# V_0 \: ,
$$
and thanks to \eqref{Dgtilde-id} we can apply Proposition~\ref{damcfstability} with $V=V_0$, $\|V_0\|(\R^n) \leq M$ $f=g_i$ and $g=\tilde{g}_i$ so that by \eqref{damcfstability5}
\begin{equation} \label{eq:deltaVtildeV_0}
 \Delta \left(\tV_{\e,\cT}(t_i),V_{\e,\cT'}(t_i)\right)= \Delta \left((\tilde{g}_{i})_\# V_0,(g_{i})_\# V_0\right) \leq M \left( c_5 \| g_i - \tilde{g}_i \|_{\xC^1} + \|J_{\cdot}g_i - J_{\cdot}\tilde{g}_i \|_{\infty} \right)
\end{equation}
and by \eqref{damcfstability2}
\begin{equation}\label{Jgitilde_gi}
\begin{split}
 \| J_{\cdot} g_i - J_{\cdot} \tilde{g}_i \|\leq  c_6 \, M^2 t_i^2 \e^{-10} \: .
\end{split}
\end{equation}
We conclude the proof of Lemma~\ref{introlemma1} thanks to \eqref{eq:deltaVtildeV_0}, \eqref{ggtildec1} and \eqref{Jgitilde_gi}.
\end{proof}

In the following lemma we compare the time-discrete approximate mean curvature flows of two given varifolds on a small time interval, one of them being defined with respect to the trivial subdivision and the other one being defined with respect to a finer subdivision of the time interval. The proof relies on Lemma \ref{introlemma1}.
\begin{lemma}\label{introlemma2}
Let $\e, \, \delta \in (0,1)$ and $M\geq1$. Let $V_0$ and $W_0$ in $V_d(\R^n)$ be two varifolds satisfying $\|V_0\|(\R^n) \leq M$ and $\|W_0\|(\R^n)\leq M$. Let $\cT=\lbrace t_i \rbrace_{i=1}^m$ be a subdivision of $[0,\delta]$ and denote by
\begin{enumerate}[label=$\bullet$]
\item $V_{\e,\cT'}(t)$ the time-discrete approximate mean curvature flow of $V_0$ w.r.t. the trivial subdivision $\cT' = \{0,\delta \}$ of $[0,\delta]$,
\item $W_{\e,\cT}(t)$ the time-discrete approximate mean curvature flow of $W_0$ w.r.t. $\cT$.
\end{enumerate}
If 
$$c_3\delta (M+1)^3 < \e^{8}$$
then, for any $i \in \lbrace 0,1 \dots m \rbrace$:
\begin{equation}\label{convergenceestimate0}
\Delta(V_{\e,\cT'}(t_i),W_{\e,\cT}(t_i))\leq \left( \Delta(V_0,W_0)  + c_{7,M} \,t_i^2 \, \e^{-n-11} \right) \exp(c_{5,M} \, t_i \,  \e^{-n-7})
\end{equation}
where $c_7 > 0$ is a constant that depends only on $n$ and $c_{7,M}=c_7 (M+1)^5$.
\end{lemma}
\begin{remk}
A particular case of \eqref{convergenceestimate0} is when $i=m$: 
\begin{equation}\label{convergenceestimate1}
 \Delta(V(\delta),W(\delta))\leq \left( \Delta(V_0,W_0)  +  c_{7,M} \, \delta^2 \, \e^{-n-11} \right) \exp(c_{5,M} \, \delta \, \e^{-n-7}).
\end{equation}
The last inequality will be useful in the sequel.
\end{remk}

\begin{proof}As previously, the constant $c_7 > 0$ is adapted throughout the proof, provided that it depends only on $n$.
The assumption $c_3 \delta (M+1)^3< \e^{8}$ implies that the involved subdivisions ($\cT'$ and $\cT$) satisfy \eqref{smallstep}, which allows to define time-discrete approximate mean curvature flows for both subdivisions. For every $ i \in \lbrace 1 \dots m \rbrace$, set:
\begin{equation*}
d_i=t_i-t_{i-1}
\end{equation*}
and define $\tilde{V}_{\e,\cT}(t)$ the auxiliary flow as in \eqref{auxflow}. We have by Lemma \ref{introlemma1} 
\begin{equation}\label{introlemma20}
\begin{split}
\Delta(V_{\e,\cT'}(t_i),W_{\e,\cT}(t_i))
& \leq  \Delta(V_{\e,\cT'}(t_i),\tilde{V}_{\e,\cT}(t_i)) +  \Delta(\tilde{V}_{\e,\cT}(t_i),W_{\e,\cT}(t_i)) 
\\& \leq c_6 M^3t_i^2 \e^{-10}+ \Delta(\tilde{V}_{\e,\cT}(t_i),W_{\e,\cT}(t_i)).
\end{split}
\end{equation}
We are left with estimating $\Delta(\tilde{V}_{\e,\cT}(t_i),W_{\e,\cT}(t_i))$. As in Lemma~\ref{introlemma1}, we introduce the notations $\tilde{f}_i= {\rm id}+d_i \, h_{\e}(\cdot,V_0)$ and $g_i = {\rm id}+ t_i \, h_{\e}(\cdot,V_0)$, and we also use the following shortened notations up to the end of the current proof: for $i \in \{0,1, \ldots, m\}$,
\begin{multline*}
    W_i = W_{\e,\cT}(t_i) \quad \text{and} \quad \tilde{V}_i = \tilde{V}_{\e,\cT}(t_i) \quad  \text{and} \quad V^\prime_i = V_{\e,\cT'}(t_i) \\ \text{and then} \quad f_{\e, W_i} = {\rm id} + d_{i+1} h_\e (\cdot, W_i ) \quad \text{and} \quad  f_{\e, \tilde{V}_i} = {\rm id} + d_{i+1} h_\e (\cdot, \tilde{V_{i}} )
\end{multline*}
so that for $l \in \{1, \ldots, m\}$,
 \begin{equation*}
  \tilde{V}_l =  (\tilde{f}_{l})_\# \tilde{V}_{l-1} \quad \text{and} \quad W_l = (f_{\e, W_{l-1}})_\# W_{l-1} \: .
 \end{equation*}
Coming back to \eqref{introlemma20} we have by triangular inequality
\begin{align}
    \Delta \left( \tilde{V}_{\e,\cT}(t_l),W_{\e,\cT}(t_l) \right) & = \Delta \left( \tilde{V}_l ,W_l \right) = \Delta \left( (\tilde{f}_{l})_\# \tilde{V}_{l-1}, (f_{\e, W_{l-1}})_\# W_{l-1} \right) \nonumber \\
    & \leq \Delta \left( (\tilde{f}_{l})_\# \tilde{V}_{l-1}, (f_{\e, \tilde{V}_{l-1}})_\# \tilde{V}_{l-1} \right) + \Delta \left( (f_{\e, \tilde{V}_{l-1}})_\# \tilde{V}_{l-1}, (f_{\e, W_{l-1}})_\# W_{l-1} \right) \: . \label{eq:VtildeW0}
\end{align}
We can directly handle the second term in the previous inequality thanks to \eqref{damcfstability6} that we apply with $V = \tilde{V}_{l-1}$, $ W=W_{l-1}$  and $\del = d_l$ so that
\begin{equation} \label{eq:VtileWrec1}
    \Delta \left( (f_{\e, \tilde{V}_{l-1}})_\# \tilde{V}_{l-1}, (f_{\e, W_{l-1}})_\# W_{l-1} \right)  \leq \left(1+c_{5,M} \, d_l \, \e^{-n-7}\right) \Delta \left( \tilde{V}_{l-1},W_{l-1} \right)  \: .
\end{equation}
We now apply \eqref{damcfstability5} and \eqref{damcfstability3} with $V=\tilde{V}_{l-1}$, $g=\tilde{f}_l$ and $\del = d_l$ 
so that together with Lemma~\ref{hepsilonstability} and \eqref{uniformbound}, we can assert that:
\begin{equation}\label{introlemma21}
\begin{split}
\Delta \left( (\tilde{f}_{l})_\# \tilde{V}_{l-1}, (f_{\e, \tilde{V}_{l-1}})_\# \tilde{V}_{l-1} \right) & \leq (M+1) \left( 
c_5 \, d_l \, \| h_{\e}(\cdot,V_0)-h_{\e}(\cdot,\tilde{V}_{l-1}) \|_{C^1} + \|J_{\cdot}\tilde{f}_l-J_{\cdot} f_{\e, \tilde{V}_{l-1}} \|_{\infty} \right) \\
& \leq (2c_4 +1) c_5 \, (M+1)^2 \, d_l \, \e^{-n-7} \Delta ( \tilde{V}_{l-1},V_0 ) \: .
\end{split}
\end{equation}
We recall that $\cT'$ being the trivial subdivision, we have $V^\prime_{l-1} = (f_{\e,V_0})_\# V_0$ with $(f_{\e,V_0})_\# V_0 = {\rm id} + t_{l-1} h_\e (\cdot, V_0)$ and we can thus apply \eqref{damcfstability5} once more, with $V = V_0$, $g = {\rm id}$ and $\Delta t = t_{l-1}$:
\begin{align*}
    \Delta(V^\prime_{l-1}, V_0) & \leq M \left( c_5 \, t_{l-1} \|h_{\e}(\cdot,V_0)\|_{C^1} 
 + \|J_{\cdot} f_{\e,V_0} - 1 \|_{\infty}\right) \nonumber \\
 & \leq c_7 \, M^2 \, t_{l-1} \, \e^{-4} \quad \text{thanks to \eqref{eq:JSfbounds} and Proposition~\ref{hepsilonbound}.} \label{introlemma23}
\end{align*}
Consequently, by triangular inequality and Lemma~\ref{introlemma1}, we obtain (with $c_3 M^3 t_{l-1} \leq \e^8$)
\begin{equation}\label{introlemma22}
\Delta( \tilde{V}_{l-1} , V_0)  \leq \Delta( \tilde{V}_{l-1} , V^\prime_{l-1}  ) + \Delta( V^\prime_{l-1} , V_0) 
 \leq c_6 \, M^3 t_{l-1}^2 \e^{-10} + \Delta( V^\prime_{l-1} , V_0) \leq c_7 \, M^3 t_{l-1} \e^{-4} 
\end{equation}
Collecting \eqref{eq:VtildeW0}, \eqref{eq:VtileWrec1} and \eqref{introlemma21}, \eqref{introlemma22} we obtain the inductive relation:
\begin{equation} \label{introlemma24}
    \Delta ( \tilde{V}_l , W_l)  \leq \underbrace{ \left(1+c_{5,M} \, d_l \, \e^{-n-7}\right) }_{= a_l} \Delta ( \tilde{V}_{l-1},W_{l-1} ) + \underbrace{c_{7,M} \,  d_l \, t_{l-1} \,  \e^{-n-11} }_{=b_l} \: .
\end{equation}
Iterating \eqref{introlemma24} for $l \in \lbrace 1, \dots i\rbrace$, we infer
\begin{equation*}
    \Delta (\tilde{V}_i ,W_i ) \leq \Delta(V_0,W_0) \, \prod\limits_{l=1}^i a_l   + \sum_{l=0}^{i-1} b_{i-l} \, \prod_{j=0}^{l-1} a_{i-j}
\end{equation*}
We first note that for $l \leq i$
\begin{equation*}
 \prod\limits_{j=0}^{l-1} a_{i-j}  \leq
 \prod\limits_{j=0}^{i-1} a_{i-j} \prod\limits_{j=1}^{i} a_{j} =  \prod\limits_{j=1}^{i}\left(1+c_{5,M} \, d_{j}\e^{-n-7}\right) \leq \exp(c_{5,M} \, t_i \, \e^{-n-7}).
\end{equation*}
where we used again that for $x \in \R$, $1 + x \leq \exp(x)$ and $\sum_{j=1}^{i} d_{j} = t_i$. As we obtain an upper bound that does not depend on $l$, then
\begin{equation*}
    \sum_{l=0}^{i-1} b_{i-l} = c_{7,M} \, \e^{-n-11} \sum_{l=0}^{i-1} d_{i-l} t_{i-l-1} \leq c_{7,M} \, \e^{-n-11} t_i \sum_{l=0}^{i-1} d_{i-l} = c_{7,M} \,  t_i^2 \, \e^{-n-11} \: ,
\end{equation*}
yields
\begin{equation*}
    \Delta ( \tilde{V}_l , W_l)  \leq c_{7,M} \,  t_i^2 \, \e^{-n-11}  \left( \Delta ( V_0 , W_0) + \exp(c_{5,M} \, t_i \, \e^{-n-7}) \right) \: .
\end{equation*}
Coming back to \eqref{introlemma20} we conclude that  
 $ \forall i \in \lbrace 0,1 \dots m \rbrace$,
\begin{equation*}
\begin{split}
  \Delta(V_{\e,\cT'}(t_i),W_{\e,\cT}(t_i))& \leq c_6M^3t_i^2\e^{-10} +c_{7,M} \,  t_i^2 \, \e^{-n-11}  \left( \Delta ( V_0 , W_0) + \exp(c_{5,M} \, t_i \, \e^{-n-7}) \right) \\
  & \leq  c_{7,M} \,  t_i^2 \, \e^{-n-11}  \left( \Delta ( V_0 , W_0) + \exp(c_{5,M} \, t_i \, \e^{-n-7}) \right) \: .
\end{split}
\end{equation*}
\end{proof}
In the following lemma, we use \ref{introlemma2} to show Proposition  \ref{stabilitysub} (stability with respect to subdivision) in the special case where the two subdivisions are nested (one included in the other).
\begin{lemma}\label{preconvergence}
 Let $\e\in(0,1)$, $M\geq 1$. Let $V_0 \in V_d(\R^n)$ satisfy $\|V_0\|(\R^n)\leq M$. Consider $\cT_1=\lbrace t_i \rbrace_{i=1}^m$ and $\cT_2=\lbrace s_j \rbrace_{j=1}^{m'}$ two subdivisions of $[0,1]$ satisfying
 \eqref{smallstep}, assume that $\cT_1 \subset \cT_2$ and set $\delta = \delta(\cT_1)$.
Then, for all $t\in[0,1]$,
\begin{equation*}
 \Delta(V_{\e,\cT_1}(t),V_{\e,\cT_2}(t)) \leq c_{7,M} \, t \, \delta \, \e^{-n-11} \exp(c_{5,M} \, t \, \e^{-n-7})
\end{equation*}
where $V_{\e,\cT_1}(t)$ denotes the time-discrete approximate mean curvature flow of $V_0$ with respect to $\cT_1$ and $V_{\e,\cT_2}(t)$ the time-discrete approximate mean curvature flow of $V_0$ with respect to $\cT_2$.
\end{lemma}
\begin{proof}
\noindent We set for $i \in \lbrace 1,\dots,m \rbrace$, $d_i = t_i - t_{i-1}$

\noindent {\bf Step 1:} We first compare the flows at time $t_i$ of the subdivision $\cT_1$, more precisely, we show that for $i \in \lbrace 0,1 \dots ,m \rbrace$,
\begin{equation}\label{convergenceestimate2}
 \Delta (V_{\e,\cT_1}(t_i),V_{\e,\cT_2}(t_i)) \leq c_{7,M} \, t_i \, \delta \e^{-n-11}\exp(c_{5,M} \, t_i \, \e^{-n-7}).
\end{equation}
Fix $i \in \lbrace 0,1 \dots ,m \rbrace$, for any $l \in \lbrace 1, \dots ,i \rbrace$, we use \eqref{convergenceestimate1} on the interval $[t_{l-1},t_l]$ of length $d_l = t_l - t_{l-1}$ with $V_0=V_{\e,\cT_1}(t_{l-1})$, $W_0=V_{\e,\cT_2}(t_{l-1})$, $\cT'$ being the trivial subdivision of $[t_{l-1},t_l]$ and $\cT= \cT_2 \cap [t_{l-1},t_l]$ to obtain
\begin{equation*}
    \Delta(V_{\e,\cT_1}(t_l),V_{\e,\cT_2}(t_l)) \leq \left[ \Delta(V_{\e,\cT_1}(t_{l-1}),V_{\e,\cT_2}(t_{l-1}))  + c_{7,M} \, d_l^2 \,  \e^{-n-11} \right] \exp(c_{5,M} d_l \e^{-n-7}).
\end{equation*}
Noting that $u_l = \exp( - c_{5,M} \, t_l \, \e^{-n-7}) \Delta(V_{\e,\cT_1}(t_l),V_{\e,\cT_2}(t_l))$ hence satisfies the inductive relation
\begin{equation*}
 u_l \leq u_{l-1} + c_{7,M} \, d_l^2 \,  \e^{-n-11}
\end{equation*}
and we infer that $u_i \leq u_0 + c_{7,M} \, \sum_{l=0}^i d_l^2 \, \e^{-n-11}  \leq u_0 + c_{7,M} \, \delta t_i \, \e^{-n-11} = c_{7,M} \, \delta t_i \, \e^{-n-11}$ since $u_0 = \Delta(V_{\e,\cT_1}(0),V_{\e,\cT_2}(0))= 0$ yielding \eqref{convergenceestimate2}.

\noindent {\bf Step 2:} We now compare the flows at time $s_k$ of the finer subdivision $\cT_2$, more precisely, for $k \in \lbrace 0,1 \dots ,m' \rbrace$,
\begin{equation}\label{convergenceestimate3}
  \Delta(V_{\e,\cT_1}(s_k),V_{\e,\cT_2}(s_k))  \leq  c_{7,M} \, s_k \delta \, \e^{-n-11} \exp(c_{5,M} \, s_k \, \e^{-n-7}).
\end{equation}
Fix $k \in \lbrace 0,1 \dots,m' \rbrace$ and for $s_k \in \cT_2$, let $i \in \lbrace 0,1 \dots m-1 \rbrace$ be such that $s_k \in [t_i,t_{i+1}]$ .
Then, applying Lemma~\ref{introlemma2} on the interval $ [t_i,t_{i+1}]$ with $V_0 = V_{\e,\cT_1}(t_i)$, $W_0= V_{\e,\cT_2}(t_i)$, $\cT'$ being the trivial subdivision of $[t_i,t_{i+1}]$ and $\cT= \cT_2 \cap [t_i,t_{i+1}]$, we obtain
\begin{align*}
\Delta & (V_{\e,\cT_1}(s_k),V_{\e,\cT_2}(s_k))  \\
& \leq \left[ \Delta(V_{\e,\cT_1}(t_i),V_{\e,\cT_2}(t_i)) + c_{7,M} \, (s_k-t_i)^2 \, \e^{-n-11}\right] \exp(c_{5,M} \, (s_k - t_i) \e^{-n-7}) \\
& \leq c_{7,M} \, t_i \, \delta \e^{-n-11} \exp(c_{5,M} \, s_k \, \e^{-n-7}) + c_{7,M} \, (s_k-t_i)^2 \, \e^{-n-11} \exp(c_{5,M} \, (s_k-t_i) \, \e^{-n-7}) \: ,
\end{align*}
thanks to \eqref{convergenceestimate2}, that leads to \eqref{convergenceestimate3} noting that $\exp(c_{5,M} (s_k-t_i) \e^{-n-7}) \leq \exp(c_{5,M}s_k\e^{-n-7})$ and
$t_i\delta +(s_k-t_i)^2 \leq \delta(t_i+s_k-t_i) \leq \delta s_k$.

\noindent {\bf Step 3:} 
Let $t \in [s_k,s_{k+1}]$ for some $k \in \lbrace 0,1, \dots , m'-1 \rbrace$, applying Proposition~\ref{damcfstability} with $V=V_{\e,\cT_1}(s_k) $, $W= V_{\e,\cT_2}(s_k)$ and $ \Delta t = t-s_k$ we obtain: 
\begin{align*}
 \Delta(V_{\e,\cT_1}(t), & V_{\e,\cT_2}(t)) \\
 \leq & \Delta \Bigl( \left({\rm id}+(t-s_k)h_{\e}(\cdot,V_{\e,\cT_1}(s_k))\right)_{\#}V_{\e,\cT_1}(s_k),\left({\rm id}+(t-s_k)h_{\e}(\cdot,V_{\e,\cT_2}(s_k))\right)_{\#}V_{\e,\cT_2}(s_k) \Bigr)\\
 \leq & (1+(t-s_k)c_{5,M}\e^{-n-7})\Delta(V_{\e,\cT_1}(s_k),V_{\e,\cT_2}(s_k)) \\
 \leq & \exp\left((t-s_k)c_{5,M} \e^{-n-7}\right)\Delta(V_{\e,\cT_1}(s_k),V_{\e,\cT_2}(s_k)).
\end{align*}
From \eqref{convergenceestimate3} we conclude that for all $t\in[0,1]$.
\begin{align*}
\Delta(V_{\e,\cT_1}(t),V_{\e,\cT_2}(t))
& \leq \exp\left(c_{5,M} \, (t-s_k) \e^{-n-7}\right) \, c_{7,M} \, s_k\delta \e^{-n-11} \exp(c_{5,M}s_k\e^{-n-7})
\\& \leq  c_{7,M} \, t \, \delta \, \e^{-n-11}\exp(c_{5,M} \, t \,  \e^{-n-7})
\end{align*}
and this ends the proof of Lemma \ref{preconvergence}.
\end{proof}

The proof of Proposition~\ref{stabilitysub} is now a straightforward consequence of Lemma~\ref{preconvergence} introducing the union of both subdivisions that is finer than each of them.

\begin{proof}[Proof of Proposition~\ref{stabilitysub}.]
Let $\cT_3= \cT_1 \cup \cT_2$ be the union of both subdivisions, so that $\cT_1, \, \cT_2 \subset \cT_3$ and let $V_{\e,\cT_3}(t)$ be the time-discrete approximate mean curvature flow with respect to $\cT_3$ starting from $V_0$. By Lemma~\ref{preconvergence} we infer that $\forall t \in [0,1]$: 
\begin{align*}
  \Delta(V_{\e,\cT_1}(t),V_{\e,\cT_2}(t)) & \leq  \Delta(V_{\e,\cT_1}(t),V_{\e,\cT_3}(t)) +  \Delta(V_{\e,\cT_3}(t),V_{\e,\cT_2}(t)) \\
  & \leq 2c_{7,M} \, t \,  \delta \, \e^{-n-11} \exp(c_{5,M} \, t _, \e^{-n-7})
\end{align*}
which concludes the proof of Proposition \ref{stabilitysub} (up to doubling $c_7>0$).
\end{proof}

We conclude the section with the following corollary. It encompasses the previous results on the stability of the time-discrete approximate mean curvature flow.

\begin{cor}\label{stab_combined}
Let $\e\in(0,1)$, $M\geq 1$. Let $V_0$, $W_0$ be two varifolds in $V_d(\R^n)$ satisfying $\|V_0\|(\R^n)\leq M$ and $\|W_0\|(\R^n)\leq M$. Let $\cT_1=\lbrace t_i \rbrace_{i=1}^m$ and $\cT_2=\lbrace s_j \rbrace_{j=1}^{m'}$ be two subdivisions of $[0,1]$ satisfying \eqref{smallstep}. Let $V_{\e,\cT_1}(t)$ (resp.$W_{\e,\cT_2}(t)$) be the time-discrete approximate mean curvature flow with respect to $\cT_1$ (resp. $\cT_2)$ starting from $V_0$ (resp. $W_0$).\\
If we set: $\delta = \max \{\delta(\cT_1),\delta(\cT_2)\}$, we have for all $t\in [0,1]$,
\begin{equation*}
 \Delta(V_{\e,\cT_1}(t),W_{\e,\cT_2}(t)) \leq  \left[ \Delta(V_0,W_0) + c_{7,M} \, t \, \delta \, \e^{-n-11} \right] \exp(c_{5,M}\, t \, \e^{-n-7}) \: .
\end{equation*}
\end{cor}
\begin{proof}
We start by setting $\cT_3= \cT_1 \cup \cT_2$. Let $V_{\e,\cT_3}(t)$ (resp. $W_{\e,\cT_3}(t)$) be the time-discrete approximate mean curvature flow with respect to $\cT_3$ and starting from $V_0$ (resp. $W_0$). By Lemma \ref{preconvergence} we infer that: for all $t\in[0,1]$,
\begin{equation*}
  \max \{ \Delta(V_{\e,\cT_1}(t),V_{\e,\cT_3}(t)) , \Delta(W_{\e,\cT_2}(t),W_{\e,\cT_3}(t)) \} \leq c_{7,M} \, t \, \delta \, \e^{-n-11} \exp(c_{5,M} \, t \, \e^{-n-7}),
\end{equation*}
and, by Proposition~\ref{stability}, we obtain
 \(\Delta(V_{\e,\cT_3}(t),W_{\e,\cT_3}(t))\leq \Delta(V_0,W_0) \exp( c_{5,M} \, t \, \e^{-n-7}),\)
and we conclude the proof of Corollary~\ref{stab_combined} using the triangular inequality.
\end{proof}
\section{Existence, uniqueness and properties of a limit approximate flow}\label{def_time_step_to_0}

For any given varifold $V_0$ of finite mass, $\e\in(0,1)$ and $\cT$ a subdivision of $[0,1]$, we constructed a time-discrete approximate mean curvature flow, denoted $\left( V_{\e,\cT}(t) \right)_{t\in[0,1]}$ (see Definition \ref{damcf}). The goal of this section is to prove that, as the subdivision is refined—regardless of how the successive finer subdivisions are chosen—there is convergence to a \underline{unique} limit flow (Theorem~\ref{damcfconvergence}) that we call the approximate mean curvature flow starting from $V_0$ and that we denote $\left( V_{\e}(t)\right)_{t\in[0,1]}$. We will exhibit some properties of this limit flow, namely, the stability with respect to initial data (Proposition~\ref{stability1}) and the mass decay (Remark~\ref{remk:massdecrease4}). In Proposition~\ref{prop:timeepsilonbrakke}, we will prove that $V_{\e}(t)$ satisfies a Brakke-type inequality (in reference to inequality \eqref{intro_eq:weak_mcf}) depending on its approximate mean curvature.

\subsection{Existence and uniqueness of a limit approximate flow}

In the following theorem, we state that the time-discrete approximate mean curvature flows starting from a given varifold $V_0$ of compact support converge, as the subdivision time step converges to $0$, to a unique limiting flow. The proof is based on the uniform boundedness of the masses (see \eqref{uniformbound}) and the stability result with respect to the subdivision of Proposition~\ref{stabilitysub}.

\begin{theo}[Convergence]\label{damcfconvergence}
Let $\e \in (0,1)$, $M\geq 1$ and let $V_0 \in V_d(\R^n)$ be a varifold of compact support and satisfying $\|V_0\|(\R^n)\leq M$. For each $j \in \N$,
\begin{itemize}
 \item let $\cT_j^{D} = \{ k \: 2^{-j} \}_{k = 0,1,\ldots, 2^j}$ be the dyadic subdivision of the interval $[0,1]$ of size $\delta (\cT_j^D) = 2^{-j} \xrightarrow[j \to \infty]{} 0$,
 \item let $V_{\e,\cT_j^D}(t)_{t\in[0,1]}$ be the time-discrete approximate mean curvature flow with respect to $\cT_j^D$ starting from $V_0$. Note that according to condition \eqref{smallstep} in Definition~\ref{damcf}, such a flow is well-defined for $j$ large enough so that $c_3 2^{-j} \leq (M+1)^{-3} \e^8$.
\end{itemize}
Then,
\begin{enumerate}[label=(\roman*)]
 \item there exists a family $(V_\e(t))_{t \in [0,1]}$ in $V_d(\R^n)$ such that for any $t \in [0,1]$: 
 \begin{enumerate}[label=\arabic*.]
  \item $\| V_\e(t) \| (\R^n) \leq M+1$,
  \item $V_{\e,\cT_j^D}(t) \xrightharpoonup[]{*} V_\e(t)$,
  \item $\Delta \left( V_{\e,\cT_j^D}(t) , V_\e (t) \right) \rightarrow 0 \quad  \text{as } j \rightarrow +\infty$.
 \end{enumerate}
 \item If $\big( \cT_j \big)_{j \in \N}$ is any other sequence of subdivisions of size $\delta \big(\cT_j \big)$ tending to $0$, then $V_{\e,\cT_j}(t)_{t\in[0,1]}$ (the time-discrete approximate mean curvature flow with respect to $\cT_j$ starting from $V_0$) converges to the same family $(V_\e(t))_{t \in [0,1]}$ as for the dyadic subdivisions: for any $t \in [0,1]$,
 \begin{equation*}
   V_{\e,\cT_j}(t) \xrightharpoonup[]{*} V_\e(t) \quad \text{and} \quad \Delta \left( V_{\e,\cT_j}(t) , V_\e (t) \right) \rightarrow 0 \quad  \text{as } j \rightarrow +\infty \: .
 \end{equation*}
\end{enumerate}
In other words, there exists a unique limit flow $(V_\e(t))_{t \in [0,1]}$ starting from $V_0$, that we call
the approximate mean curvature flow of $V_0$.
\end{theo}
\begin{proof}
Let $\e \in (0,1)$, $M\geq 1$ and let $V_0 \in V_d(\R^n)$ be a varifold with compact support and satisfying $\|V_0\|(\R^n)\leq M$.\\
We start with the proof of $(i)$.
For $j\in \N$, let $\cT_j^D$ and $V_{\e, \cT_j^D} \:$ be as in the statement of the theorem. Let $t \in [0,1]$.
By construction, we know that for any $t\in[0,1]$: 
\begin{equation*}
 \|V_{\e,\cT_j^D}(t)\|(\R^n) \leq M+1.
\end{equation*}
Thus, by Banach-Alaoglu's theorem, there exists a subsequence $a_t(j)$ (depending on $t$) for which the sequence $V_{\e,\cT_{a_t(j)}^D}(t)$  converges weakly-* to a certain limit denoted by $V_\e(t)$. Note that up to this point, such a limit could depend on the extraction $a_t$ and on the specific choice of the dyadic subdivisions $(\cT_j^D)_j$. We first show that the whole sequence $V_{\e,\cT_j^D}(t)$ (and not only the extracted one) converges to $V_\e(t)$ as $j \to \infty$.

\noindent As $V_0$ has compact support, there exists $R_0 > 0$ such that $\supp V_0 \subset B_{R_0}(0) \times \G$. Then, thanks to Remark~\ref{rem:uniformbound}, 
all the varifolds we are considering hereafter are supported in the common bounded set $B \left(0, R_0 + c_1 (M+1) \e^{-2} \right) \times \G$.
Applying Proposition~\ref{metrization1}, we can deduce that,
\begin{equation*}
\Delta(V_{\e}(t),V_{\e,\cT_{a_t(j)}^D}(t))\xrightarrow[j\rightarrow\infty]{} 0 \: .                                                                                                                                                                \end{equation*}
Note that $a_t(j) \geq j$ and therefore the dyadic subdivision $\cT_{a_t(j)}^D$ is finer than $\cT_j^D$. For $j$ large enough so that $c_3 2^{-j} \leq (M+1)^{-3} \e^{8}$, we can apply Lemma~\ref{preconvergence} with $\cT_j^D \subset \cT_{a_t(j)}^D$ and obtain
\begin{equation*}
 \Delta( V_{\e,\cT_{a_t(j)}^D}(t) , V_{\e,\cT_j^D}(t) ) \leq c_{7,M} \, t \,  \delta_j \, \e^{-n-11} \exp(c_{5,M} \, t \, \e^{-n-7}) \quad \text{with } \delta_j = 2^{-j} \: .
\end{equation*}
This implies
\begin{equation*}
 \Delta(V_{\e}(t),V_{\e,\cT_j^D}(t))\leq \Delta(V_{\e}(t),V_{\,\cT_{a_t(j)}^D}(t)) + \Delta(V_{\e,\cT_{a_t(j)}^D}(t),V_{\e,\cT_j^D}(t)) \xrightarrow[j\rightarrow\infty]{} 0.
\end{equation*}
Thus, the full sequence $V_{\e,\cT_j^D}(t)$ converges to $V_{\e}(t)$ for each $t\in[0,1]$ in the bounded Lipschitz topology and thus in the weakly-* topology (again thanks to Proposition~\ref{metrization1}).\\
We now prove $(ii)$. Let $\big( \cT_j \big)_{j \in \N}$ be a sequence of subdivisions of size $\delta \big(\cT_j \big)$ tending to $0$. For $j$ large enough so that $c_3 \delta(\cT_j) \leq (M+1)^{-3} \e^8$, let $V_{\e,\cT_j}(t)_{t\in[0,1]}$ be the time-discrete approximate mean curvature flow with respect to $\cT_j$ starting from $V_0$. Let $t \in [0,1]$ and set $\tilde{\delta}_j=\max \{\delta(\cT_j), \delta(\cT_j^D) \}$; we apply Proposition \ref{stabilitysub} and obtain 
\begin{equation*}
 \Delta (V_{\e,\cT_j^D}(t),V_{\e,\cT_j}(t)) \leq  c_{7,M} \, t \, \tilde{\delta}_j \, \e^{-n-11}\exp(c_{5,M} \, t \, \e^{-n-7}) \xrightarrow[j\rightarrow \infty]{}0.
\end{equation*}
This implies that for any $t\in[0,1]$, $V_{\e,\cT_j}(t)$ converges to $V_{\e}(t)$ both in the bounded Lipschitz topology and thus in the weak-* topology (again thanks to Proposition~\ref{metrization1}).

\noindent We conclude that independently of how the time step goes to $0$, the limit flow exists and is \textit{unique}, we call it the approximate mean curvature flow and we will denote it by $(V_{\e}(t))_t$.
\end{proof}

Given $\e \in (0,1)$ and a subdivision $\cT$ of $[0,1]$, we proposed in Remark~\ref{remk:PWflow} an alternative definition $V_{\e,\cT}^{pc}$ of time-discrete approximate mean curvature flow: we recall that the difference with $V_{\e,\cT}$ lies in the way the flow is extended from the points $t_0, t_1, \ldots, t_m \in \cT$ of the subdivision to any $t \in [0,1]$. While $V_{\e,\cT}$ is defined through a kind of linear interpolation between the flow at time $t_i$ and $t_{i+1}$, $V_{\e,\cT}^{pc}$ is set to be constant in between such subdivision times.
In the following proposition, we derive an error term estimate between both extensions and infer that they lead to the same definition of limit flow $(V_\e(t))_{t \in [0,1]}$.

\begin{prop}\label{coincidence}
 Let $\e\in(0,1)$. Let $\cT=\lbrace t_i \rbrace_{i=0}^m$ be a subdivision of $[0,1]$ satisfying \eqref{smallstep}. Let $V_0 \in V_d(\R^n)$ of compact support and  $V_{\e,\cT}(t)$ the time-discrete approximate mean curvature flow with respect to $\cT$ starting from $V_0$. Let $V_{\e,\cT}^{pc}(t)$ be the associated piecewise constant flow with respect to $\cT$ (Remark~\ref{remk:PWflow}).
 Then,
 \begin{equation*}
  \Delta(V_{\e,\cT}(t),V_{\e,\cT}^{pc}(t)) \leq c_8\, (M+1)^2 \, \delta(\cT) \, \e^{-4},\quad \forall t \in [0,1],
 \end{equation*}
 where $c_8>0$ is a constant that depends only on $n$. As a consequence, when the step of the subdivision goes to $0$, $V_{\e,\cT}^{pc}(t)$ converges to $V_{\e}(t)$ (defined in Theorem \ref{damcfconvergence}):
 for any $t \in [0,1]$,
 \begin{equation} \label{eqPWcv}
   V_{\e,\cT_j}^{pc}(t) \text{ converges weakly-* to } V_\e(t) \quad \text{and} \quad \Delta \left( V_{\e,\cT_j}^{pc}(t) , V_\e (t) \right) \rightarrow 0 \quad  \text{as } j \rightarrow +\infty \: .
 \end{equation}
\end{prop}
\begin{proof}
As $\e$ and $\cT$ are fixed, we denote $V_{\e,\cT}(t)$ and $V_{\e,\cT}^{pc}(t)$ by $V(t)$ and $V(t)^{pc}$ throughout the proof. Let $i\in\lbrace 0, \dots,m-1\rbrace$ be such that $t \in [t_{i},t_{i+1})$. Introducing $f={\rm id}+(t-t_{i})h_{\e}(\cdot,V(t_{i}))$ and  $g={\rm id}$, we have
\begin{equation} \label{eq:VtVtpc}
\Delta(V(t),V^{pc}(t)) =  \Delta(V(t),V(t_i))
= \Delta \Bigl( f_{\#}V(t_{i}),g_{\#}V(t_{i})) \Bigr).
\end{equation}
Using \eqref{hepsilonbound2} we can check that (noting that $\|V(t_i)\|\leq M+1$)
\begin{equation*}
\| Df - Dg \|_{\infty} = (t-t_i) \|  Dh_{\e}(\cdot,V(t_i)) \|_{\infty} \leq 2c_1(t-t_i)(M+1) \e^{-4} \leq 2c_1 \, \delta(\cT) (M+1) \e^{-4} \leq \frac12
\end{equation*}
and we then can apply \eqref{damcfstability5} with $V=V(t_{i})$ so that
\begin{equation} \label{eq:VtVtpc-2}
 \Delta \Bigl(f_{\#}V(t_{i}), g_{\#} V(t_{i})  \Bigr) 
 \leq   (M+1) \left( c_5 (t-t_{i})\|h_{\e}(\cdot,V(t_{i}))\|_{C^1} + \| J_{\cdot}f-1 \|_{\infty} \right).
\end{equation}
We conclude the proof of Proposition~\ref{coincidence} thanks to \eqref{eq:VtVtpc}, \eqref{eq:VtVtpc-2}, Proposition~\ref{hepsilonbound} and \eqref{eq:JSfbounds}. 
\end{proof}

Thereafter $(V_{\e}(t))_{t \in [0,1]}$ denotes the approximate mean curvature flow starting from $V_0$ as defined in Theorem~\ref{damcfconvergence}. We now investigate the properties of this flow, starting with the stability with respect to the initial varifold.
\begin{prop}\label{stability1}
 Let $\e\in(0,1), M \geq 1$. Let $V_0$, $W_0$ be two varifolds in $V_d(\R^n)$ satisfying $\|V_0\|(\R^n) \leq M$ and $\|W_0\|(\R^n) \leq M$, both compactly supported. Then, for all $t \in [0,1]$,
 \begin{equation*}
\Delta(V_{\e}(t),W_{\e}(t)) \leq  \Delta(V_0,W_0) \exp(c_{5,M} \, t \, \e^{-n-7}) \: ,
\end{equation*}
where $V_\e$ (resp. $W_\e$) is the approximate mean curvature flow starting from $V_0$ (resp. $W_0$).
\end{prop}
\begin{proof}
We fix a sequence of subdivisions $\left( \cT_j \right)_{j \in \N}$ with time step $\delta_j \rightarrow 0$ as $j\rightarrow \infty$ (we can take the dyadic subdivisions for instance). Let $\left( V_{\e, \cT_j}(t) \right)_{t \in [0,1]}$ (resp. $\left( W_{\e,\cT_j}(t) \right)_{t \in [0,1]}$) be the time-discrete approximate mean curvature flow with respect to $\cT_j$ starting from $V_0$ (resp. $W_0$). 
Let $j$ be large enough so that \eqref{smallstep} holds: $c_3 \delta_j \leq (M+1)^{-3} \e^8$, then we can simply
apply Proposition~\ref{stability} and obtain for all $t \in [0,1]$,
\begin{equation*}
  \Delta(V_{\e, \cT_j}(t),W_{\e, \cT_j}(t)) \leq  \Delta(V_0,W_0) \exp(c_{5,M} \: t \: \e^{-n-7}) \: ,
\end{equation*}
and therefore, by the triangle inequality and Theorem~\ref{damcfconvergence}, letting $j$ tend to $\infty$, we can conclude that
\begin{align*}
  \Delta(V_{\e}(t),W_{\e}(t)) &
\leq \Delta(V_{\e}(t),V_{\e, \cT_j}(t)) + \Delta(V_{\e, \cT_j}(t),W_{\e, \cT_j}(t))+\Delta(W_{\e, \cT_j}(t),W_{\e}(t)) \\
& \leq \Delta(V_0,W_0) \exp(c_{5,M} \, t \, \e^{-n-7}) \: .
\end{align*}
\end{proof}

\begin{remk} \label{remkPointCloudApprox}
Note that if we reformulate Proposition~\ref{stability1} in the case where a compactly supported varifold $V_0$ is approximated by a sequence of varifolds $W_k \xrightharpoonup[k \rightarrow \infty]{\ast} V_0$, $W_k$ being successive discretizations of $V_0$ for instance. Then, considering 
\begin{itemize}
 \item $(V_{\e}(t))_{t\in[0,1]}$, the approximate mean curvature flow starting from $V_0$,
 \item $\forall k$, $((W_k)_{\e}(t))_{t\in[0,1]}$, the approximate mean curvature flow starting from $W_k$,
\end{itemize}
we have that
\begin{equation*}
\Delta(V_{\e}(t),(W_k)_{\e}(t)) \leq  \Delta(V_0,W_k)\exp(t c_{5,M} \e^{-n-7}) \xrightarrow[k \to \infty]{} 0 \: ,
\end{equation*}
provided that $\supp W_k$ are contained in a common compact set and $\sup \| W_k \|(\R^n) < \infty$.
In order to state some convergence result, one can choose $\e_k \to 0$ such that $\Delta(V_0,W_k)\exp(t c_{5,M} \e^{-n-7}) \rightarrow 0$, however, we would need a convergence property of $V_\e (t)$ as $\e \to 0$ to some limit flow $V(t)$. This issue is addressed in Section~\ref{sec:spacetime-Brakke-flows}.
\end{remk}

\subsection{Equality à la Brakke}
In Proposition~\ref{prop:timeepsilonbrakke}, we show a Brakke-type equality for the approximate flow $\left( V_{\e}(t) \right)_{t\in[0,1]}$ with respect to its approximate mean curvature. 
The proof consists of taking the limit in inequality \eqref{timeepsilonbrakke0} which results from the expansion of the push-forward varifold formula \eqref{appbrakke} and Theorem \ref{damcfconvergence}. We conclude the section with the decay property of mass $t \mapsto \| V_{\e}(t)\|(\R^n)$, which follows directly from \eqref{timeepsilonbrakke} and \eqref{massdecrease1}.

\noindent We first introduce the following lemma on the regularity of the weighted first variation with respect to the varifold.
\begin{lemma}\label{lem:wfirstvar}
Let $\phi \in \xC^2(\R^n, \R_+)$, $X \in \xC^2(\R^n, \R^n)$ and let $V$, $W \in V_d(\R^n)$ be two varifolds of finite mass. Then
\begin{equation*}
 \left| \delta (V,\phi)(X) - \delta (W,\phi)(X)  \right| \leq n \| \phi \|_{\xC^2} \| X \|_{\xC^2} \Delta(V,W)\: .
\end{equation*}
\end{lemma}
\begin{proof}
Let $\phi \in \xC^1(\R^n, \R_+)$, $X \in \xC^2(\R^n, \R^n)$ and set $\Theta(x,S) := \phi(x) \mdiv_SX(x) + \nabla\phi(x)\cdot X(x)$, we recall that $\mdiv_SX = \tr(S\circ DX)$.
From Definition \eqref{wfirstvar} one has
\begin{equation} \label{eq:deltaVWphi}
 \left| \delta (V,\phi)(X) - \delta (W,\phi)(X)  \right| \leq \max \{ \| \Theta \|_\infty, \lip(\Theta) \} \Delta(V,W).
\end{equation}
First note that for $A, B\in\mathcal{M}_n$, one has
\begin{equation*}
|\tr (A \circ B)| \leq n | A \circ B |_\infty \leq n \| A \circ B \| \leq \| A \| \| B \| \: .
\end{equation*}
In particular, for $x \in \R^n$, $S\in\G$, one has 
\begin{equation} \label{eq:trAB}
|\mdiv_S X(x)| = |\tr(S\circ DX(x)) |\leq n\| DX(x) \|.
\end{equation}
and consequently,
\begin{equation} \label{eq:ThetaSup}
 \| \Theta \|_\infty \leq n \| \phi \|_\infty \| D X \|_\infty + \| \nabla \phi \|_\infty \| X \|_\infty \leq n \| X \|_{\xC^1}  \| \phi \|_{\xC^1} \: .
\end{equation}
We are left with estimating $\lip(\Theta)$.

\noindent Let $x,y \in \R^n$ and $(S,T) \in \G$, thanks to \eqref{eq:trAB}, we have 
\begin{align*}
&   | \phi(x) \mdiv_SX(x)  - \phi(y) \mdiv_TX(y) | \\
 & \leq  | \phi(x) - \phi(y)| | \mdiv_SX(x) |  
 + | \phi(y)|\, | \mdiv_SX(x) - \mdiv_SX(y)|  + | \phi(y)| |  \mdiv_SX(y) -  \mdiv_TX(y) | \\
  \leq & n \| \nabla \phi \|_\infty \| DX \|_\infty |x-y| + n \| \phi \|_\infty \| DX(x) - DX(y) \| + n \| \phi \|_\infty \| DX \|_\infty \| S - T \| \\
  \leq & n \left( \| \nabla \phi \|_\infty \| DX \|_\infty + \| \phi \|_\infty \| D^2 X \|_\infty \right) |x-y| + n \| \phi \|_\infty \| DX \|_\infty \| S -T \| \: .
 \end{align*}
Furthermore,
$
 \lip(\nabla\phi\cdot X) \leq \| \nabla^2 \phi \|_{\infty} \|X\|_{\infty} + \|\nabla\phi \|_{\infty}\|DX\|_{\infty}$.
Therefore,
\begin{equation} \label{eq:ThetaLip}
 | \Theta(x, S) - \Theta(t,y) | \leq n \|  \phi \|_{\xC^2} \| X \|_{\xC^2} \: ;
\end{equation}
We conclude the proof of Lemma~\ref{lem:wfirstvar} thanks to \eqref{eq:deltaVWphi}, \eqref{eq:ThetaSup} and \eqref{eq:ThetaLip}. 
\end{proof}
\begin{prop} \label{prop:timeepsilonbrakke}
Let $\e \in(0,1)$ and $M\geq 1$. Let $V_0$ in $V_d(\R^n)$ be a varifold with compact support such that $\|V_0\|(\R^n)\leq M$. 
\noindent Let $(V_{\e}(t))_{t \in [0,1]}$ be the approximate mean curvature flow starting from $V_0$.  For any $\phi \in \xC^1([0,1] \times \R^n,\R_+)$  and $0 \leq a \leq b \leq 1$ we have 
\begin{equation}\label{timeepsilonbrakke}
    \|V_{\e}(b)\|(\phi(\cdot,b)) - \|V_{\e}(a)\|(\phi(\cdot,a)) = \int_{a}^{b} \delta (V_{\e}(t),\phi(\cdot,t)) (h_{\e}(\cdot,V_{\e}(t))) \: dt +\int_{a}^{b} \|V_{\e}(t)\|(\partial_t \phi(\cdot,t)) \: dt.
\end{equation}

\end{prop}
\begin{proof}
Let $\e \in(0,1)$ and $M\geq 1$. Let $V_0$ in $V_d(\R^n)$ be a varifold with compact support such that $\|V_0\|(\R^n)\leq M$. 
Let $\cT=\lbrace t_i \rbrace_{i=1}^m$ be a  uniform subdivision of $[0,1]$ of time step $\del = \frac{1}{m}$ satisfying \eqref{smallstep}. Let $(V^{pc}_{\e,\cT}(t))_{t\in[0,1]}$ be the piecewise constant approximate mean curvature flow with respect to to $\cT$ starting from $V_0$.
We first prove that \eqref{timeepsilonbrakke} holds for $(V^{pc}_{\e,\cT}(t))_{t\in[0,1]}$ up to an error term of order $\Delta t$. More precisely, we prove in Steps~$1$ and $2$ that $\exists C > 0$ (only depending on $n$ and $M$) such that for any $\phi \in \xC^2([0,1] \times \R^n,\R_+)$ and $0 \leq a \leq b \leq 1$ we have 
\begin{multline}\label{timeepsilonbrakke0}
 \left| \|V^{pc}_{\e,\cT}(b)\|(\phi(\cdot, b)) - \|V^{pc}_{\e,\cT}(a)\|(\phi(\cdot, a)) - \int_{a}^{b} \delta (V^{pc}_{\e,\cT}(t),\phi ( \cdot, t) )(h_{\e}(\cdot, V^{pc}_{\e,\cT}(t))) \: dt \right. \\ \left. - \int_{a}^{b} \int_{\R^n}  \partial_t \phi (\cdot, t)  \, d \| V^{pc}_{\e,\cT}(t) \| dt \right| \leq C \|\phi\|_{\xC^2} \del \e^{-8} \: .
\end{multline}
In Step~$3$, recalling that $(V^{pc}_{\e,\cT}(t))_{t\in[0,1]}$ converges to $(V_\e(t))_{t\in[0,1]}$ when considering subdivisions $\cT$ whose size tends to $0$, we take the limit in \eqref{timeepsilonbrakke0} and establish \eqref{timeepsilonbrakke} for $\phi$ of regularity $\xC^2$. We conclude the proof of Proposition~\ref{prop:timeepsilonbrakke} by density of $\xC^2 ([0,1] \times \R^n, \R_+)$ in $\xC^1 ([0,1] \times \R^n, \R_+)$ in Step 4. Remark \ref{rem:uniformbound} states that $\displaystyle \bigcup_{t\in[0,1]} \supp V_{\e,\cT}^{pc}(t)$ is contained in a compact set $K = K_{\e}$ that does not depend on the subdivision $\cT$, hence $\displaystyle\bigcup_{t\in[0,1]} \supp V_{\e}(t)$ is contained in $K$ as well. In the proof, the $\xC^k$-norms ($k\in\lbrace 1,2 \rbrace$) of the test functions $\psi$ and $\phi$ are implicitly taken with respect to the set $K$ (in particular they are finite).
%
In both Steps~$1$ and $2$, $\cT$ and $\e$ are fixed and we denote for simplicity $V(t) := V^{pc}_{\e,\cT}(t)$. Throughout the proof, $C > 0$ is a generic constant that depends only on $n$ and $M$ and may change from line to line.

\noindent {\bf Step 1:} We prove the inequality \eqref{timeepsilonbrakke0} for $a,b \in \cT$: 
let $\ell \in \{ 0, 1, \ldots, m-1 \}$ and $\phi \in \xC^2([0,1] \times \R^n, \R_+)$.\\
We can apply \eqref{appbrakke} in Proposition~\ref{mcfmotion1} to a spatial test function $\psi \in \xC^2(\R^n, \R_+)$.
We recall that the piecewise constant mean curvature flow coincides with the time discrete approximate mean curvature flow at the points of the subdivision and furthermore, $f_\# V(t_\ell) = V(t_{\ell + 1})$ for $f = {\rm id} + \del h_\e (\cdot, V(t_\ell))$ as in \eqref{appbrakke}. Therefore,
\begin{align*}
\Big| \|V( t_{\ell + 1} )\|(\psi) - \|V(t_\ell)\|(\psi)- \del \,  \: \delta (V(t_\ell),\psi)(h_{\e}(\cdot,V(t_\ell))) \Big|
 & \leq c_3 (M+1)^3 \|\psi\|_{\xC^2}(\del)^2 \e^{-8} \\
 &\leq C \|\psi\|_{\xC^2}(\del)^2 \e^{-8} \: .
\end{align*}
We now recall that $V(t)$ is piecewise constant and thus, for all $t \in (t_\ell, t_{\ell + 1})$, $V(t)$ = $V(t_\ell)$ and
\begin{equation*}
\int_{t_\ell}^{t_{\ell+1}} \delta (V(t),\psi)(h_{\e}(\cdot,V(t))) \: dt =  \del \, \delta(V(t_\ell),\psi)(h_{\e}(\cdot,V(t_\ell)))
\end{equation*}
so that taking $\psi = \phi( \cdot, t_{\ell +1})$, we obtain 
\begin{multline} \label{eq:BrakkeEpsStep1:-1}
\left| \|V( t_{\ell + 1} )\|(\phi( \cdot, t_{\ell +1})) - \|V(t_\ell)\|(\phi( \cdot, t_{\ell +1}))- \int_{t_\ell}^{t_{\ell+1}} \delta (V(t),\phi( \cdot, t_{\ell +1}))(h_{\e}(\cdot,V(t))) \: dt \right|
\\ \leq C \|\phi( \cdot, t_{\ell +1})\|_{\xC^2} (\del)^2 \e^{-8} \leq C \|\phi \|_{\xC^2} (\del)^2 \e^{-8} \: .
\end{multline}
Applying the mean value theorem to $\phi$ and $\nabla \phi$:
\[
\int_{t_\ell}^{t_{\ell + 1}} \| \phi ( \cdot, t) - \phi(\cdot, t_{\ell +1}) \|_{\xC^1}  \: dt \leq \int_{t_\ell}^{t_{\ell + 1}}  \| \phi \|_{\xC^2} |t - t_{\ell+1} | \: dt \leq \| \phi \|_{\xC^2} (\del)^2  \: ,
\]
and then using Remark~\ref{remk:wfirstvar} and  Proposition~\ref{hepsilonbound}, we have
\begin{align}
 & \Big|  \int_{t_\ell}^{t_{\ell+1}} \delta (V(t),  \phi ( \cdot, t_{\ell + 1}) )(h_{\e}(\cdot,V(t))) \: dt - \int_{t_\ell}^{t_{\ell+1}} \delta (V(t),\phi ( \cdot, t) )(h_{\e}(\cdot,V(t))) \: dt \Big| \nonumber \\
 \leq &  \int_{t_\ell}^{t_{\ell + 1}} \left| \delta (V(t),\phi ( \cdot, t) - \phi(\cdot, t_\ell))(h_{\e}(\cdot,V(t))) \right| \: dt \nonumber \\
  \leq & \int_{t_\ell}^{t_{\ell + 1}} n M \| h_\e (\cdot, V(t)) \|_{\xC^1} \| \phi ( \cdot, t) - \phi(\cdot, t_\ell) \|_{\xC^1} \:  dt \nonumber \\
  \leq & C \|\phi \|_{\xC^2}(\del)^2 \e^{-4} \: . \label{eq:BrakkeEpsStep1:0}
\end{align}
Writing for $x \in \R^n$, $\displaystyle \phi(x, t_{\ell+1}) - \phi(t,x_\ell) = \int_{t_\ell}^{t_{\ell +1}} \partial_t \phi(t,x) \: dt$ and integrating with respect to $\| V (t_\ell) \| = \| V(t) \|$ for all $t \in (t_\ell, t_{\ell + 1})$,
\begin{multline} \label{eq:BrakkeEpsStep1:00}
 \| V (t_{\ell }) \| ( \phi(\cdot, t_{\ell + 1}) )  - \| V (t_{\ell }) \| ( \phi(\cdot, t_{\ell}) )  = \int_{x \in \R^n} \int_{t_\ell}^{t_{\ell + 1}} \partial_t \phi (x ,t) \: dt \: d \| V(t_{\ell}) \|(x) \\
= \int_{t_\ell}^{t_{\ell + 1}} \int_{x \in \R^n} \partial_t \phi (x ,t)  \: d \| V(t_{\ell}) \|(x) \: dt = \int_{t_\ell}^{t_{\ell + 1} } \partial_t \phi (\cdot , t) \: d \| V (t) \| \: dt 
\: .
\end{multline}
We can now combine \eqref{eq:BrakkeEpsStep1:-1}, \eqref{eq:BrakkeEpsStep1:0} and \eqref{eq:BrakkeEpsStep1:00} to obtain that for $\ell \in \{0, 1, \ldots, m - 1 \}$,
\begin{align*}
 & \left| \|V(t_{\ell +1 })\|(\phi(\cdot, t_{\ell + 1})) - \|V(t_\ell)\|(\phi(\cdot, t_\ell)) - \int_{t_\ell}^{t_{\ell +1}} \delta (V(t),\phi ( \cdot, t) )(h_{\e}(\cdot,V(t))) \: dt \right. \\
 & \qquad \left. - \int_{t_\ell}^{t_{\ell +1}} \| V(t) \| ( \partial_t \phi (\cdot, t) ) \: dt \right| \\
= & \left| \|V(t_{\ell + 1})\|(\phi(\cdot, t_{\ell+1}))  - \|V(t_\ell)\|(\phi(\cdot, t_{\ell +1})) -  \int_{t_\ell}^{t_{\ell+1}} \delta (V(t),\phi ( \cdot, t) )(h_{\e}(\cdot,V(t))) \: dt \right|   \\
 & \leq C \|\phi\|_{\xC^2} (\Delta t)^2 \: \e^{-8} \: ,
\end{align*}
which implies \eqref{timeepsilonbrakke0} for $a = t_\ell$ and $b = t_{\ell + 1}$ since $\del \leq 1$. Let now $a = t_p \leq t_q = b$,
note that if $p = q$, \eqref{timeepsilonbrakke0} is trivial and otherwise, summing up the previous estimates for $\ell \in \left\lbrace p, \dots , q -1 \right\rbrace$ and using $(q-p) \del = t_q - t_p \leq 1$ leads to the inequality \eqref{timeepsilonbrakke0}, which concludes the proof of Step 1 (case where $a$, $b \in \cT$).
%

\noindent {\bf Step 2:}
We now recover the approximate Brakke-type equality \eqref{timeepsilonbrakke0} for any arbitrary $a,b$. Let $0 \leq a < b \leq 1$ and $\phi \in \xC^2([0,1] \times \R^n,\R_+)$. Let the points $t_p$, $t_q \in \cT$ be such that $t_{p} \leq a < t_{p+1}$ and  $t_{q} \leq b < t_{q+1}$. We then have $|a - t_p| < \del$ and $| t_q - b| < \del $, and recalling that $V (t)$ is piecewise constant on intervals of the form $[t_{\ell}, t_{\ell+1})$, we also have $V(a) = V(t_p)$ and $V(b) = V(t_q)$ so that
\begin{align} \label{eq:BrakkeEpsStep2:1}
 \Big| \|V(b)\|(\phi(\cdot,b)) & - \|V(a)\|(\phi(\cdot, a)) - \|V(t_{q})\|(\phi(\cdot,t_q)) + \|V(t_p)\|(\phi(\cdot, t_p)) \Big| \nonumber \\
 & = \Big| \|V(t_q)\|(\phi(\cdot,b) - \phi(\cdot,t_q)) - \|V(t_p)\|(\phi(\cdot, a) - \phi(\cdot, t_p))  \Big| \nonumber \\
 & \leq 2 (M+1) \| \phi \|_{\xC^1} \del
\end{align}
thanks to the mean value theorem applied to $\phi$.\\
Furthermore, for all $t \in [0,1]$, using Remark~\ref{remk:wfirstvar} and Proposition~\ref{hepsilonbound}, we have
\begin{align*}
 \left| \delta (V(t),\phi(\cdot,t) )(h_{\e}(\cdot,V(t))) \right| \leq n \| h_\e (\cdot,V(t))) \|_{\xC^1} (M+1) \| \phi (\cdot, t)\|_{\xC^1} \leq n c_1 (M+1)^2 \| \phi \|_{\xC^1} \: \e^{-4}
\end{align*}
and therefore
\begin{multline} \label{eq:BrakkeEpsStep2:2}
 \left| \int_a^b \delta (V(t),\phi(\cdot,t) )(h_{\e}(\cdot,V(t)))  \: dt - \int_{t_p}^{t_q} \delta (V(t),\phi(\cdot,t) )(h_{\e}(\cdot,V(t)))  \: dt \right| \\ \leq ( \underbrace{|t_p - a| + |t_q -b|}_{\leq 2 \del} ) \sup_{t \in [0,1]} \left| \delta (V(t),\phi)(h_{\e}(\cdot,V(t))) \right| \leq C \| \phi \|_{\xC^1} \del \: \e^{-4} \: .
\end{multline}
We are left with estimating
\begin{align} \label{eq:BrakkeEpsStep2:3}
 \left| \int_a^b \| V(t) \| (\partial_t \phi(\cdot,t) ) \: dt - \int_{t_p}^{t_q} \| V(t) \| ( \partial_t \phi(\cdot,t) ) \: dt \right| \leq 2 (M+1)  \| \phi \|_{\xC^1} \del \: .
\end{align}
We can complete the proof of Step 2 and establish \eqref{timeepsilonbrakke0} by combining \eqref{eq:BrakkeEpsStep2:1} and \eqref{eq:BrakkeEpsStep2:2} and \eqref{eq:BrakkeEpsStep2:3} together with Step 1.

\noindent{\bf Step 3:} We show \eqref{timeepsilonbrakke} restricted to $\xC^2$ test functions.

\noindent We first recall that the approximate mean curvature flow $(V_\e(t))_{t\in[0,1]}$ starting from $V_0$ can be obtained as the limit flow ($j \to \infty$) of any time-discrete approximate mean curvature flow $(V_{\e, \cT_j}(t))_t$ for subdivisions $\cT_j$ of size $\delta(\cT_j)$ tending to $0$, as stated in Theorem~\ref{damcfconvergence}. We can thus choose a sequence of uniform subdivisions $\cT_j$ of size $\Delta t_j := \delta(\cT_j) \xrightarrow[j \to \infty]{} 0$, and we fix the subdivisions $\cT_j = \{t_{\ell,j}\}_{\ell = 0, \ldots, m_j}$ hereafter, we will write $t_\ell$ instead of $t_{\ell,j}$ in the proof in order to lighten notations. We additionally recall that the piecewise constant flow $(V_{\e, \cT_j}^{pc} (t))_t$ introduced in Remark~\ref{remk:PWflow} converges as well to $(V_{\e}(t))_t$ thanks to Proposition~\ref{coincidence} and as it is more convenient in this proof, we work with $(V_{\e, \cT_j}^{pc} (t))_t$. For the sake of lightening the notation we will denote $V_j(t):=V_{\e, \cT_j}^{pc}(t)$ hereafter.

\noindent We carry on with the proof of Step 3, let $\phi \in \xC^2 ([0,1] \times \R^n, \R_+)$ and $0 \leq a < b \leq 1$, from \eqref{timeepsilonbrakke0} we have for any $j\in\N$, 
\begin{multline}\label{timeepsilonbrakke2}
 \left| \|V_j(b)\|(\phi(\cdot, b)) - \|V_j(a)\|(\phi(\cdot, a)) - \int_{a}^{b} \delta (V_j(t),\phi ( \cdot, t) )(h_{\e}(\cdot,V_j(t))) \: dt \right. \\ \left. - \int_{a}^{b} \| V_j(t) \| ( \partial_t \phi (\cdot, t) ) \: dt \right| \leq C \|\phi\|_{\xC^2} \del_j \: \e^{-8} \: .
\end{multline}
When $j$ goes to $\infty$ we know thanks to \eqref{eqPWcv} in Proposition~\ref{coincidence} that for all $t \in [0,1]$,
\begin{equation*}
 V_{j}(t) \text{ converges weakly-* to } V_\e(t) \quad \text{and} \quad \Delta(V_j(t),V_{\e}(t)) \xrightarrow[j\rightarrow\infty]{} 0 \: .
\end{equation*}
As a first consequence, we obtain that 
\begin{equation} \label{eq:BrakkeEpsStep3:1}
 \| V_j(b) \| (\phi (\cdot,b)) \xrightarrow[j \to \infty]{} \| V_\e(b) \| (\phi (\cdot,b)) \quad \text{and} \quad \| V_j(a) \| (\phi (\cdot,a)) \xrightarrow[j \to \infty]{} \| V_\e(a) \| (\phi (\cdot,a)) \: .
\end{equation}
We recall that $\| V_\e(t) \|(\R^n) \leq M+1$ and thanks to Proposition~\ref{hepsilonbound} and Lemma \ref{hepsilonstability}, we have for all $t \in [0,1]$,
\begin{equation*}
\left\| h_{\e}(\cdot,V_\e(t)) \right\|_{\xC^2} \leq 3 c_1 (M+1) \e^{-6} \text{ and }
 \|h_{\e}(\cdot,V_j(t))-h_{\e}(\cdot,V_{\e}(t))\|_{\xC^1}
 \leq 2c_4 (M+1) \e^{-n-7} \Delta(V_j(t),V_{\e}(t))
 \end{equation*}
and therefore, applying Remark~\ref{remk:wfirstvar}, Lemma~\ref{lem:wfirstvar} and Lemma~\ref{hepsilonstability} (with $V = V_j(t)$ and $W = V_\e(t)$) we infer
\begin{align*}
& \Big| \delta (V_{j}(t),\phi(\cdot,t))(h_{\e}(\cdot,V_{j}(t))) -  \delta (V_{\e}(t),\phi(\cdot,t))(h_{\e}(\cdot,V_{\e}(t))) \Big| \\
\leq & n (M+1) \| \phi \|_{\xC^1} \left\| h_{\e}(\cdot,V_j(t))-h_{\e}(\cdot,V_{\e}(t)) \right\|_{\xC^1} + n \| \phi \|_{\xC^1} \left\| h_{\e}(\cdot,V_{\e}(t))  \right\|_{\xC^2} \Delta(V_j(t),V_{\e}(t)) \\
\leq &  C \| \phi \|_{\xC^1} \: \e^{-n-7} \: \Delta(V_j(t),V_{\e}(t)) \: .
\end{align*}
Integrating the previous inequality, we obtain by dominated convergence, noting that for all $t$, $\Delta(V_j(t),V_{\e}(t)) \leq 2(M+1)$:
\begin{equation} \label{eq:BrakkeEpsStep3:2}
 \left| \int_{a}^b \delta (V_{j}(t),\phi(\cdot,t))(h_{\e}(\cdot,V_{j}(t))) \: dt - \int_a^b \delta (V_{\e}(t),\phi(\cdot,t))(h_{\e}(\cdot,V_{\e}(t))) \: dt \right|  \xrightarrow[j \to \infty]{} 0 \:. 
\end{equation}


\noindent It remains to let $j$ tend to $\infty$ in the following term
\begin{equation} \label{eq:BrakkeEpsStep3:3}
 \lim_{j \to \infty} \int_a^b \| V_j(t) \| (\partial_t \phi(\cdot,t) ) \: dt  = \int_a^b \| V_\e(t) \| (\partial_t \phi(\cdot,t) ) \: dt  
\end{equation}
where the convergence holds by dominated convergence: for each $t \in [0,1]$ the weakly-* convergence of $V_j(t)$ to $V_\e(t)$ implies that $\lim_{j \to \infty} \| V_j(t) \| (\partial_t \phi(\cdot,t) ) = \| V_\e(t) \| (\partial_t \phi(\cdot,t) )$ and the integrands are uniformly bounded by the constant $\| \phi \|_{\xC^1} (M+1)$.\\
We can eventually let $j$ tend to $+\infty$ in \eqref{timeepsilonbrakke2} and conclude the proof of Step 3 (i.e. \eqref{timeepsilonbrakke} for all $\xC^2$ test function $\phi$) combining the convergence of the $3$ terms given by \eqref{eq:BrakkeEpsStep3:1}, \eqref{eq:BrakkeEpsStep3:2} and \eqref{eq:BrakkeEpsStep3:3} in the l.h.s. while the r.h.s. tends to $0$.

\noindent {\bf Step 4:} It remains to check that we can pass from $\xC^2$ to $\xC^1$ test functions $\phi$ in \eqref{timeepsilonbrakke} to conclude the proof of Proposition~\ref{prop:timeepsilonbrakke}.\\
Let $\phi \in \xC^1([0,1] \times \R^n, \R_+)$ and apply density of $\xC^2 ([0,1] \times \R^n, \R_+)$ in $\xC^1 ([0,1] \times \R^n, \R_+)$ to have a sequence of functions $(\phi_k)_{k \in \N}$ such that for all $k$, $\phi_k \in \xC^2 ([0,1] \times \R^n, \R_+)$ and
\begin{equation*}
 \| \phi - \phi_k \|_{\xC^1} \xrightarrow[k \to \infty]{} 0 \: .
\end{equation*}
Thanks to Step 3, we know that \eqref{timeepsilonbrakke} holds for $\xC^2$ functions and thus for all $k \in \N$,
\begin{equation} \label{eq:BrakkeEpsStep4:0}
 \|V_{\e}(b)\|(\phi_k (\cdot,b)) - \|V_{\e}(a)\|(\phi_k (\cdot,a)) = \int_{a}^{b} \delta (V_{\e}(t),\phi_k(\cdot,t)) (h_{\e}(\cdot,V_{\e}(t))) \: dt +\int_{a}^{b} \|V_{\e}(t)\|(\partial_t \phi_k(\cdot,t)) \: dt
\end{equation}
and it remains to check that we can let $k$ tend to $+\infty$, which basically follows from the fact that each term involved in \eqref{timeepsilonbrakke} is linearly continuous with respect to $\phi \in \xC^1 ([0,1] \times \R^n, \R_+)$ endowed with $\| \cdot \|_{\xC^1}$. Indeed, we first note that for any $t \in [0,1]$,
\begin{equation*}
\Big| \|V_{\e}(t)\|(\phi_k (\cdot,t)) -  \|V_{\e}(t)\|(\phi (\cdot,t)) \Big| \leq (M+1) \| \phi_k (\cdot,t) - \phi (\cdot,t) \|_\infty \leq (M+1) \| \phi - \phi_k \|_{\xC^1} 
\end{equation*}
and consequently,
\begin{equation} \label{eq:BrakkeEpsStep4:1}
\lim_{k \to \infty}  \|V_{\e}(b)\|(\phi_k (\cdot,b)) =  \|V_{\e}(b)\|(\phi (\cdot,b)) \quad \text{and} \quad \lim_{k \to \infty}  \|V_{\e}(a)\|(\phi_k (\cdot,a)) =  \|V_{\e}(a)\|(\phi (\cdot,a)) \: .
\end{equation}
Similarly, for all $t \in [0,1]$,
\begin{equation*}
\Big| \|V_{\e}(t)\|(\partial_t \phi_k (\cdot,t)) -  \|V_{\e}(t)\|(\partial_t \phi (\cdot,t)) \Big| \leq (M+1) \| \partial_t \phi_k (\cdot,t) - \partial_t \phi (\cdot,t) \|_\infty \leq (M+1) \| \phi - \phi_k \|_{\xC^1} 
\end{equation*}
and consequently,
\begin{equation} \label{eq:BrakkeEpsStep4:2}
\lim_{k \to \infty} \int_a^b \|V_{\e}(t)\|(\partial_t \phi_k (\cdot,t)) \: dt  =  \int_a^b \|V_{\e}(t)\|(\partial_t \phi (\cdot,t)) \: dt \: .
\end{equation}
We can apply Remark~\ref{remk:wfirstvar} and Proposition~\ref{hepsilonbound} to the remaining term:
\begin{multline*}
\Big| \delta (V_{\e}(t),\phi_k(\cdot,t)) (h_{\e}(\cdot,V_{\e}(t))) - \delta (V_{\e}(t),\phi(\cdot,t)) (h_{\e}(\cdot,V_{\e}(t))) \Big| \\ \leq n (M+1) \left\| h_{\e}(\cdot,V_{\e}(t)) \right\|_{\xC^1}  \| \phi_k (\cdot,t) - \phi (\cdot,t) \|_{\xC^1} \leq C \e^{-4}  \| \phi - \phi_k \|_{\xC^1} 
\end{multline*}
and consequently,
\begin{equation} \label{eq:BrakkeEpsStep4:3}
\lim_{k \to \infty} \int_a^b \delta (V_{\e}(t),\phi_k(\cdot,t)) (h_{\e}(\cdot,V_{\e}(t))) \: dt  =  \int_a^b \delta (V_{\e}(t),\phi(\cdot,t)) (h_{\e}(\cdot,V_{\e}(t))) \: dt \: .
\end{equation}
We eventually conclude the proof of Step 4, hence of the Proposition~\ref{prop:timeepsilonbrakke} thanks to \eqref{eq:BrakkeEpsStep4:0}, \eqref{eq:BrakkeEpsStep4:1}, \eqref{eq:BrakkeEpsStep4:2} and \eqref{eq:BrakkeEpsStep4:3}.
\end{proof}
We conclude this section by noting a straightforward though important consequence of the Brakke-type equality we established in Proposition~\ref{prop:timeepsilonbrakke}: the mass $t \mapsto \| V_\e(t) \| (\R^n)$ is non-increasing along the flow.

\begin{remk}[Mass decay]
\label{remk:massdecrease4}
Let $\e \in(0,1)$ and let $V_0 \in V_d(\R^n)$ be a compactly supported varifold satisfying $\|V_0\|(\R^n)\leq M$. Let $(V_{\e}(t))_{t \in [0,1]}$ be the approximate mean curvature flow starting from $V_0$. Then, for all $0 \leq a < b \leq 1$,
\begin{equation*}
    \|V_{\e}(b)\|(\R^n) - \|V_{\e}(a)\|(\R^n) = - \int_a^b \int_{\R^n} \frac{|(\Phi_{\e}\ast\delta V_\e(t))(y)|^2}{(\Phi_{\e}\ast \| V_\e(t) \| )(y)+\e} \: dy \: dt \leq 0 \: .
\end{equation*}

\noindent In particular, 
\begin{equation}\label{remk:massdecrease4_eq2}
 \|V_{\e}(t)\|(\R^n)  \leq \|V_0\|(\R^n), \,\, \forall t \in [0,1],
\end{equation}
and 
\begin{equation*}
 \int_0^1 \int_{\R^n} \frac{|(\Phi_{\e}\ast\delta V_\e(t))(y)|^2}{(\Phi_{\e}\ast \| V_\e(t) \| )(y)+\e} \: dy \: dt \leq \|V_0\|(\R^n).
\end{equation*}
\end{remk}
Indeed, let $0 \leq a < b \leq 1$, applying Proposition~\ref{prop:timeepsilonbrakke} with the constant test function $\phi :(t,y) \mapsto 1 \in \xC^1([0,1] \times \R^n, \R_+)$ we have $\partial_t \phi = 0$, $\delta (V_\e(t) , \phi( \cdot, t)) = \delta V_\e(t)$ (see \eqref{wfirstvar}) and then using \eqref{massdecrease1} in Proposition~\ref{mcfmotion1}, we infer
\begin{equation*}
 \|V_{\e}(b)\|(\R^n) - \|V_{\e}(a)\|(\R^n) = \int_{a}^{b} \delta V_{\e}(t) \left(h_{\e}(\cdot,V_{\e}(t))\right) \: dt = - \int_a^b \int_{\R^n} \frac{|(\Phi_{\e}\ast\delta V_\e(t))(y)|^2}{(\Phi_{\e}\ast\|V_\e(t)\|)(y)+\e} \: dy \: dt.
\end{equation*}

%
\section{Convergence of approximate mean curvature flows to spacetime Brakke flows}\label{sec:spacetime-Brakke-flows}

Up to this point, given a compactly supported varifold $V_0\in V_d(\R^n)$, we constructed an approximate mean curvature flow $\left(V_{\e}(t)\right)_{t\in[0,1]}$ (Theorem \ref{damcfconvergence}) obtained as the limit of time-discrete approximate mean curvature flows (Definition \ref{damcf}) when the subdivision time step converges to $0$, the approximation scale $\e$ being fixed. In this section, we investigate the behavior of  $\left(V_{\e}(t)\right)_{t\in[0,1]}$ in the limit $\e \rightarrow 0$.

Following the works of Brakke and Kim \& Tonegawa \cite{brakke,kt}, one can prove the convergence of $\left(\| V_{\e}(t)\| \right)_{t\in[0,1]}$ to a limit measure $\mu(t)$ for all $ t \in [0,1]$, up to an extraction that is independent of $t$. In~\cite{brakke,kt}, the convergence is first established for a common sequence $( \e_j)_j$ for all dyadic numbers of $[0,1]$ thanks to the uniform boundedness of the mass and a diagonal extraction argument. Then, it is extended to all $t\in [0,1]$ using the continuity property of $t \mapsto \| V_{\e}(t)\| $. The previous procedure does not work for the measures $\left(V_{\e}(t)\right)_{t\in[0,1]}$ because of the lack of a continuity property. This issue is very common when studying the compactness of Brakke-type flows (cf. for instance \cite[Theorem 3.7]{ton}). To circumvent this issue, we consider the tensor product measure $\dt \otimes V_{\e}(t)$.

In section \ref{chap3_sec_st}, we introduce the notions of spacetime mean curvature and  spacetime Brakke flow, and we list some of their properties. In section \ref{chap3_sec_conv}, we prove that, up to an extraction, the measure  $\dt \otimes V_{\e}(t) $ converges, thanks to the uniform boundedness of the mass in $\e$, to a limit measure, denoted by $\lambda$. We prove that $\lambda$ has a bounded spacetime mean curvature in $L^2$, and that it is a spacetime Brakke flow provided its $\R^n \times \G$-component is rectifiable.

We introduce after~\cite{kt} classes of test functions and vector fields that are suitable for studying the behaviour of the approximate mean curvature flows. For $ j \in \N$ we define 
\begin{equation*}
 \mathcal{A}_j := \lbrace \phi \in C^2(\R^n;\R_+) : \phi(x) \leq  1 , |\nabla\phi(x)|\leq j \phi(x) , \|\nabla^2\phi(x)\|\leq j \phi(x) \,  \text{for all} \,  x \in \R^n \rbrace,
\end{equation*}
\begin{equation*}
  \mathcal{B}_j := \lbrace g \in C^2(\R^n;\R^n) : |g(x)|\leq j , \|\nabla g(x)\|\leq j, \|\nabla^2 g(x)\|\leq j \,  \mbox{for all} \,  x \in\R^n \, \mbox{and} \, \|g\|_{L^2(\R^n)}\leq j \rbrace.
\end{equation*}

\subsection{Spacetime mean curvature and spacetime Brakke flows}\label{chap3_sec_st}

As we explained, our purpose is to investigate the limit of $\lambda_\e = \dt \otimes V_\e(t)$ when $\e$ tends to $0$. Let us start by recalling the framework that allows to consider such a generalized tensor product of measures. 
Following \cite[Section 2.5]{afp}, we say that {\it a family of $d$--varifolds $(W(t))_{t \in [0,1]}$ is $\dt$--measurable} if for any Borel set $B \subset \R^n \times \G$, the map $t \in [0,1] \mapsto W(t)(B) \in \R_+$ is $\dt$--measurable. In this case, the measure $\dt \otimes W(t)$ is well-defined as follows:
\begin{equation*}
    \forall B \subset [0,1] \times \R^n \times \G \text{  Borel set,} \quad \left[ \dt \otimes W(t) \right] (B) = \int_{[0,1]} \left( \int_{\R^n \times \G} \one_B \: d W(t) \right) \dt
\end{equation*}
and satisfies that for any Borel function $f : [0,1] \times \R^n \times \G \rightarrow \R_+$,
\begin{equation*}
    \int f \:  d \left[ \dt \otimes W(t) \right] = \int_{t \in [0,1]} \left( \int_{(y,S) \in \R^n \times \G} f(t,y,S) \: d W(t) \right) \dt \: .
\end{equation*}
Note that the above assertion then holds for any $f \in L^1(\dt \otimes W(t))$. 
We point out that in such a case the measure $\dt \otimes W(t)$ is not a Radon measure in general. Indeed, given a compact set $K \subset \R^n$, $K^\prime = [0,1] \times K \times \G$ is compact and
\[
\left[\dt \otimes W(t) \right](K^\prime) = \int_0^1 \| W(t) \|(K) \: dt \: .
\]
As any compact set of $[0,1] \times \R^n \times \G$ can be included in a compact set of the form $K^\prime$ we conclude that $\dt \otimes W(t)$ is a Radon measure if and only if for all compact sets $K \subset \R^n$,
\[
\int_0^1 \| W(t) \|(K) \: dt < \infty \: .
\]
We now introduce the \emph{spacetime mass} and \emph{first variation} of a Radon measure $\beta$ defined in $[0,1] \times \R^n \times \G$. Such definitions naturally arise when requiring that they are consistent with the usual mass and first variation when $\beta = \dt \otimes W(t)$ (see Remark~\ref{rmkSpaceTimeMC}).

\begin{dfn}[Spacetime mean curvature]\label{spacetimeMC}
 Let $\beta$ be a Radon measure on $[0,1] \times \R^n \times \G$. 
 \begin{itemize}
 \item We introduce the \emph{spacetime mass} $\|\beta\|$ of $\beta$: $\|\beta\| := \Pi_\# \beta$ is a Radon measure in $[0,1] \times \R^n$, where $\Pi : [0,1] \times \R^n \times \G \rightarrow [0,1] \times \R^n$ is the canonical projection.
 \item We define the \emph{first variation of $\beta$}: for any vector field $X\in \xC_c^1([0,1] \times \R^n,\R^n)$,
 \begin{equation*}
  \delta \beta(X) := \int_0^1\int_{\R^n\times \G} \mdiv_S X(t,y) d\beta(t,y,S).
 \end{equation*}
 \item If in addition, $\delta \beta$ is bounded in the $\xC_c([0,1] \times \R^n,\R^n)$--topology, we say that $\beta$ has locally bounded first variation, and then, by the Riesz representation theorem, $\delta \beta$ identifies with a vector-valued Radon measure (also denoted by $\delta \beta$)
\[
 \delta \beta (X) = \int_{[0,1] \times \R^n} X \cdot \: d \delta \beta \quad \forall X \in \xC_c ([0,1] \times \R^n, \R^n) \: .
\]
It follows from the Radon--Nikodym decomposition theorem that there exist a Radon measure $(\delta\beta)_s$ singular with respect to $\|\beta\|$, and $h=h(\cdot,\cdot,\beta) \in L^1([0,1] \times \R^n, \R^n, \|\beta\|)$ that we call \emph{the spacetime mean curvature}, satisfying $\delta \beta = - h \| \beta \| + (\delta\beta)_s$, i.e.:
\begin{equation}\label{eq:spacetimeMC}
 \delta \beta(X) = -\int_0^1 \int_{\R^n} h(t,y,\beta) \cdot X(t,y) d\|\beta\|(t,y) + (\delta\beta)_s(X),
\end{equation}
for any $X\in \xC_c ([0,1] \times \R^n,\R^n)$.
\end{itemize}
\end{dfn}
\begin{remk} \label{rmkSpaceTimeMC}
Let $(W(t))_{t \in [0,1]}$ be a family of $d$--varifolds that is $\dt$--measurable
and such that $\beta = \dt \otimes W(t)$ is a Radon measure. Then, $\dt \otimes \| W(t) \|$ is a well-defined Radon measure and $\| \beta \| = \dt \otimes \| W(t) \|$ ; and for any $X \in \xC_c^1([0,1] \times \R^n, \R^n)$,
\begin{equation*}
    \delta \beta (X) = \int_0^1 \delta W(t) (X) \: dt \: .
\end{equation*}
Assume moreover that for a.e $t \in [0,1]$, $W(t)$ has locally bounded first variation and $\delta W(t) \ll \| W(t) \|$:
\[
 \delta W(t) = - h( \cdot , W(t)) \| W(t) \|  \text{ with } h(\cdot , W(t)) \in L_{loc}^1(\R^n, \|W(t)\|) \: . 
\]
Then, $\beta$ has locally bounded first variation in the sense of Definition~\ref{spacetimeMC} if and only if
$$
\int_{[0,1] \times \R^n} | h(y , W(t)) | \: d W(t) \: dt < \infty
$$
and in this case,
\begin{multline} \label{eqConsistencySpaceTimeMC}
    \delta \beta = - \dt \otimes h( \cdot , W(t)) \: d W(t) \quad \text{that is}\\
    \delta \beta \ll \|\beta\| \quad \text{and} \quad h(t,y,\beta) = h(y,W(t))  \quad \text{for $\|\beta\|$--a.e} \, (t,y)\in [0,1] \times \R^n.
\end{multline}
Indeed, for $X \in \xC_c^1([0,1] \times \R^n, \R^n)$
 \begin{equation*}\begin{split}
  \delta\beta(X) & = \int_{t=0}^1 \left( \int_{\R^n\times\G} \mdiv_S X(t,y) \, d W(t)(y,S) \right) dt   \\
  &  = -\int_0^1 \int_{\R^n} h( \cdot , W(t)) \cdot X(t, \cdot) \,d \|W(t)\| \: dt \\
& = -\int_{[0,1] \times \R^n} h(y,W(t)) \cdot X(t,y) d \|\beta\| (t,y) \:,
\end{split}\end{equation*}
and we infer \eqref{eqConsistencySpaceTimeMC}.
\end{remk}

We now give the definition of spacetime Brakke flows.

\begin{dfn}(Spacetime Brakke flows)\label{def:spacetimebf}
Let $\lambda$ be a Radon measure on $[0,1] \times \R^n \times \G$.
The measure $\lambda$ is called a spacetime Brakke flow if: 
\begin{itemize}
\item[$(i)$] there exists a family of $d$--varifolds $(V(t))_{t \in [0,1]}$ such that $\lambda = \dt \otimes V(t)$,
 \item[$(ii)$] $\lambda$ has locally bounded first variation, $\delta \lambda \ll \| \lambda \|$, and there exists $h = h(\cdot, \cdot, \lambda) \in L^2_{loc}([0,1] \times \R^n , \R^n, \|\lambda\|)$ such that
 \begin{equation*}
     \delta \lambda = - h \: \| \lambda \| \: .
 \end{equation*}
 \item[$(iii)$] (Integral Brakke inequality). For any $\phi \in \xC_c^1([0,1] \times \R^n, \R_+)$ and for all $0 \leq t_1 \leq t_2 \leq 1$ we have
 \begin{equation}\label{intbrakkeineq}
  \begin{split}
 \|V(t_2)\|(\phi(t_2, \cdot))& -  \|V(t_1)\|(\phi(t_1, \cdot))  \leq -\int_{t_1}^{t_2}\int_{\R^n} \phi(t, \cdot )|h(t, \cdot ,\lambda)|^2 \, d \|V(t)\| \: dt \\
+ &\int_{t_1}^{t_2}\int_{\R^n \times \G } S^{\perp}(\nabla_y \phi(t,y)\cdot h(t,y,\lambda) \, d\lambda(t,y,S) + \int_{t_1}^{t_2}\int_{\R^n} \partial_t \phi(t, \cdot)\, d \|V(t)\| \: dt
\end{split}
 \end{equation}
\end{itemize}
We furthermore say that $\lambda$ starts from $V(0)$.
\end{dfn}
\begin{remk}
We give two important consequences of the definition of the spacetime Brakke flow in the case where $h \in L^1 (\| \lambda \|)$:
\begin{itemize}
 \item[$(i)$] {\it Mass decay.} For all $ 0 \leq t_1 \leq t_2 \leq 1$,
 \begin{equation*}
  \|V(t_2)\|(\R^n) \leq \|V(t_1)\|(\R^n) \leq \|V(0)\|(\R^n) \: .
\end{equation*}
 \item[$(ii)$]{\it $L^2$ bound on the mean curvature.} If $V(0)$ has finite mass $\|V(0)\| < \infty$ then $h \in L^2(\|\lambda\|)$ and more precisely,
 \begin{equation*}
\| h \|_{L^2(\|\lambda\|)} = \int_{0}^1 \int_{\R^n} |h(t,y,\lambda)|^2 \, d\|V(t)\|(y) \: dt \leq \|V(0)\|(\R^n) \: .
\end{equation*}
\end{itemize}
\end{remk}
\begin{proof}[Proof of $(i)$ and $(ii)$]
Given a cutoff function $\chi \in \xC_c^1(\R^n,\R_+)$ satisfying $0 \leq \xi \leq 1$, $\chi = 1$ in $B_1(0)$ and $\supp \chi \subset B_3(0)$,  
we introduce for $r > 0$, $\phi_r \in \xC_c^1([0,1] \times \R^n, \R_+)$ as: 
\[
\forall (t,x) \in [0,1] \times \R^n, \,\, \phi_r(t,x) = \chi \left( \frac{x}{r} \right) \: .
\]
Note that $0 \leq \phi_r \leq 1$ satisfies $\supp \phi_r \subset [0,1] \times B_{3/r}(0)$, $\partial_t \phi = 0$ and $\| \nabla_y \phi_r \|_{\infty} \leq r^{-1}$. Plugging $\phi_r$ in \eqref{intbrakkeineq} we obtain (denoting $h := h(\cdot, \cdot,\lambda)$ for simplicity)  
\begin{align*}
 \|V(t_2)\|(\phi_r) & \leq \|V(t_1)\|(\phi_r) -\int_{t_1}^{t_2}\int_{\R^n} \phi_r|h|^2 \, d\mu(t)dt
 +\int_{t_1}^{t_2} \int_{\R^n \times \G } S^{\perp}(\nabla\phi_r)\cdot h \, d\lambda \\
 & \leq \|V(t_1)\|(\R^n) -\int_{t_1}^{t_2} \int_{B_{1/r}(0)} |h|^2 \: d\|V\|(t) \: dt  + \frac{1}{r} \| h \|_{L^1 (\| \lambda\|)} .
 \end{align*}
Since $\|V(t_2)\|(\phi_r) \geq \|V(t_2)\|(B_\frac{1}{r}(0)) \xrightarrow[r \to 0]{} \| V(t_2)\|(\R^n)$ and $h \in L^1(\|\lambda\|)$ we let $r\rightarrow +\infty$ and obtain
\begin{equation*}
 \|V(t_2)\|(\R^n) + \int_{t_1}^{t_2} \int_{\R^n} |h|^2 \: d\|V\|(t) \: dt  \leq \| V(t_1) \| (\R^n)
\end{equation*}
thus proving $(i)$.
Then specifying $t_1 = 0$ and $t_2 = 1$ allows to conclude the proof of $(ii)$.
\end{proof}
\begin{remk}[Brakke flows and spacetime Brakke flows]
Let $(V(t))_{t\in[0,1]}$ be a Brakke flow (see Definition~\ref{intro_def_brakke-V2}), then $\lambda= \dt \otimes V(t)$ is a spacetime Brakke flow in the sense of Definition~\ref{def:spacetimebf}.
\end{remk}

\subsection{Convergence of approximate mean curvature flows, properties of the limit} \label{chap3_sec_conv}
We now study the convergence of the approximate mean curvature flow  $\left( V_{\e}(t) \right)_{t \in [0,1]}$ introduced in Section~\ref{def_time_step_to_0}. As previously mentioned, it is possible to find an extracted sequence $(\e_j)_j$ that does not depend on $t \in [0,1]$ and a family of Radon measures $(\mu(t))_{t \in [0,1]}$ such that for all $t \in [0,1]$, $\| V_{\e_j} \|$ converges to $\mu(t)$, as shown in Proposition~\ref{limitmass}. However, we are not able to obtain a similar result concerning the whole varifold structure: it is possible to have an extracted sequence $(\e_j^t)_j$ such that for each time $t \in [0,1]$, $V_{\e_j^t}(t)$ converges to some limit varifold, but we are not able to remove the time dependency of the extracted sequence $(\e_j^t)_j$. Therefore, in Theorem~\ref{chap3_theo:cv}, we investigate the limit of the measures $\dt \otimes V_\e (t)$ when $\e$ tends to $0$ instead.

\begin{theo}(Convergence)\label{chap3_theo:cv}
Let $V_0 \in V_d(\R^n)$ of compact support and finite mass. Given $\e \in (0,1)$, $(V_{\e}(t))_{t\in[0,1]}$ denotes the approximate mean curvature flow starting from $V_0$. We have:
\begin{enumerate}[label=(\roman*)]
 \item For $\e \in (0,1)$, the measure $\dt \otimes V_e(t)$ is well-defined and is a Radon measure in $[0,1] \times \R^n \times \G$.
 \item There exists a sequence $(\e_j)_j \xrightarrow[j\rightarrow \infty]{} 0$ such that 
 \begin{equation*}
  \dt \otimes V_{\e_j}(t)  \xrightarrow[j\rightarrow \infty]{}  \lambda = \dt \otimes V(t)  \quad \text{and} \quad \forall t \in [0,1], \,\, \| V_{\e_j} (t) \| \xrightarrow[j\rightarrow \infty]{}  \|V(t)\|,
 \end{equation*}
where $(V(t))_{t \in [0,1]}$ is a family of $d$--varifolds.

\item $\delta \lambda$ is locally bounded, $\delta \lambda = - h (\cdot, \cdot, \lambda) \, \| \lambda \|$ and $\|h(\cdot,\cdot,\lambda)\|_{L^2(\|\lambda)\|}^2 \leq \|V_0\|(\R^n)$.

\item If we assume that $V(t)$ is a rectifiable varifold for a.e $ t\in [0,1]$ then $\lambda$ is a spacetime Brakke flow.
\end{enumerate}
\end{theo}

We emphasize that, in this theorem, we have the convergence of the mass $\| V_{\e_j} (t) \| \xrightarrow[j\rightarrow \infty]{}  \|V(t)\|$ for all $t \in [0,1]$ for a sequence $\e_j \to 0$ that does not depend on $t$, however, we do not know whether $V_{\e_j}(t)$ converges to $V(t)$.

\begin{remk} The proof of Theorem~\ref{chap3_theo:cv} is crucially based on Propositions 5.4 and 5.5 in~\cite{kt}. Note however that they are stated in~\cite{kt} for codimension $1$ varifolds, while we work with varifolds of arbitrary codimension in $\R^n$. Actually, the proofs of Propositions 5.4 and 5.5 in~\cite{kt} do not require that the varifolds have codimension $1$ because they are based on results due to Brakke, which are valid for general varifolds.
\end{remk}

For simplicity, we split the proof of Theorem~\ref{chap3_theo:cv} into several steps: it will follow from Propositions~\ref{prop:kt5455} to \ref{spacetimebrakkeinequality}.
The next proposition combines the generalized estimates (generalized in terms of codimension) of \cite[Proposition 5.4]{kt} and \cite[Proposition 5.5]{kt}.

\begin{prop} \label{prop:kt5455}
There exists some $\epsilon_\ast = \epsilon_\ast(n,M) > 0$ (depending only on $n$ and $M$) such that for any varifold $W \in V_d(\R^n)$ such that $\| W \|(\R^n) \leq M$; if $0 < \e \leq \e_\ast$ and $\e^{-\frac{1}{6}} \geq 2m$, we have 
\begin{itemize}
\item for any $\psi \in \mathcal{A}_m$
\begin{equation}\label{eqKTprop54}
\left| \delta W \left(\psi h_{\e}(\cdot, W) \right) + \int_{\R^n} \frac{\psi |\Phi_{\e}\ast\delta W|^2 }{\Phi_{\e}\ast \|W\|+\e} \: dx \right| \leq \e^\frac{1}{4} \left( 1 + \int_{\R^n} \frac{\psi |\Phi_{\e}\ast\delta W|^2 }{\Phi_{\e}\ast \|W\|+\e} \: dx \right)
\end{equation}
and
\begin{equation}\label{eq2KTprop54}
 \int_{\R^n} \psi|h_{\e}(\cdot,W)|^2 \, d\W \leq \int_{\R^n} \frac{\psi |\Phi_{\e} \ast \delta W|^2}{\Phi_{\e}\ast \W + \e}(1+\e^{\frac14}) \,dx + \e^{\frac14},
\end{equation}
 \item for any $X\in  \mathcal{B}_m$
\begin{equation}\label{eqKTprop55}
 \Big| \int_{\R^n}h_{\e} \cdot X \, d\|W\| + \delta W(X) \Big| \leq \e^{\frac14} + \e^{\frac14} \left( \int_{\R^n} \frac{|\Phi_{\e} \ast \delta W|^2}{\Phi_{\e} \ast \|W\| + \e} \, dx\right)^{\frac12}. 
\end{equation}
\end{itemize}
\end{prop}
%


%
\begin{prop}[Existence of a limit for the mass measure]\label{limitmass}
Let $V_0 \in V_d(\R^n)$ of bounded support, for any $\e\in(0,1)$, let $V_{\e}(t)$ be the approximate mean curvature flow starting from $V_0$.\\
There exists a sequence $(\e_{j})_{j=1}^{\infty}$ (not depending on $t$) converging to $0$ as $j \rightarrow \infty$, and a family of Radon measures $(\mu(t))_{t\in[0,1]}$ on $\R^n$ such that: 
\begin{equation}\label{meconeq}
\|V_{\e_{j}}(t)\| \xrightharpoonup[j\rightarrow \infty]{\ast} \mu(t)
\end{equation} 
for \underline{all} $t\in [0,1]$.
\end{prop}
\begin{proof}
 Proposition \ref{limitmass} states that, up to an extraction independent of $t$, the mass measure $\| V_{\e}(t) \|$ converges as $\e$ goes to $0$ to a limit time-dependent measure $\mu(t)$. The proof is a direct adaptation of the proof of \cite[Proposition 6.4 $(1)$]{kt}, which is itself based on the results of Section 5 in~\cite{kt} and two estimates \cite[$(6.3)$]{kt}, \cite[$(6.5)$]{kt}. In short, the proof of~\cite{kt} is based on the following arguments: a limit measure $\mu(t)$ is defined for a countable, dense collection $D$ of times using an extraction argument, the extension of $\mu(t)$ to almost all times, more precisely the complement of a countable set of  times, follows from a continuity property of $\mu(t)$ on $D$, the convergence of the subsequence $\| V_{\e}(t) \|$ to $\mu(t)$ for all these times $t$ follows from an approximate continuity property satisfied by $\| V_{\e}(t) \|$, and a last extraction is used to recover $\mu(t)$ for all $t$.
 The adaptation of the proof to our framework uses the following arguments:
\begin{itemize}
 \item $V_0 \in V_{d}(\R^n)$ is not necessarily an open partition in our case, but the proof in~\cite{kt} remains valid in this case;
 \item the proof of~\cite{kt} is based on the results of \cite[Section 5]{kt} which are valid for varifolds of any codimension (and not only for codimension $1$ varifolds which are the subject of \cite[Proposition 6.4$(1)$]{kt});
 \item the weight function $\Omega$ is identically equal to $1$ in our case because we work with varifolds of finite mass;
 \item estimate \cite[Inequality $(6.3)$]{kt} is replaced by the decay property of the mass stated in Remark \ref{remk:massdecrease4};
 \item estimate \cite[Inequality $(6.5)$]{kt} is replaced by \eqref{timeepsilonbrakke} (with test functions depending only on the space variable). 
 \end{itemize}
\noindent {\bf Step 1:} Let $\e \in (0,1), \, V \in V_d(\R^n)$. We prove the following inequality  
\begin{equation*}
 \delta(V,\phi)(h_{\e}(\cdot,V)) \leq 2\e^{\frac14} + \frac12 \int_{\R^n} \frac{|\nabla \phi|^2}{\phi} \,d\|V\|,
\end{equation*}
for any $\phi \in \cA_j$ and $2j\leq \e^{-\frac16}.$ 

\noindent Indeed, let $\phi \in \cA_j$ with $2j\leq \e^{-\frac16}$ and set for simplicity $\displaystyle b :=\int_{\R^n} \frac{\phi|\Phi_{\e}\ast\delta V|^2}{\Phi_{\e}\ast\|V\|+\e} dx $, we have 
\begin{align*}
\begin{split}
    &\delta\left( V,\phi \right) \left( h_{\e}(\cdot,V) \right) = \delta V(\phi h_{\e}) + \int_{\R^n} S^{\perp}(\nabla\phi)\cdot h_{\e}(\cdot,V) \, d V(x,S)\\&
    \leq -b+\e^{\frac{1}{4}}b +\e^{\frac{1}{4}}+ \frac{1}{2}\int_{\R^n} \phi |h_{\e}(\cdot,V)|^2 \, d\|V\|+ \frac{1}{2}\int_{\R^n} \frac{|\nabla\phi|^2}{\phi} \, d\|V\| \,\, \text{by \eqref{eqKTprop54}}
    \\& \leq -b+\e^{\frac{1}{4}}b +\e^{\frac{1}{4}}+ 
     \frac{1}{2}b (1+\e^{\frac{1}{4}}) +\frac{1}{2}\e^{\frac{1}{4}} + \frac{1}{2} \int_{\R^n} \frac{|\nabla\phi|^2}{\phi} \, d\|V\| \,\, \text{by \eqref{eq2KTprop54}} 
     \\& \leq \frac{1}{2}(-1+3\e^{\frac{1}{4}})b + 2\e^{\frac{1}{4}} + \frac12\int_{\R^n} \frac{|\nabla\phi|^2}{\phi} \,d\|V\|
     \\& \leq 2\e^{\frac{1}{4}} + \frac12\int_{\R^n} \frac{|\nabla\phi|^2}{\phi} \, d\|V\|.
\end{split}
\end{align*}
\noindent {\bf Step 2:} We define the limit measure $\mu(t)$ for a.e $t\in [0,1]$.

\noindent Let $D\cap[0,1]$ be the set of dyadic numbers in $[0,1]$. Let $(\e_j)_{j\in \N}$ be a sequence converging to $0$ such that, $\|V_{\e_j}(t)\|$ converges for any $t\in D$, denote the limit $\mu(t)$. The previous claim stems from the Banach-Alaoglu theorem and the uniform boundedness of the mass \eqref{remk:massdecrease4_eq2}.\\
Let $Z:=\left( \phi_q \right)_{q\in\N}$ be a countable subset of $C_{c}^2(\R^n,\R_+)$ which is dense in $C_{c}(\R^n,\R_+)$. Take $\phi_{q} \in Z$ and assume without loss of generality that $\phi_q \sleq 1$. Then, for any $i \in \N$ large enough, we have $\phi_q+i^{-1} \in \mathcal{A}_m$ for any $m \geq m_0$, where $m_0$ depends on $i$ and $\phi_q$. We apply \eqref{timeepsilonbrakke} with $\phi(\cdot,t)=\phi_q$, together with Step~1 to obtain
\begin{equation*}
 \|V_{\e_j}(b)\|(\phi_q+i^{-1}) - \|V_{\e_j}(a)\|(\phi_q+i^{-1}) \leq (b-a)2\e_j^{\frac14} + \frac12\int_{a}^{b}\int_{\R^n} \frac{|\nabla(\phi_q+i^{-1})|^2}{\phi_q+i^{-1}} \, d\|V_{\e_j}(t)\|dt,
\end{equation*}
for any $a,b \in [0,1], a \leq b$ and $2m_0\leq \e_j^{-\frac{1}{6}}$. We obtain from \cite[Lemma 3.1]{ton}
\begin{equation*}
 \frac{|\nabla(\phi_q+i^{-1})|^2}{\phi_q+i^{-1}} \leq  \frac{|\nabla \phi_q |^2}{\phi_q} \leq 2\|\nabla^2 \phi_q \|_{\infty}.
\end{equation*}
Therefore, for any $a,b \in [0,1], a \leq b$
\begin{equation}\label{eq:monotone_prf1}
 \|V_{\e_j}(b)\|(\phi_q+i^{-1}) - \|V_{\e_j}(a)\|(\phi_q+i^{-1}) \leq (b-a)2\e_j^{\frac14} + (b-a)\|\nabla^2 \phi_q \|_{\infty}\|V_0\|(\R^n).
\end{equation}
We let $j \rightarrow \infty$, we deduce for $a,b \in D, a \leq b$ that 
\begin{equation*}
\mu(b)(\phi_q+i^{-1}) - \mu(a)(\phi_q+i^{-1}) \leq (b-a)\|\nabla^2 \phi_q \|_{\infty}\|V_0\|(\R^n).
\end{equation*}
We let $i \rightarrow \infty$, using the uniform boundedness of the mass \eqref{remk:massdecrease4_eq2}, we deduce for $a,b \in D, a \leq b$ that 
\begin{equation*}
\mu(b)(\phi_q) - \mu(a)(\phi_q) \leq (b-a)\|\nabla^2 \phi_q \|_{\infty}\|V_0\|(\R^n).
\end{equation*}
The previous inequality tells us that the map $ g_q : t \mapsto \mu(t)(\phi_q)-t(b-a)\|\nabla^2\phi_q\|_{\infty}\|V_0\|(\R^n)$ is nonincreasing for $t\in D$. Define
\begin{equation*}
 \cC := \lbrace t \in [0,1], \text{for some $q\in\N$} \, \lim_{s \rightarrow t^-} g_q(s) > \lim_{s \rightarrow t^+} g_q(s) \rbrace.
\end{equation*}
By the monotonicity property of $g_q$, $\cC$ is a countable set in $[0,1]$, and $\mu(t)(\phi_q)$ may be defined continuously on the complement of $\cC$ uniquely from the values on $D$; then, one can  define the measure $\mu(t)$ for a.e $t\in[0,1]$ by density of $Z$ in $C_c(\R^n,\R_+)$.

\smallskip
\noindent {\bf Step 3:} We prove that for any $t\in [0,1] \setminus \cC$ 
\begin{equation*}
\|V_{\e_j}(t)\| \xrightharpoonup[]{*} \mu(t). 
\end{equation*}
\noindent Let $t\in [0,1] \setminus \cC$ and $s\in D,  t<s$. From \eqref{eq:monotone_prf1} we have 
\begin{equation*}
\|V_{\e_j}(s)\|(\phi_q+i^{-1}) \leq \|V_{\e_j}(t)\|(\phi_q+i^{-1}) + O(s-t).
\end{equation*}
We let $j\rightarrow $ so that
\begin{equation*}
 \mu(s)(\phi_q+i^{-1})  \leq \liminf_j  \|V_{\e_j}(t)\|(\phi_q+i^{-1}) + O(s-t).
\end{equation*}
Then we take the limit in $i$ to obtain
\begin{equation*}
 \mu(s)(\phi_q)  \leq \liminf_j  \|V_{\e_j}(t)\|(\phi_q) + O(s-t).
\end{equation*}
We now let $s \rightarrow t^-$ and use the continuity of $g_q$ at $t$ so that 
\begin{equation*}
 \mu(t)(\phi_q)  \leq \liminf_j  \|V_{\e_j}(t)\|(\phi_q).
\end{equation*}
The same reasoning for $s<t$ gives $\mu(t)(\phi_q)  \geq  \limsup_j  \|V_{\e_j}(t)\|(\phi_q)$; hence
\begin{equation*}
 \lim_j \|V_{\e_j}(t)\|(\phi_q) =  \mu(t)(\phi_q), \quad \forall \phi_q \in Z \,\, \text{and} \, \, \forall t \in [0,1] \setminus \cC;
\end{equation*}
 we conclude the proof of step 3 by density of $Z$ in $C_c(\R^n,\R_+)$.

\vspace{0.2cm}
The set $\cC$ is countable, hence, by further extraction of the sequence $(\e_j)_j$, we can define $\mu(t)$ on $[0,1]$ entirely and ensure the convergence for \underline{all} $ t\in[0,1]$.
\end{proof}

\begin{prop}[Limiting flow]\label{limitflow}
Let $V_0 \in V_d(\R^n)$ be of compact support. 
\begin{enumerate}[label=(\roman*)]
    \item For any $\e\in(0,1)$,
let $V_{\e}(t)$ be the approximate mean curvature flow starting from $V_0$. Then $(V_\e(t))_{t\in [0,1]}$ is $\dt$--measurable and $\lambda_\e = \dt \otimes V_{\e}(t)$ is a Radon measure, satisfying for all $\phi \in \xC_c([0,1] \times \R^n \times \G, \R)$,
 \[
  \int \phi \: d \lambda_{\e} = \int_{0}^1 \left( \int_{(x,S) \in \R^n \times \G} \phi (t,x,S) \: d V_\e(t) \right) \: dt. 
 \]
    \item There exists a sequence $( \e_j )_j \rightarrow 0$ and a Radon measure $\lambda$ in $[0,1] \times \R^n \times \G$ such that
\begin{equation} \label{eqCVlambdaj}
 \lambda_{\e_j} \xrightharpoonup[j \to \infty]{\ast} \lambda
\end{equation}
and furthermore, $\lambda = \dt \otimes V(t)$ where $(V(t))_{t\in [0,1]}$ is a family of $d$--varifolds.
    \item The mass of the varifolds converges: for all $t \in [0,1]$, $\|V_{\e_j}(t)\|  \xrightharpoonup[j \to \infty]{\ast} \|V(t)\|$.
\end{enumerate}
\end{prop}
\begin{proof}
Let $B \subset \R^n \times \G$ be a Borel set and $f : [0,1] \rightarrow \R_+$, $t \mapsto V_\e(t) (B)$. We check that $f$ is $\dt$--measurable. According to~\cite[Prop 2.26]{afp}, it is sufficient to check that $f$ is measurable in the case where $B$ is an open set, if so,
\begin{align*}
 f(t) = V_\e(t) (B) & = \sup \left\lbrace \int \phi \: d V_\e(t) \: : \: \phi \in \xC_c (B, \R_+), \, 0 \leq \phi \leq 1 \right\rbrace \\
 & = \sup \left\lbrace \int \phi \: d V_\e(t) \: : \: \phi \in \mathcal{D}, \, 0 \leq \phi \leq 1 \right\rbrace 
\end{align*}
where $\mathcal{D}$ is countable and dense in $\xC_c (B, \R_+)$ and the measurability of $f$ reduces to the measurability of $g : t \mapsto \int \phi \: d V_\e(t)$ given any $\phi \in \mathcal{D} \subset \xC_c (B, \R_+)$.
Recall that for any $t \in [0,1]$, $V_\e(t)$ is the weak star limit of $V_{\e, \cT_j^D}^{pc}(t)$ (where $\cT_j^D$ refers to dyadic subdivisions, see Theorem~\ref{damcfconvergence} and \eqref{eqPWcv}) so that 
\begin{equation*}
 g(t) = \int \phi \: d V_\e(t) = \lim_{j \to \infty} g_j(t) \quad \text{with} \quad g_j(t) = \int \phi \: d V_{\e, \cT_j^D}^{pc}(t) \: .
\end{equation*}
The functions $g_j$ are measurable as they are piecewise constant since $t \mapsto V_{\e, \cT_j^D}^{pc}(t)$ is piecewise constant by definition (see Remark~\ref{remk:PWflow}), and therefore $f$ is measurable and the family $(V_\e(t))_{t \in [0,1]}$ is $\dt$--measurable.
Moreover, thanks to \eqref{remk:massdecrease4_eq2}
\begin{equation} \label{eq:BanachAlaogluBound}
 \lambda_\e ([0,1] \times \R^n \times \G) =  \int_0^1 \| V_\e(t) \| (\R^n) \: dt \leq  \| V_0 \| (\R^n) < \infty
\end{equation}
and $\lambda_\e$ is a finite Radon measure, which concludes the proof of $(i)$.

We carry on with the proof of $(ii)$ and $(iii)$. We know by Proposition~\ref{limitmass} that there exists a sequence $(\e_j )_j \rightarrow 0$ for which $\|V_{\e_j}(t)\|$ converges to a limit measure $\mu(t)$ for all $t\in[0,1]$. Using \eqref{eq:BanachAlaogluBound},
we can assert by Banach-Alaoglu's compactness theorem that,
up to a further extraction, $\lambda_{\e_j} \xrightharpoonup[]{\ast} \lambda$, where $\lambda$ is a finite Radon measure on $[0,1] \times \R^n \times \G$.
\\
Denoting by $\Pi$ the canonical projection $(t,x,S) \in [0,1] \times \R^n \times \G \mapsto (t,x) \in [0,1] \times \R^n$ we now show that $\Pi_\# \lambda = \dt \otimes \mu(t)$. Indeed, on the one hand, as a consequence of \eqref{eqCVlambdaj}, we have 
\begin{equation} \label{eqPiLambda1}
 \Pi_\# \lambda_{\e_j} \xrightharpoonup[j \to \infty]{\ast} \Pi_\# \lambda \: .
\end{equation}
On the other hand, by definition of push-forward measure, we have for $\varphi \in \xC_c([0,1] \times \R^n,\R)$,
\begin{align}
 \Pi_\# \lambda_{\e_j} (\varphi)  = & \int_{[0,1] \times \R^n \times \G} \varphi \circ \Pi \: d \lambda_{\e_j} = \int_{(t,x,S) \in [0,1] \times \R^n \times \G} \varphi ( \Pi(t,x,S) ) \: d V_{\e_j}(t) \: dt  \nonumber \\
 = & \int_{t=0}^1 \int_{x \in \R^n} \varphi (t,x) \: d \| V_{\e_j}(t) \| \: dt \nonumber \\
 \xrightarrow[j \to \infty]{} & \int_{t=0}^1 \int_{x \in \R^n} \varphi (t,x) \: d \mu(t) \: dt. \label{eqPiLambda2} 
\end{align}
where the convergence follows from dominated convergence and for all $t\in [0,1]$, $\|V_{\e_j} \| \xrightharpoonup[j \to \infty]{\ast} \mu(t)$.
From \eqref{eqPiLambda1} and \eqref{eqPiLambda2} we obtain
\begin{equation*}
  \Pi_{\#} \lambda = \dt \otimes \mu(t) \: .
\end{equation*}
It follows by Young's disintegration theorem~\cite[Theorem 2.28]{afp} that there exists a family of probability measures $\left( \nu_{(t,x)} \right)_{(t,x)} $ on $\G$ (defined up to a $\dt \otimes \mu(t)$--null set), such that: 
\begin{equation*}
 \lambda = \Pi_\# \lambda \otimes \nu_{(t,x)} =   \dt \otimes \mu(t)\otimes \nu_{(t,x)},
\end{equation*}
which completes the proof (denoting $V(t) = \mu(t)\otimes \nu_{(t,x)}$).
\end{proof}
\begin{prop}\label{propSpaceTimeL2MC}
 Let $(\e_j)_j$ be the extracted sequence introduced in Proposition~\ref{limitflow}.
 Let $\lambda_{\e_j}  = \dt \otimes V_{\e_j} (t)$ and $\lambda$ be the limit measure defined in Proposition~\ref{limitflow}. We have: 
 \begin{equation} \label{eqSpaceTimeFVcv}
\forall X \in \xC_c^1 ([0,1] \times \R^n, \R^n), \quad \delta \lambda_{\e_j} (X) \xrightarrow[j \to \infty]{} \delta \lambda (X) \: .
 \end{equation}
Moreover, $\lambda$ has bounded first variation (i.e. $\delta \lambda$ is a finite Radon measure) and
\begin{equation} \label{eqSpaceTimeL2MC}
 \delta \lambda  = - h(\cdot ,\cdot,\lambda) \|\lambda\| = - h(t,\cdot,\lambda) \: \dt \otimes \| V(t)\| \quad \text{and} \quad \| h \|_{L^2(\|\lambda \|)}^2 \leq \|V_0\|(\R^n) < \infty \: ,
\end{equation}
Furthermore, for all bounded $\psi \in \xC ([0,1] \times \R^n, \R_+)$, for all $0 \leq t_1 < t_2 \leq 1$,
\begin{equation} \label{eqSpaceTimeL2LSC}
\int_{t_1}^{t_2} \int_{\R^n} \psi(t,y) |h(t,y,\lambda)|^2 \: d \|V(t)\| \: dt \leq \liminf_{j \to \infty}  \int_{t_1}^{t_2} \int_{\R^n} \psi(t,y) \frac{|\Phi_{\e_j}\ast\delta V_{\e_j}(t)|^2 }{\Phi_{\e_j}\ast \|V_{\e_j}(t)\|+\e_j} \: dy \: dt.
\end{equation}
\end{prop}
\begin{proof}
Let us first check that \eqref{eqSpaceTimeFVcv} is a consequence of the convergence $\lambda_{\e_j} \xrightharpoonup[j \to \infty]{\ast} \lambda$. Indeed, let $X \in \xC_c^1 ([0,1] \times \R^n, \R^n)$, then $g : (t,y,S) \mapsto \mdiv_S X(y) \in \xC_c ([0,1] \times \R^n \times \G)$ and thus
\[
 \delta \lambda_{\e_j} (X) = \int g \:  d\lambda_{\e_j} \xrightarrow[j \to \infty]{} \int g \: d \lambda = \delta \lambda (X) \: .
\]
Let us now consider the sequences $(\mu_j)_{j \in \N}$ of Radon measures in $[0,1] \times \R^n$ and $(f_j)_{j \in \N}$ of functions in $\xC^\infty([0,1] \times \R^n, \R^n)$ defined as
\[
 \mu_j = \dt \otimes (\Phi_{\e_j} \ast \|V_{\e_j}(t)\| + \e_j) {\rm d} y  \quad \text{and} \quad f_j (t, \cdot) = \frac{\Phi_{\e_j}\ast \delta V_{\e_j}(t)}{(\Phi_{\e_j} \ast \|V_{\e_j}(t)\| + \e_j)} \quad \text{for all } j \in \N.
\]
Let $\varphi \in \xC_c ([0,1] \times \R^n)$.
First note that by definition, 
\[
\int_{[0,1] \times \R^n} \varphi f_j \: d \mu_j = \int_0^1 \int_{\R^n} \varphi(t,y) \left( \Phi_{\e_j}\ast \delta V_{\e_j}(t) \right) dy \: dt \quad \Longrightarrow \quad f_j \mu_j = \dt \otimes (\Phi_{\e_j}\ast \delta V_{\e_j}(t)) {\rm d} y \: .
\]
We obtain by standard arguments (bicontinuity of distributional bracket to be more specific) that $\mu_j$ converge to $\| \lambda \| = \dt \otimes \| V(t)\|$ as Radon measures and $f_j \mu_j$ converges to $\delta \lambda$ as distributions of order $1$. Indeed, as $\Phi_\e$ is a mollifier, we recall that for all $t \in [0,1]$, $\| \varphi (t,\cdot) \ast \Phi_{\e_j} - \varphi(t, \cdot) \|_\infty  \xrightarrow[j \to \infty]{} 0$, and $\| \varphi (t,\cdot) \ast \Phi_{\e_j} - \varphi(t, \cdot) \|_{\xC^1} \xrightarrow[j \to \infty]{} 0$ if $\phi$ is additionally $\xC^1$.
Therefore, by dominated convergence and $\| V_{\e_j} \| (\R^n) \leq \| V_0 \|(\R^n)$,
\[
 \int_0^1 \left| \int_{\R^n} \varphi \: \left( d \Phi_{\e_j} \ast \| V_{\e_j} \|(t) - d \| V_{\e_j} \|(t) \right) \right| \: dt \leq \| V_0 \|(\R^n) \int_0^1 \| \varphi (t,\cdot) \ast \Phi_{\e_j} - \varphi(t, \cdot) \|_{\infty} \: dt \xrightarrow[j \to \infty]{} 0
\]
so that recalling that $\| \lambda_{\e_j} \| = \dt \otimes \| V_{\e_j}(t) \| $ converges to $ \| \lambda \|$,
\begin{align*}
\left| \int \varphi \: d \mu_j - \int \varphi \: d \| \lambda \| \right| & \leq \int_0^1 \left| \int_{\R^n} \varphi \: \left( d \Phi_{\e_j} \ast \| V_{\e_j} \|(t) - d \| V_{\e_j} \|(t) \right) \right| \: dt \\
& + \left| \int \varphi \: \left( d \| \lambda_{\e_j} \| - d \| \lambda \| \right) \right|
+ \e_j \left| \int \varphi \: dy \: dt \right| \xrightarrow[j \to \infty]{} 0 \: .
\end{align*}
We hence checked that $\mu_j \xrightharpoonup[j \to \infty]{\ast} \| \lambda \|$. In a very similar way, we can check that for all $X \in \xC_c^1 ([0,1] \times \R^n, \R^n)$, $\int X \cdot f_j \: d \mu_j \xrightarrow[j \to \infty]{} \delta \lambda (X)$ since
\begin{align*}
 \left|  \Phi_{\e_j}\ast \delta V_{\e_j}(t) (X) - \delta V_{\e_j}(t) (X)  \right| & = \left|   \delta V_{\e_j}(t) (\Phi_{\e_j}\ast X) - \delta V_{\e_j}(t) (X)  \right| \leq \| V_{\e_j}(t) \|(\R^n) \|\Phi_{\e_j}\ast X - X \|_{\xC^1} \\
 & \leq \| V_0 \| (\R^n) \|\Phi_{\e_j}\ast X - X \|_{\xC^1} \xrightarrow[j \to \infty]{} 0 
\end{align*} 
where we used Remark~\ref{remk:massdecrease4}, recalling \eqref{eqSpaceTimeFVcv} we obtain the desired distributional convergence
\begin{align*}
\left| \int X \cdot f_j \: d \mu_j - \delta \lambda(X) \right| 
& \leq \int_0^1 \left|  \Phi_{\e_j}\ast \delta V_{\e_j}(t) (X) - \delta V_{\e_j}(t) (X)  \right| \: dt + \left| \delta (V_{\e_j}(t) \otimes dt) (X) - \delta \lambda (X) \right|\\
& \xrightarrow[j \to \infty]{} 0 \: .
\end{align*}
Let $\psi \in \xC ([0,1] \times \R^n,\R_+)$ be bounded and consider $F : ((t,y), q) \mapsto \psi(t,y) |q|^2$, then $F$ is non-negative continuous, and with respect to $q$, it is convex and has superlinear growth, hence satisfying the assumptions of \cite{hut1} 
4.1.2. We additionally have by Remark~\ref{remk:massdecrease4} that for all $j$,
\begin{align*}
 \int_{[0,1] \times \R^n} F( (t,y) , f_j (t,y) ) \: d \mu_j (t,y) & = \int_{[0,1] \times \R^n} \psi (t,y) \frac{\left| \Phi_{\e_j}\ast \delta V_{\e_j}(t) \right|^2}{(\Phi_{\e_j} \ast \|V_{\e_j}(t)\| + \e_j)} \:  dy \: dt \\
 & \leq \| \psi \|_\infty \| V_{\e_j}(0) \|(\R^n) = \| \psi \|_\infty \| V_0 \|(\R^n) < \infty
\end{align*}
and we can apply
\cite{hut1} 
4.4.2(i) and (ii) (see also 2.36 in \cite{afp}): there exists $f \in L^1 ([0,1] \times \R^n, \R^n, \| \lambda \|)$ such that, up to extraction, the sequence of vector measures $f_j \mu_j$ converge to $f \| \lambda \|$ and
\begin{align}
\int_0^1 \int_{\R^n} \psi(t,y) |f(t,y)|^2 \: d \| \lambda \|   = & \int F((t,y), f(t,y)) \: d  \| \lambda \| \nonumber \\ 
\leq \liminf_{j \to \infty} & \int F( (t,y) , f_j (t,y) ) \: d \mu_j (t,y) \nonumber \\
\leq \liminf_{j \to \infty} & \int_{[0,1] \times \R^n} \psi (t,y) \frac{\left| \Phi_{\e_j}\ast \delta V_{\e_j}(t) \right|^2}{(\Phi_{\e_j} \ast \|V_{\e_j}(t)\| + \e_j)} \:  dy \: dt \leq  \| \psi \|_\infty \| V_0 \|(\R^n)  \label{eqSpaceTimeL2LSC01} \: .
\end{align}
Thanks to the uniqueness of the distributional limit: $f \| \lambda \| = \delta \lambda$ so that $(\delta \lambda)_s = 0$ and $f = - h( \cdot, \cdot, \lambda)$, and we obtain \eqref{eqSpaceTimeL2MC} plugging $\psi = 1$ in \eqref{eqSpaceTimeL2LSC01}.

\noindent We are left with proving \eqref{eqSpaceTimeL2LSC} for $0 \leq t_1 \leq t_2 \leq 1$, and we can take an affine cut--off approximating $\one_{[t_1, t_2]}$ from below: for $k$ large enough with respect to $t_2 - t_1$, let $\chi_k$ be a continuous piecewise-affine function satisfying $ \one_{[t_1 + \frac{1}{k}, t_2 - \frac{1}{k}]} \leq \chi_k \leq \one_{[t_1, t_2]}$, then applying  \eqref{eqSpaceTimeL2LSC01} to $\chi_k \psi$ gives
\begin{equation*}
\int_0^1 \int_{\R^n} \chi_k \psi |h|^2 \: d \| \lambda \| 
\leq \liminf_{j \to \infty} \int_{[0,1] \times \R^n} \chi_k \psi  \frac{\left| \Phi_{\e_j}\ast \delta V_{\e_j}(t) \right|^2}{(\Phi_{\e_j} \ast \|V_{\e_j}(t)\| + \e_j)} \:  dy \: dt
\end{equation*}
and we can take the limit $k \to \infty$ in the l.h.s. by dominated convergence while we use $\chi_k \leq \one_{[t_1, t_2]}$ in the r.h.s. to conclude the proof of \eqref{eqSpaceTimeL2LSC}, and hence the current proof.
\end{proof}
\begin{prop}[Spacetime Brakke inequality for the limit flow]\label{spacetimebrakkeinequality}
Let $\lambda = \dt \otimes V(t)$ be the limit measure defined in Proposition~\ref{limitflow}.
We assume that $V(t)$ is a rectifiable $d$--varifold for a.e $t\in[0,1]$: 
\[
V(t) = \theta_t \cH^d_{| \cM_t} \otimes \delta_{T_\cdot \cM_t}= \| V(t) \| \otimes  \delta_{T_\cdot \cM_t} \: .
\]
Then, for any $\phi \in \xC_c^1([0,1] \times \R^n, \R_+)$ and $0\leq t_1 \leq t_2 \leq 1$,
\begin{equation*}
\begin{split}
 \|V(t_2) \|(\phi(t_2,\cdot)) - & \|V(t_1)\|(\phi(t_1,\cdot))  \leq -\int_{t_1}^{t_2}\int_{\R^n} \phi(t,y)|h(t,y,\lambda)|^2 \: d \| V(t) \| \: dt
\\ +\int_{t_1}^{t_2} &\int_{\R^n \times \G } T_y\cM_t^{\perp}(\nabla\phi(t,y))\cdot h(t,y,\lambda) \: d \| V(t) \| \: dt + \int_{t_1}^{t_2} \int_{\R^n}\partial_t\phi(t,\cdot) \: d \| V(t) \| \: dt \: .
\end{split}
\end{equation*}
\end{prop}
\begin{proof}
Denote $\lambda_{\e_j}= \dt \otimes V_{\e_j}(t) $ and choose (as in Proposition~\ref{limitflow}) a sequence $\left( \e_j \right)_j$ satisfying:
\begin{equation*}
\lim\limits_{j \to \infty} \lambda_{\e_j} = \lambda = \dt \otimes V(t) \quad \text{and} \quad   \lim\limits_{j \to \infty}\|V_{\e_j}(t)\|= \| V(t) \| \: .
\end{equation*}
Consider $\phi \in \xC_c^1([0,1] \times \R^n,\R_+)$ , $t_1,t_2$ such that $0 \leq t_1 \leq t_2 \leq 1$. The inequality we are seeking to prove is linear in $\phi$, without loss of generality we assume $\phi < 1$,  and for all sufficiently large $i \in \N$ we define $\phi_i:= \phi + i^{-1} < 1$. We can plug $\phi_i$ in \eqref{timeepsilonbrakke}, also recalling \eqref{wfirstvar}, we obtain:
\begin{multline} \label{eqStep0BrakkeIneq}
\|V_{\e_j}(t_2)\|(\phi_i(t_2,\cdot)) - \|V_{\e_j}(t_1)\|(\phi_i(t_1,\cdot)) 
- \int_{t_1}^{t_2} \|V_{\e_j}(t)\|(\partial_t \phi_i(t,\cdot)) \: dt \\
= \int_{t_1}^{t_2} \int_{\R^n\times \G} S^{\perp}(\nabla_x \phi_i) \cdot h_{\e_j}(\cdot,V_{\e_j}(t)) \,d V_{\e_j}(t) \,dt + \int_{t_1}^{t_2} \delta (V_{\e_j}(t)) \Big[ \phi_i (t,\cdot) h_{\e}(\cdot,V_{\e}(t)) \Big] \,dt  
\end{multline}
and the proof now consists in taking the limit, first in $j$ and then in $i$.

\medskip
{\bf Step 1:} We take the limit in the l.h.s. of \eqref{eqStep0BrakkeIneq}, that is, we prove
\begin{multline} \label{eqStep1BrakkeIneq}
\|V_{\e_j}(t_2)\|(\phi_i(t_2,\cdot)) - \|V_{\e_j}(t_1)\|(\phi_i(t_1,\cdot)) 
- \int_{t_1}^{t_2} \|V_{\e_j}(t)\|(\partial_t \phi_i(t,\cdot)) \: dt \\
\xrightarrow[i,j \to \infty]{} \|V(t_2)\|(\phi(t_2,\cdot)) - \|V(t_1)\|(\phi(t_1,\cdot)) 
- \int_{t_1}^{t_2} \|V(t)\|(\partial_t \phi(t,\cdot)) \: dt.
\end{multline}

\noindent First note that $\partial_t \phi_i = \partial_t \phi$ and recall that for all $t \in [0,1]$, 
\[
 \| V_{\e_j} \|(t) \xrightharpoonup[j \to \infty]{\ast} \|V(t)\| \quad \Longrightarrow \quad \left\lbrace \begin{array}{l}
 \|V_{\e_j}(t)\|(\partial_t \phi_i(t,\cdot)) = \|V_{\e_j}(t)\|(\partial_t \phi (t,\cdot))  \xrightarrow[j \to \infty]{} \|V(t)\| (\partial_t \phi(t,\cdot)) \\
 \|V_{\e_j}(t)\|(\phi(t,\cdot))  \xrightarrow[j \to \infty]{} \|V(t)\|( \phi(t,\cdot))
 \end{array}
 \right.
\]
and since $\|V_{\e_j}(t)\|(\partial_t \phi (t,\cdot)) \leq \| \partial_t \phi \|_\infty \|V_0 \|(\R^n)$ by the decay of the mass (Remark \ref{remk:massdecrease4}), we infer by dominated convergence that for any $i$,
\[
 \int_{t_1}^{t_2} \|V_{\e_j}(t)\|(\partial_t \phi_i(t,\cdot)) \: dt  \xrightarrow[j \to \infty]{} \int_{t_1}^{t_2} \|V(t)\|(\partial_t \phi(t,\cdot)) \: dt \: .
\]
Using again the decay of the mass and $\phi_i = \phi + i^{-1}$, we obtain
\[
 \left| \|V_{\e_j}(t)\|(\phi_i(t,\cdot))  - \|V(t)\|( \phi(t,\cdot))  \right| \leq i^{-1} \|V_0\|(\R^n) + \left| \|V_{\e_j}(t)\|(\phi(t,\cdot))  - \|V(t)\|( \phi(t,\cdot))  \right| \xrightarrow[i,j \to \infty]{} 0
\]
and with $t = t_1, t_2$ we can conclude the proof of \eqref{eqStep1BrakkeIneq} (Step 1).
We now deal with the two terms involving the mean curvature. 

{\bf Step 2:}
We now prove that
\begin{equation}\label{brakkestep4}
 \limsup_{i \to \infty} \limsup_{j \to \infty} \int_{t_1}^{t_2} \delta V_{\e_j}(t)(\phi_ih_{\e_j}(\cdot,V_{\e_j}(t))) dt \leq - \int_{t_1}^{t_2} \int_{\R^n} \phi(t,y)|h(t,y,\lambda)|^2 \: d \|V(t)\| (y) \: dt.
\end{equation}
First, we note that $\phi_i = \phi + i^{-1}$ and then, there exists $m_{i,\phi} \in \N $ (large enough, depending on $i$ and $\phi$: e.g. $m_{i,\phi} \geq i \| \nabla_x \phi \|_\infty$) such that for all $m \geq m_{i,\phi}$, $\phi_i \in \mathcal{A}_m$.
We apply \eqref{eqKTprop54} with $\e = \e_j$ and $\psi = \phi_i$ whence, for fixed $\phi$ and $i$, one has to take $j$ large enough to ensure $\e_j \leq \e_{*}$ and $\e_j^{-\frac{1}{6}} \geq 2m_{i,\phi}$: this is the reason why we have to take $\lim_{j \to \infty}$ before $\lim_{i \to \infty}$ hereafter. Concerning the varifold, we apply \eqref{eqKTprop54} with $W = V_{\e_j}(t)$ (for $t \in [0,1]$) and $M = \|V_0\|(\R^n)$ since we know that $\|V_{\e_j}(t)\|(\R^n) \leq \|V_0\|(\R^n) \leq M$. 
We obtain, for all $t \in [0,1]$ and for all $j$ large enough,
\begin{equation*} 
\left| \delta V_{\e_j}(t) \left(\phi_i h_{\e_j}(\cdot, V_{\e_j}(t)) \right) + \int_{\R^n} \frac{\phi_i |\Phi_{\e_j}\ast\delta V_{\e_j}(t)|^2 }{\Phi_{\e_j}\ast \|V_{\e_j}(t)\|+\e_j} \: dx \right| \leq \e_j^\frac{1}{4} \left( 1 + \int_{\R^n} \frac{\phi_i |\Phi_{\e_j}\ast\delta V_{\e_j}(t)|^2 }{\Phi_{\e_j}\ast \|V_{\e_j}(t)\|+\e_j} \: dx \right)
\end{equation*}
which we integrate between $t_1$ and $t_2$ so that using $0 \leq \phi_i \leq 1$ and Remark~\ref{remk:massdecrease4},
\begin{multline*} 
\int_{t_1}^{t_2} \left| \delta V_{\e_j}(t) \left(\phi_i h_{\e_j}(\cdot, V_{\e_j}(t)) \right) + \int_{\R^n} \frac{\phi_i |\Phi_{\e_j}\ast\delta V_{\e_j}(t)|^2 }{\Phi_{\e_j}\ast \|V_{\e_j}(t)\|+\e_j} \: dx \right|  dt 
\\ \leq \e_j^\frac{1}{4} \left( 1 + \int_{[0,1] \times \R^n} \frac{\phi_i |\Phi_{\e_j}\ast\delta V_{\e_j}(t)|^2 }{\Phi_{\e_j}\ast \|V_{\e_j}(t)\|+\e_j} \: dx \: dt \right)
\leq \e_j^\frac{1}{4} \left( 1 + \|V_0\|(\R^n) \right) \xrightarrow[j \to \infty]{} 0 \: .
\end{multline*}
We infer:
\begin{equation*}\begin{split}
 \limsup_{j \to \infty} \int_{t_1}^{t_2} \delta V_{\e_j}(t)(\phi_ih_{\e_j}(\cdot,V_{\e_j}(t))) dt &= -\liminf_{j \to \infty} \int_{t_1}^{t_2} \int_{\R^n} \frac{\phi_i|\Phi_{\e_j}\ast\delta V_{\e_j}(t)|^2 }{\Phi_{\e_j}\ast \|V_{\e_j}(t)\|+\e_j} \: dx \: dt \\
 \leq - \int_{t_1}^{t_2} \int_{\R^n} & \phi_i|h(\cdot, \cdot, \lambda)|^2 d \|V(t)\| \: dt \quad \text{by \eqref{eqSpaceTimeL2LSC} in Proposition~\ref{propSpaceTimeL2MC}}
\end{split}\end{equation*}
and the proof of \eqref{brakkestep4} (Step 2) follows from $-\phi_i \leq -\phi$. 

{\bf Step 3:} As $\nabla_x \phi_i = \nabla_x \phi$, we are left with the proof of
\begin{equation}\label{brakkestep5}
\limsup_{j \to \infty} \int_{t_1}^{t_2}\int_{\R^n \times \G } S^{\perp}(\nabla_y \phi) \cdot h_{\e_j}(y,V_{\e_j}(t)) \: d V_{\e_j}(t) \: dt
 \leq  \int_{t_1}^{t_2} \int_{\R^n \times \G } S^{\perp}(\nabla_y \phi)\cdot h(t,y,\lambda) \, d\lambda \: .
\end{equation}

\noindent We fix an arbitrary function $g\in \xC_c^2(\R^n\times[0,1],\R^n)$, there exists $m_{g} \in\N$ such that $g(t,\cdot))\in \mathcal{B}_{m} \, \forall m \geq m_{g} $ and  $ \forall t\in[0,1]$, this is due to the compactness of $[0,1]$.
We apply \eqref{eqKTprop55} with $\e = \e_j$ and $X=g$ and we take $j$ large enough to ensure that $\e_j \leq \e_*$ and $\e_j^{-\frac16} \geq 2m_{g}$. Concerning the varifold, we apply \eqref{eqKTprop55} with $W= V_{\e_j}(t)$ (for $t\in [0,1]$) and $M=\| V_0\|(\R^n)$ since we know that $\|V_{\e_j}(t)\|(\R^n) \leq \|V_0\|(\R^n) \leq M$ by Remark \ref{remk:massdecrease4}. We obtain, for $j$ large enough 
\begin{equation*}
  \bigg| \int_{\R^n} h_{\e_j}(\cdot,V_{\e_j}(t))   \cdot g(t,\cdot) \, d \|V_{\e_j}(t)\| + \delta V_{\e_j}(t)(g(t,\cdot)) \bigg| \leq  \e_j^{\frac14} + \e_j^{\frac14} \left( \int_{\R^n} \frac{|\Phi_{\e} \ast \delta V_{\e_j}(t)|^2}{\Phi_{\e} \ast \|V_{\e_j}(t)\| + \e} \, dx\right)^{\frac12}
\end{equation*}
which we integrate between $t_1$ and $t_2$ so that
\begin{align}
\int_{t_1}^{t_2} \bigg| \int_{\R^n} h_{\e_j}(\cdot,V_{\e_j}(t)) &  \cdot g(t,\cdot) \, d \|V_{\e_j}(t)\| + \delta V_{\e_j}(t)(g(t,\cdot)) \bigg| \: dt \nonumber  \\
\leq & \int_{t_1}^{t_2}\e_j^{\frac14} \: dt + \e_j^{\frac14}  \int_{t_1}^{t_2} \left(\int_{\R^n} \frac{|\Phi_{\e_j}\ast\delta V_{\e_j}(t)|^2 }{\Phi_{\e_j}\ast \|V_{\e_j}(t)\|+\e_j}dx\right)^{\frac12} \: dt \nonumber \\
\leq & \e_j^{\frac14}  + \e_j^{\frac14} \left( \int_0^1 \int_{\R^n} \frac{|\Phi_{\e_j}\ast\delta V_{\e_j}(t)|^2 }{\Phi_{\e_j}\ast \|V_{\e_j}(t)\|+\e_j} \: dx \: dt \right)^{\frac12} \nonumber \\
\leq & \e_j^{\frac14} \left( 1 + \|V_0\|(\R^n)^{\frac12} \right) \xrightarrow[j \to \infty]{} 0 \: , \label{hgdeltag}
\end{align}
where we used Jensen inequality and Remark \ref{remk:massdecrease4}.

\noindent We now observe that the map $\displaystyle g:(t,y) \mapsto (T_{y}\cM_t)^{\perp}(\nabla\phi)$ is $\dt \otimes \-V(t)\|$--measurable and belongs to  $L^2(\dt \otimes \-V(t)\|)$ ($\|V(t)\|$ is finite and $\phi \in \xC_c^1$, hence $g$ is bounded by $\|\nabla \phi \|$), we can assert that, for any $\eta \in (0,1)$, there exists a map $g_{\eta}\in \xC_c^2(\R^n\times[t_1,t_2],\R^n)$ such that:
\begin{equation}\label{brakkestep5prf1}
 \int_{t_1}^{t_2} \int_{\R^n} \big| (T_{y}\cM_t)^{\perp}(\nabla\phi(y)) - g_{\eta}(t,y) \big|^2 d\|V(t)\|(y) \: dt < \eta^2.
\end{equation}
Now we compute
\begin{equation}\label{brakkestep5prf2}
 \begin{split}
  \int_{t_1}^{t_2} & \int_{\R^n\times \G} S^{\perp}(\nabla\phi)\cdot h_{\e_j}(\cdot,V_{\e_j}(t)) \, d\lambda_{\e_j}
   =  \int_{t_1}^{t_2} \int_{\R^n \times \G} \left(S^{\perp}(\nabla\phi) -g_{\eta} \right) \cdot h_{\e_j}(\cdot,V_{\e_j}(t)) \, d\lambda_{\e_j}
  \\& + \left( \int_{t_1}^{t_2} \int_{\R^n} g_{\eta}\cdot h_{\e_{j}}(\cdot,V_{\e_j}(t)) \, d\|V_{\e_j}(t)\|dt + \int_{t_1}^{t_2}\delta V_{\e_j}(t)(g_{\eta}) \, dt \right)
  \\&  -\int_{t_1}^{t_2} \delta V_{\e_j}(t)(g_{\eta}) \,dt + \delta \lambda(g_{\eta})
  \\& + \int_{t_1}^{t_2} \int_{\R^n} h(\cdot,\cdot,\lambda)\cdot \left( g_{\eta} - (T_{y}\cM_t)^{\perp}(\nabla\phi) \right) d \|V(t)\| \: dt
  \\& + \int_{t_1}^{t_2} \int_{\R^n} (T_{y}\cM_t)^{\perp}(\nabla\phi(t,y))  \cdot h(t,y,\lambda) \,d \|V(t)\| \: dt.
 \end{split}
\end{equation}
By the varifold convergence, \eqref{brakkestep5prf1} and \eqref{eq2KTprop54} we have 
\begin{equation}\label{brakkestep5prf3}
\begin{split}
 &\limsup_{j}
 \int_{t_1}^{t_2} \int_{\R^n \times \G} \left(S^{\perp}(\nabla\phi) -g_{\eta} \right) \cdot h_{\e_{j}}(\cdot,V_{\e_j}(t)) \,d\lambda_{\e_j} 
 \\& \leq \left(\int_{t_1}^{t_2} \int_{\R^n \times \G} |S^{\perp}(\nabla\phi) -g_{\eta} |^2 \, d\lambda\right)^\frac{1}{2}  \left( \limsup_j \int_{t_1}^{t_2}\int_{\R^n}  |h_{\e_{j}}(\cdot,V_{\e_j}(t)) |^2 \,d\|V_{\e_j}(t)\|\right)^\frac{1}{2}
 \\& = \left(\int_{t_1}^{t_2} \int_{\R^n} |(T_{y}\cM_t)^{\perp}(\nabla\phi) -g_{\eta} |^2 \, d \|V(t)\| \: dt\right)^\frac{1}{2}\left(\limsup_j \int_0^1 \int_{\R^n} \frac{|\Phi_{\e_j}\ast\delta V_{\e_j}(t)|^2 }{\Phi_{\e_j}\ast \|V_{\e_j}(t)\|+\e_j} \: dx \: dt \right)^{\frac12}
 \\& \leq  \eta \left( \|V_0\|(\R^n) \right)^{\frac12}
\xrightarrow[\eta \rightarrow 0]{}0.
 \end{split}
\end{equation}
By \eqref{hgdeltag} we have:
\begin{equation}\label{brakkestep5prf4}
 \limsup_{j}  \left( \int_{t_1}^{t_2} \int_{\R^n} g_{\eta} \cdot h_{\e_{j}}(\cdot,V_{\e_j}(t)) \, d\|V_{\e_j}\|(t) + \int_{t_1}^{t_2}\delta V_{\e_j}(t)(g_{\eta}) \, dt \right) =0.
\end{equation}
By the varifold convergence we have:
\begin{equation}\label{brakkestep5prf5}
  \limsup_{j} |\int_{t_1}^{t_2}\delta V_{\e_j}(t)(g_{\eta}) \, dt- \delta \lambda(g_{\eta}) | =0,
\end{equation}
and finally, the Cauchy-Schwarz inequality and \eqref{brakkestep5prf1} imply 
\begin{equation}\label{brakkestep5prf6}
\begin{split}
 & \int_{t_1}^{t_2} \int_{\R^n} h(\cdot,\cdot,\lambda)\cdot \left( g_{\eta} - (T_{y}\cM_t)^{\perp}(\nabla\phi) \right) d \|V(t)\| \: dt \\& \leq \left(\int_{t_1}^{t_2} \int_{\R^n} |h(\cdot,\cdot,\lambda)|^2 \, d \|V(t)\| \: dt \right)^{\frac{1}{2}} \left( \int_{t_1}^{t_2} \int_{\R^n } |g_{\eta} - (T_{y}\cM_t)^{\perp}(\nabla\phi)|^2 \, d \|V(t)\| \: dt \right)^{\frac{1}{2}}
 \\&
 \leq (\|V_0\|(\R^n))^{\frac12} \eta\xrightarrow[\eta \rightarrow 0]{}0.
 \end{split}
\end{equation}
From \eqref{brakkestep5prf2}-\eqref{brakkestep5prf6} we deduce \eqref{brakkestep5} (Step 3), this concludes the proof of Proposition \ref{spacetimebrakkeinequality}. 
\end{proof}
\begin{remk}
From the proof, writing $V(t) = \|V(t)\| \otimes \nu^{(t,x)}$, we notice that assuming that $\nu^{(t,x)}$
is a Dirac measure is sufficient, as we do not use the fact that $T_{\cdot} \cM_t$ is the tangent space nor the properties of $\|V(t)\|$ as a rectifiable measure.
\end{remk}

\begin{remk}[Non uniqueness of the limit spacetime Brakke flows]

\noindent We recall that the limit measure $\lambda$ in Theorem \ref{chap3_theo:cv} depends on the choice of the subsequence $\left(\e_j\right)_{j\in \N}$, hence is not unique (in general). This, somehow, is related to the non-uniqueness of Brakke flows, as Brakke flows themselves are spacetime Brakke flows when tensored with the measure $dt$.
\end{remk}
\appendix

\section{Some useful lemmas in matrix algebra}
\label{appendix}

We collect in this appendix several results in linear algebra used in the proofs of Section~\ref{def_time_discrete}.
The following result follows directly from the triangle inequality, but since it is used multiple times in the proof of Proposition~\ref{damcfstability}, we think it is useful to have it as a lemma.
\begin{lemma} \label{lem:intro_trivial_lemma}
Let $A,B\in \cM_n$ and $S,T \in \cM_{d,n}$  be such that  $\|S\|= \|T\|=1$, one has
\begin{equation}\label{intro_trivial_lem}
 \| S A S^t - T B T^t \| \leq \left( \|A\| +\|B\| \right) \|S-T\| + \|A-B\|.
\end{equation}
\end{lemma}
\begin{proof}
Using the triangle inequality, and the fact that $\|M^t\|=\|M\|$ for any matrix $M$ we infer that 
\begin{equation*}
\begin{split}
\| S A S^t - T B T^t \| & \leq \|S A S^t-S A T^t \| + \|S A T^t- S B T^t \| + \|S B T^t - T B T^t \|
\\& \leq \|S\|  \|A \| \| S^t- T^t \| + \| S\| \| A-B\| \|T^t\| + \|S-T\| \|B\| \|T^t\|
\\& \leq \left( \|A\| +\|B\| \right) \|S-T\| + \|A-B\| \: .
\end{split}
\end{equation*}
This concludes the proof of \eqref{intro_trivial_lem}.
\end{proof}

The following lemma contains several properties on  determinant's expansions, mainly used to prove Propositions \ref{mcfmotion1} and \ref{damcfstability}. The matrix product is denoted by $\circ$ for clarity.

%
\begin{lemma}\label{detprop}
There exists a constant $c_2 \geq 1$ only depending on $n$ such that the following estimates hold:
\begin{enumerate}
 \item Let $1 \leq k \leq n$ and $Q \in \mathcal{M}_k$ be such that $| Q |_\infty \leq 1$, then
\begin{equation} \label{eq:detprop1}
 \left| \det(I_k + Q) - \det(I_k) \right| \leq c_2 | Q |_\infty \quad \text{and} \quad \left| \det(I_k+Q) - \det(I_k) - \tr(Q) \right| \leq  c_2 (| Q |_\infty)^2 \: .
\end{equation}

\item Let $1\leq d \leq n$, $L\in \mathcal{M}_{d,n}$ , $R \in \mathcal{M}_n$ be such that $ L \circ L^t = I_d$ and $| R |_\infty \leq 1$. Then
\begin{equation} \label{detprop3}
((I_n + R)\circ L^t)^t\circ ((I_n +R) \circ L^t) = I_d + Q \quad \text{with} \quad \left|  Q \right|_\infty \leq c_2   | R |_\infty \: .
\end{equation}
Moreover, if we assume that $c_2  |R|_\infty \leq 1$, then
\begin{equation}\label{detprop4}
 \Big| \det\left( (( I_n +R)\circ L^t)^t\circ ((I_n +R) \circ L^t)  \right)^{\frac12} -1  \Big| \leq c_2 |R|_\infty,
\end{equation}
and 
\begin{equation}\label{detprop5}
 \Big| \det\left( (( I_n +R)\circ L^t)^t\circ ((I_n +R) \circ L^t)  )\right)^{\frac12} -1 - \tr( R \circ L^t \circ L)  \Big| \leq c_2 |R|_\infty^2.
\end{equation}
\end{enumerate}
\end{lemma}
\begin{proof}
As long as $c_2$ only depends on $n$, we may increase it whenever needed throughout the proof.
We consider the normed space $(\cM_{p,q}, | \cdot |_\infty)$ and we recall that for $M \in \cM_{p,q}$ and $N \in \cM_{q,r}$, 
\begin{equation} \label{eq:multInfMatrixNorm}
| M N |_\infty \leq q \: |M|_\infty | N |_\infty \: .
\end{equation}
Let $Q \in \cM_k$ be such that $|Q |_\infty \leq 1$ and let $B = \{ M \in \cM_k \: : \: | I_k - M |_\infty \leq 1 \}$ be the closed unit ball centered at $I_k$. By compactness of $B$, we can introduce
\begin{equation*}
 c_{2,k} = \max \left\lbrace 1 ,  \max\limits_{M \in B}\|D\det (M)\|_\infty , \frac{1}{2} \max\limits_{M \in B} \|  D^2\det (M)\|_\infty \right\rbrace \:,
\end{equation*}
where $\| \cdot \|_\infty$ and $|\!\| \cdot |\!\|_\infty$ denote, respectively, the linear and bilinear operator norms associated with $(\cM_{p,q}, | \cdot |_\infty)$. Note that $c_{2,k}$ depends on $k$ (since $B$ depends on $k$) though this can be avoided by defining
$c_2 = \max_{1\leq k \leq n} c_{2,k}$. Recall that the differential of the determinant map $\det$ at $I_k$ is the trace map, i.e. $D \det (I_k) = \tr$, therefore,
the application of the Taylor-Lagrange inequality for $\det$ on the line segment $[I_k, I_k +Q] \subset B$ yields \eqref{eq:detprop1}.

Let $L\in \mathcal{M}_{d,n}$ , $R \in \mathcal{M}_n$ be such that $ L \circ L^t = I_d$ and $| R |_\infty \leq 1$ and let us use the notation $Q = L \circ \left( R^t + R \right)\circ L^t + L\circ R^t \circ R \circ L^t \in \cM_d$ hereafter, so that $((I_n + R)\circ L^t)^t\circ ((I_n +R) \circ L^t) = I_d + Q$.
First note that $| L |_\infty \leq 1$, indeed, the assumption $L \circ L^t = I_d$ can be reformulated as follows: the columns of $L^t$ (i.e. the rows of $L$) constitute an orthonormal family $(v_1, \ldots, v_d)$ of $\R^n$ so that $| L |_\infty = \max_{ij} |L_{ij}| = \max_{ij} |v_i \cdot e_j | \leq 1$.
Using \eqref{eq:multInfMatrixNorm} and $| M^t |_\infty = | M |_\infty$,  we have
\begin{equation*}
\begin{split}
| Q |_\infty & = \left| L\circ \left( R^t + R \right)\circ L^t +L\circ R^t \circ R \circ L^t \right|_\infty  \leq n^2 | R + R^t |_\infty | L |_\infty^2 + n^3 | L |_\infty^2 | R |_\infty^2 \\
& \leq ( 2 n^2 + n^3 | R |_\infty ) | R |_\infty 
\leq c_2 | R |_\infty \: ,
 \end{split}
\end{equation*}
that is \eqref{detprop3}. 

We now
assume $c_2 | R |_\infty \leq 1$ and consequently $(( I_n +R) \circ L^t)^t\circ ((I_n +R) \circ L^t)  = I_d + Q$ with $| Q |_\infty \leq 1$
so that the first part of \eqref{eq:detprop1} gives
\begin{equation} \label{eq:detpropz}
 \left| \det\left( I_d + Q  \right) -1  \right|
  \leq c_2 | Q |_\infty \leq c_2^2 | R |_\infty \leq 1 \quad \text{and in particular} \quad  \det\left( I_d + Q  \right) \geq 0 \: .
\end{equation}
We infer \eqref{detprop4} applying $|a-1| \leq |a^2-1|$ with $a = \det\left( I_d + Q  \right)  \geq 0$. We are left with the proof of \eqref{detprop5}. We now apply the second inequality in \eqref{eq:detprop1} to obtain
\begin{equation}\label{detproof2}
\left| \det\left( I_d + Q  \right) -1 - \tr Q \right|
 \leq c_2 |Q|_\infty^2 \leq c_2^3 |R|_\infty^2 \: .
\end{equation}
Furthermore, using $\tr(A) = \tr(A^t)$ and 
$\tr(AB) = \tr(BA)$ when both products make sense, we have
\[
\tr( L \circ R^t \circ L^t) = \tr(L \circ R \circ L^t) =\tr (R \circ L^t \circ L) 
\]
and thus, by definition of $Q$ and \eqref{eq:multInfMatrixNorm},
\begin{align} \label{eq:detproof3}
\left| \tr Q - 2 \tr (R \circ L^t \circ L) \right| & = 
  \left| \tr( L\circ R^t \circ R \circ L^t) \right| \leq d | L\circ R^t \circ R \circ L^t |_\infty \leq d n^3 | L |_\infty^2 | R |_\infty^2 \nonumber \\
  & \leq n^4 | R |_\infty^2 \: .
\end{align}
From \eqref{detproof2} and \eqref{eq:detproof3} we obtain 
\begin{equation} \label{eq:detproof4}
\left| \det \left( I_d + Q \right)
  -1 - 2\tr(R \circ L^t \circ L)
 \right|\leq (c_2^3 + n^4) |R|_\infty^2
\end{equation}
We now apply the following inequality, valid for any $z \geq -1$,
\[
 \left| \sqrt{1+z} - 1 - \frac 12 z \right| \leq \frac12 z^2
\]
with $\displaystyle z = \det \left( I_d + Q \right) - 1$,  from \eqref{eq:detpropz} we know that $-1 \leq z \leq c_2^2 | R |_\infty$, we hence obtain
\begin{align*}
\left| \sqrt{\det(I_d + Q)} - 1 -  \tr(R \circ L^t \circ L) \right| & \leq \left| \sqrt{1 + z} - 1 - \frac 12 z \right| + \left| \frac{1}{2} z - \tr(R \circ L^t \circ L) \right| \\
& \leq \frac12 z^2 + \frac{1}{2} (c_2^3 + n^4) |R|_\infty^2 \text{ thanks to \eqref{eq:detproof4},} \\
& \leq \frac{1}{2} (c_2^4 + c_2^3 + n^4) |R|_\infty^2 
\end{align*}
which concludes the proof of \eqref{detprop5}.
\end{proof}
The following lemma stems directly from Lemma \ref{detprop}.
\begin{lemma}
Let $P, N \in \mathcal{M}_d$ and assume that $P$ is invertible, then 
\begin{equation}\label{claimcondition}
\| P^{-1} \| \: \| P-N \| \leq 1 \quad \Rightarrow \quad \left| \det(P) - \det(N) \right| \leq c_2 | \det(P) | \| P^{-1} \| \: \| P - N \| \:.
\end{equation}
\end{lemma}
\begin{proof}
Indeed, first note that
\begin{equation*}
\left| \det(P) - \det(N) \right| = | \det(P) | \: \big| 1 - \det( P^{-1}N) \big| \quad \text{and} \quad  P^{-1}N = I_d + P^{-1}(N-P) \: .
\end{equation*}
Furthermore
$\displaystyle
\left| P^{-1}(N-P) \right|_\infty \leq  \| P^{-1}(N-P) \| \leq  \| P^{-1} \| \: \| P - N \| \leq 1
$
so that applying \eqref{eq:detprop1} with $k=d$ and $Q = P^{-1}(N-P)$ we can assert that
\begin{equation*}
\big| 1 - \det( P^{-1}N) \big| \leq c_2 \left| P^{-1}(N-P) \right|_\infty \leq c_2 \: \| P^{-1} \| \: \| N-P \|  
\end{equation*}
which concludes the proof of \eqref{claimcondition}.
\end{proof}
The following is a crucial step to prove Proposition~\ref{damcfstability}. 
%
\begin{lemma}\label{StildetoS}
Let $S,T\in \G$, there exist $\tilde{S}=\left( s_1 | \dots | s_d \right)^t, \tilde{T}=\left( t_1 | \dots | t_d \right)^t \in \mathcal{M}_{d,n}$  where $\{s_i\}_{i=1}^d$ and $\{t_i\}_{i=1}^d$ are two orthonormal bases of $S$ and $T$ such that 
\begin{equation*}
 \| \tilde{S} - \tilde{T} \| \leq 2 \| S - T \|.
\end{equation*} 
\end{lemma}
\begin{proof}
Let $\theta$ be the largest principal angle between the subspaces $S$ and $T$, which can be characterized by:  
$$  \sin(\theta) = \max_s \min_t \sqrt{ 1 - \langle s,t \rangle^2 }, \,\, s \in S, t \in T \quad \text{and} \,|s|=|t|=1.$$ 
We infer from \cite[Proposition III.29]{amari} that 
$\| S - T \| = \sin(\theta)$, furthermore, there exists a rotation $r$ of $\R^n$ such that $r(S)=T$, with 
\begin{equation*}
 \| r - I_n \| \leq 2 \sin(\theta/2).
\end{equation*}
Let $\tilde{S}=\left( s_1 | \dots | s_d \right)^t\in \mathcal{M}_{d,n}$, with $\{s_i\}_{i=1}^d$ being an orthonormal basis of $S$, the matrix $\tilde{T}=r \circ \tilde{S}\in \mathcal{M}_{d,n}$ can be written as $\left( t_1 | \dots | t_d \right)^t $  where $\{t_i\}_{i=1}^d$ is an orthonormal basis of and $T$. We have then, using that $\| \tilde{S}\| =1$ 
\begin{equation*}
 \| \tilde{S} - \tilde{T} \|= \| \tilde{S}- r \circ \tilde{S} \|  \leq \| I_n - r \|\|\tilde{S}\| \leq 2 \sin(\theta/2);
\end{equation*}
the result follows from noting that $ 2 \sin(\theta/2) \leq 2\sin(\theta)$ as $\theta \in [0, \pi/2]$. 
\end{proof}

\bibliographystyle{abbrv}

 \end{document}